\documentclass[a4paper,fleqn,11pt]{article}

\usepackage[spanish, english]{babel}
\usepackage[spanish]{babel}
\usepackage[latin1]{inputenc}
\usepackage{epic}
\usepackage{amssymb}
\usepackage{amsmath}

\usepackage{calc}
\usepackage[authoryear]{natbib}
\usepackage{graphicx}
\usepackage{amsthm}
\usepackage{verbatim}
\usepackage{multirow}
\usepackage{tikz}
\usepackage{fancyhdr}
\usepackage{datetime}
\usepackage{pgffor}
\usepackage{dsfont}
\usepackage{enumerate}
\usepackage{setspace}
\usepackage{pdfpages}
\usepackage[normalem]{ulem}
\usepackage{rotating}

\usepackage[ruled,vlined]{algorithm2e}
\newcommand{\nota}[1]{}

\newcommand{\com}[1]{}

\newcommand{\formula}[1]{\begin{align}#1 & \notag\end{align}}
\newcommand{\formulaN}[1]{\begin{align}#1 &         \end{align}}

\newcommand{\hhp}{\hspace{-0.6cm}}

\newcommand{\Hp}{\hspace{-0.6cm}\hspace{1.23cm}}

\newcommand{\esp}{@{\hspace{0.2cm}}}

\newcommand{\TP}[1]{\textcolor[rgb]{0.00,0.00,1.00}{\textsc{#1}}\newline} %Títulos de parrafos (para saber en en una lectura rápida de que va la cosa)

\theoremstyle{plain}
\newcounter{LemmaCounter}
\newtheorem{lemma}[LemmaCounter]{Lemma}

\newtheorem{proposition}[LemmaCounter]{Proposition}
\newtheorem{property}[LemmaCounter]{Property}

\theoremstyle{definition}

\newcounter{examplecounter}

\newtheorem{example}[examplecounter]{Example}

\theoremstyle{remark}
\newcounter{RemarkCounter}

\newtheorem{remark}[RemarkCounter]{Remark}

\usepackage{anysize}
\usepackage[colorlinks=true,citecolor=black,linkcolor=black,urlcolor=black,colorlinks=true]{hyperref}

\usepackage{etex}
\reserveinserts{18}
\usepackage{morefloats}

\marginsize{1.5cm}{1.5cm}{1.5cm}{1.5cm} %(izq, der, sup, inf)
\renewcommand{\TP}[1]{}                         %Activar para quitar los títulos auxiliares de párrafo

\begin{document}
%=======================================================================================================================
\title{A modeling framework for\\ Ordered Weighted Average Combinatorial Optimization}

\author{Elena Fernández \footnote{Department of Statistics and Operational Research. BarcelonaTech. E-mail: e.fernandez@upc.edu} \\
Miguel A. Pozo \footnote{Department of Statistics and Operational Research. University of Seville. E-mail: miguelpozo@us.es} \\
Justo Puerto \footnote{Department of Statistics and Operational Research. University of Seville. E-mail: puerto@us.es} \\
}
\date{\today}

\maketitle
%=======================================================================================================================
\begin{abstract}
Multiobjective combinatorial optimization deals with problems considering more than one viewpoint or scenario. The problem of aggregating multiple criteria to obtain a globalizing objective function is of special interest when the number of Pareto solutions becomes considerably large or when a single, meaningful  solution is required. Ordered Weighted Average or Ordered Median operators are very useful when preferential information is available and objectives are comparable since they assign importance weights not to specific objectives but to their sorted values. In this paper, Ordered Weighted Average optimization problems are studied from a modeling point of view. Alternative integer programming formulations for such problems are presented and their respective domains studied and compared. In addition, their associated polyhedra are studied and some families of facets and new families of  valid inequalities presented. The proposed formulations are particularized for two well-known combinatorial optimization problems, namely, shortest path and minimum cost perfect matching, and the results of computational experiments presented and analyzed.
These results indicate that the new formulations reinforced with appropriate constraints can be effective for efficiently solving medium to large size instances.

\end{abstract}
%=======================================================================================================================

%=======================================================================================================================
\section{Introduction}

%\TP{Multiobjective optimization and aggregation functions}
Multiobjective combinatorial optimization deals with problems considering more than one viewpoint or scenario. They inherit the complexity difficulty of their scalar counterparts, but incorporate additional difficulties derived from dealing with partial orders in the objective function space. The standard solution concept is the set of Pareto solutions. However, the number of Pareto solutions can grow exponentially with the size of the instance and the number of objectives. A first approach to overcome this difficulty focuses on a specific subset of the Pareto set, such as, for instance, the supported Pareto solutions (see, for instance, \citealp{Ehrgott2005}). Those are the Pareto solutions that optimize linear scalarizations of the different objectives. However, it is possible to exhibit instances for which even the number of supported solutions grows exponentially with the size of the instance. Furthermore, focusing on supported Pareto solutions a priori excludes compromise solutions that could be preferred by the decision maker. For the above reasons, more involved decision criteria have been proposed in the field of multicriteria decision making \citep{Perny2003}. These include objectives focusing on one particular compromise solution, which, for tractability and decision theoretic reasons, seem to be better suited when an appropriate aggregation operator is available.

%\TP{OWA}
In some cases, a particularly important Pareto solution related to a weighted ordered average aggregating function is sought. Provided that some imprecise preference information on the objectives is available, and that they are comparable, an averaging operator can be used to aggregate the vector of objective values of feasible solutions. The Ordered Median (OM) objective function is very useful in this context since it assigns importance weights not to specific objectives but to their sorted values. OM operators have been successfully used for addressing various types of combinatorial problems (see, for instance, \citealp{Ogryczak2003}; \citealp{Nickel2005};  \citealp{Puerto2005}; \citealp{Boland2006}; or, \citealp{Fernandez2012}). \\

When applied to values of different objective functions in multiobjective problems, the OM operator is called in the literature Ordered Weighted Average (OWA) (\citealp{Yager1988,Yager1997}). It  assigns importance weights to the sorted values of the objective function elements in a multiple objective optimization problem. The OWA has been also used in the literature under the name of Choquet optimization to address continuous problems (\citealp{Schmeidler1986}) and more recently it has been applied to some combinatorial optimization problems, like the minimum spanning tree and 0-1 knapsack (\citealp{Galand2012}).
The OWA is, however, a very broad operator, which, depending on the cases, can be seen as an Ordered Median or as Vector Assignment Ordered Median \citep{Lei2012}, and which can be applied to any combinatorial optimization problem. We therefore believe that its full potential within combinatorial optimization is worth being exploited. This naturally leads to a thorough study of its modeling properties and alternatives, which is the focus of this paper.\\

From a modeling point of view, the OWA operator can be formulated with a combination of discrete and continuous decision variables linked by several families of linear constraints. Since the domain of combinatorial optimization problems can be characterized with ad hoc discrete variables and linear constraints, it becomes clear that any combinatorial optimization problem with an OWA objective can be formulated as a linear integer programming problem, by suitably relating the two sets of variables and constraints. Of course, not all formulations are equally useful. Moreover, it is not even clear that the best formulation for the domain of the combinatorial object should be preferred, because its ``integration'' with the formulation of the OWA may imply additional difficulties. In this work we propose three alternative basic formulations for a combinatorial object with an OWA objective. Each basic formulation uses a different set of decision variables to model the OWA objective. We study properties yielding to alternative formulations, which preserve the set of optimal solutions, and we also compare the formulations among them. In addition we propose various families of facets and valid inequalities, which can be used (independently or in combination) to reinforce the basic formulations.
In the final part of the paper, we focus on two classical optimization problems: shortest path and minimum cost perfect matching. For these two problems we analyze the empirical performance of the alternative basic formulations and their possible reinforcements and variations.  From our computational experience we can not conclude that any of the formulations is superior to the others since the behavior of the proposed formulation varies with the different combinatorial object to be considered (see Section  \ref{Sec:Computational}).\\

The paper is structured as follows. Section \ref{Sec:Definition} gives the formal definition of the OWA operator and shows that it has as particular cases both the Ordered Median and the Weighted Assignment Ordered Median. Section \ref{Sec:ModelingtheOWAP} presents the three basic formulations, and their variations, for a combinatorial problem with an OWA objective, studies their properties and compares them, whereas Section \ref{Sec:ValIneq} presents diferent families of valid inequalities and possible reinforcements. Sections \ref{Sec:ShortestPath} and \ref{Sec:PerfectMatching} respectively present the formulation of the combinatorial object that we use in our empirical study of the shortest path and minimum cost perfect matching problems with an OWA objective. Finally, Section \ref{Sec:Computational} describes the computational experiments that we have run and presents and analyzes the obtained numerical results. The paper ends in Section \ref{Sec:Conclusions} with some comments and possible avenues for future research.
%=======================================================================================================================

%=======================================================================================================================
\section{The Ordered Weighted Average Optimization}\label{Sec:Definition}
%\TP{OWA: Definition}

The Ordered Weighted Average (OWA) operator is defined over a feasible set $Q \subseteq \mathbb{R}^n$. Let $C\in \mathbb{R}^{p\times n}$ be a given matrix, whose rows, denoted by $C^i$, are associated with the cost vectors of $p$ objective functions. The index set for the rows of $C$ is denoted by $P=\{1, \dots, p\}$. For $x\in Q$, the vector $y\in \mathbb{R}^p$ with $y=Cx$ is referred to as the outcome vector relative to $C$.
For a given $y=Cx$, with $x\in Q$, let $\sigma$ be a permutation of the indices of $i\in P$ such that $y_{\sigma_1}\geq\hdots\geq y_{\sigma_p}$. Let also $\mathbf{\omega}\in\mathbb{R}^{p+}$ denote a vector of non-negative weights.
Feasible solutions $x\in Q$ are evaluated with an operator defined as $OWA_{(Q, C, \omega)}(x)=\omega ' y_\sigma$. The OWA optimization Problem (OWAP) is to find $x\in Q$ of minimum value with respect to the above operator.\\

%---------------------------------------
\begin{example}\label{ex1}
Consider
$$Q=\left\{x\in\{0,1\}^3:x_1+x_2+x_3=2\right\} \, , \,
C=\left(\begin{array}{ccc}
      1 & 4&  1 \\ \hline
      1 & 1 & 3 \\ \hline
      5 & 1 & 2 \\
    \end{array}\right) \textrm{ and }
\omega ' =\left(\begin{array}{ccc}1 & 2 & 4\end{array}\right).$$

Table \ref{Tabla1} illustrates, for each feasible $x\in Q$, the values of $y=Cx$,  $y_{\sigma}$ and $OWA_{(Q, C, \omega)}(x)=\omega ' y_\sigma$. The optimal value to the OWAP is $ \min_{x\in Q}OWA_{(Q, C, \omega)}(x)=23$.

 \begin{table}[h!]
\renewcommand{\esp}{@{\hspace{0.15cm}}}
\begin{center}
\begin{tabular}{|@{}\esp c \esp |c\esp |c\esp |c\esp  @{}|}
\hline
$x$ & $y$ &  $y_\sigma$ & $OWA_{(Q, C, \omega)}(x)=\omega ' y_\sigma$\\
\hline											
$\left(\begin{array}{ccc}1 & 1 & 0\end{array}\right)'$	& $\left(\begin{array}{ccc}5 & 2 & 6\end{array}\right)'$		 & $\left(\begin{array}{ccc}6 & 5 & 2\end{array}\right)'$	 & 24\\
\hline
$\left(\begin{array}{ccc}1 & 0 & 1\end{array}\right)'$	& $\left(\begin{array}{ccc}2 & 4 & 7\end{array}\right)'$		 & $\left(\begin{array}{ccc}7 & 4 & 2\end{array}\right)'$	 & 23\\
\hline
$\left(\begin{array}{ccc}0 & 1 & 1\end{array}\right)'$ & $\left(\begin{array}{ccc}5 & 4 & 3\end{array}\right)'$		 & $\left(\begin{array}{ccc}5 & 4 & 3\end{array}\right)'$	& 25\\
\hline
\end{tabular}
\caption[]{Solutions $x\in Q$, values $y=Cx$,  sorted values $y_{\sigma}$ and $OWA_{(Q, C, \omega)}(x)$ for Example \ref{ex1}.}
\label{Tabla1}
\end{center}
\end{table}		

\end{example}
%-----------------------------------------------------------------------------------------------------------------------
%\newpage
%-----------------------------------------------------------------------------------------------------------------------
The OWA operator is a very general function which, as we see below, has as particular cases well-known objective functions. We next describe some of them.

%\begin{itemize}
%\item \textbf{The ordered median objective function (OM).}
\subsection{The ordered median objective function (OM).}
The OM objective \citep{Nickel2005} minimizes a weighted sum of ordered elements. It is a well known function that unifies many location problems as the $p$-median problem, the $p$-center problem, etc.

Let $Q\subseteq \mathbb{R}^n$ denote the feasible domain for an optimization problem and let $d\in \mathbb{R}^n$ be a cost vector and $\omega\in \mathbb{R}^n$ a given weights vector. For $x\in Q$, let $\sigma$ denote a permutation of the indices of $x$, such that $d_{\sigma_j}x_{\sigma_j}\ge d_{\sigma_{j+1}}x_{\sigma_{j+1}}$, $j\in\{1, 2, \dots, n-1\}$.
The OM operator is $OM_{(Q,d,\omega)}(x)= \sum\limits_{j\in P}\omega_jd_{\sigma_j}x_{\sigma_j}$.\\

To cast the OM operator as an OWA operator, we only need to set the rows of the $C$ matrix as \((C^i)'=d_i \mathbf{e^i}\), $i\in\{1, \dots, n\}$, where \(\mathbf{e^i}\in \mathbb{R}^n\) is the $i$-th unit vector of the canonical basis of $\mathbb{R}^n$. Let $Diag(d)$ denote   the diagonal matrix whose diagonal entries are the components of the vector $d$, thus, $C=Diag(d)$. Then $OM_{(Q,d,\omega)}(x)=OWA_{(Q, Diag(d), \omega)}(x)$.

%---------------------------------------
\begin{example}\label{ex2}
Consider
$$Q=\left\{x\in\{0,1\}^3:x_1+x_2+x_3=2\right\} \, , \,
d=\left(\begin{array}{ccc} 5 & 1 & 2 \\ \end{array}\right)' \textrm{ and }
\omega =\left(\begin{array}{ccc}1 & 2 & 4\end{array}\right)'$$

Table \ref{Tabla2} illustrates, for each feasible $x\in Q$, the values of $(d_jx_j)_{j\in P}$, $(d_{\sigma_j}x_{\sigma_j})_{j\in P}$, and $OM_{(Q,d,\omega)}(x)= \sum_{j\in P}\omega_jd_{\sigma_j}x_{\sigma_j}$. The optimal OM value is $\min_{x\in Q}OM_{(Q, d, \omega)}(x)=4$.

\begin{table}[h!]
\renewcommand{\esp}{@{\hspace{0.15cm}}}
\begin{center}
\begin{tabular}{|@{}\esp c \esp |c\esp |c\esp |c\esp  @{}|}
\hline
$x$ & $(d_jx_j)_{j\in P}$ &  $(d_{\sigma_j}x_{\sigma_j})_{j\in P}$ & $OM_{(Q,d,\omega)}(x)= \sum_{j\in P}\omega_jd_{\sigma_j}x_{\sigma_j}$\\
\hline											
$\left(\begin{array}{ccc}1 & 1 & 0\end{array}\right)'$& $\left(\begin{array}{ccc}5 & 1 & 0\end{array}\right)'$	& $\left(\begin{array}{ccc}5 & 1 & 0\end{array}\right)'$	& 7\\
\hline
$\left(\begin{array}{ccc}1 & 0 & 1\end{array}\right)'$& $\left(\begin{array}{ccc}5 & 0 & 2\end{array}\right)'$&  $\left(\begin{array}{ccc}5 & 2 & 0\end{array}\right)'$	& 9\\
\hline
$\left(\begin{array}{ccc}0 & 1 & 1\end{array}\right)'$& $\left(\begin{array}{ccc}0 & 1 & 2\end{array}\right)'$& $\left(\begin{array}{ccc}2 & 1 & 0\end{array}\right)'$	& 4\\\hline
\end{tabular}
\caption[]{Solutions $x\in Q$, values $d_jx_j$,  sorted $d_{\sigma_j}x_{\sigma_j}$ and $OM_{(Q,d,\omega)}(x)$ for the OM of Example \ref{ex2}.}\label{Tabla2}
\end{center}
\end{table}

To cast the OM operator as an OWA operator, we only need to set the rows of the $C$ matrix as
$$C=Diag(d)=\left(\begin{array}{ccc}
      5 & 0 & 0 \\ \hline
      0 & 1 & 0 \\ \hline
      0 & 0 & 2 \\
    \end{array}\right).$$

The values of $y=Cx$, $y_{\sigma}$ and $OWA_{(Q, C, \omega)}(x)=\omega ' y_\sigma$ are shown in Table \ref{Tabla3}. The optimal OWA value is $\min_{x\in Q}OWA_{(Q, Diag(d), \omega)}(x)=\min_{x\in Q}OM_{(Q, d, \omega)}(x)=4$.

 \begin{table}[h!]
\renewcommand{\esp}{@{\hspace{0.15cm}}}
\begin{center}
\begin{tabular}{|@{}\esp c \esp |c\esp |c\esp |c\esp  @{}|}
\hline
$x$ & $y$ &  $y_\sigma$ & $OWA_{(Q, C, \omega)}(x)=\omega ' y_\sigma$\\
\hline											
$\left(\begin{array}{ccc}1 & 1 & 0\end{array}\right)'$& $\left(\begin{array}{ccc}5 & 1 & 0\end{array}\right)'$	& $\left(\begin{array}{ccc}5 & 1 & 0\end{array}\right)'$	& 7\\
\hline
$\left(\begin{array}{ccc}1 & 0 & 1\end{array}\right)'$& $\left(\begin{array}{ccc}5 & 0 & 2\end{array}\right)'$	& $\left(\begin{array}{ccc}5 & 2 & 0\end{array}\right)'$	& 9\\
\hline
$\left(\begin{array}{ccc}0 & 1 & 1\end{array}\right)'$& $\left(\begin{array}{ccc}0 & 1 & 2\end{array}\right)'$	& $\left(\begin{array}{ccc}2 & 1 & 0\end{array}\right)'$	& 4\\
\hline
\end{tabular}
\caption[]{The OM instance of Example \ref{ex2} as an OWA: $y=Cx$, $y_{\sigma}$ and $OWA_{(Q, C, \omega)}(x)$.}\label{Tabla3}
\end{center}
\end{table}		

\end{example}
%-----------------------------------------------------------------------------------------------------------------------
%\newpage
%-----------------------------------------------------------------------------------------------------------------------
%\item \textbf{The vector assignment ordered median objective function.}
\subsection{The vector assignment ordered median objective function.}
The Vector Assignment Ordered Median (VAOM) problem was recently introduced by \citet{Lei2012} in the context of discrete location-allocation problems. In this context, the VAOM generalizes both OM and Vector Assignment Median \citep{Weaver1985}. As we see below the OWA generalizes the VAOM as well. First, we briefly introduce the VAOM.\\

The main decisions in location-allocation problems are the set of facilities to open, and the assignment of customers to open facilities so as to satisfy their demand. Consider a given set of customers $P=\{1, \dots, p\}$, where each customer is also a potential location for a facility, and let $q\le p$ denote the number of facilities to open. Associated with each customer $i\in P$ there is a demand $a_i$. A unit of demand at customer $i$ served from facility $k$  incurs a cost $d_{ik}$. We will use \(\mathbf{d}^i\) to denote the $p$ dimensional vector of the distances associated with customer $i$. Usual objectives focus on service cost minimization.

Many location-allocation models allow splitting the demand at customers among several facilities, so allocating customer $i$ to facility $k$ means that some positive fraction of $a_i$ is served from facility $k$. However, without any further incentive or constraint, in optimal solutions customers will be allocated to one single facility, the closest one among those that are open. Since such solutions often exhibit privileged customers, equity measures have been proposed to balance out the service level of the customers.
This is the case of the VAOM that imposes  the specific fractions of the demand at each customer to be served from the various open facilities. Let  $\gamma_{i\ell}$ denote the fraction of $a_i$ that must be served from the $\ell$-th closest facility to customer $i$ where $\ell\in I=\{1,...,q\}$. To measure the service level of customer $i$ in a given solution, the distances from $i$ to the different open facilities are ordered and weighted with the values $\gamma_{i\ell}$ according to their rank in the sorted list of distances. This invites to characterize solutions by means of binary decision variables $x_{k\ell}^i$, $i, k\in P$, $\ell\in I$, where $x_{k\ell}^i$ is equal to 1 if $i$ is allocated to facility $k$ as the $\ell$-th closest facility.
Now, the service cost of customer $i$ can be computed as $s_i=\sum_{k\in P}\sum_{\ell\in I}a_i\gamma_{i\ell}d_{ik}x_{k\ell}^i$. Note that $s_i$ can be expressed in a compact way as \(s_i=\overline{C}^i \mathbf{x}^i\), where \(\mathbf{x}^i\) is the vector of decision variables $(x_{k\ell}^i)_{k\in P, \ell\in I}=(x_{11}^i,x_{12}^i,...,x_{21}^i,x_{22}^i,...)'$, and  $(\overline{C}^i)^\prime=(a_i\gamma_{il}d_{ik})_{k\in P, \ell\in I}$.\\

The VAOM operator is computed as a weighted sum of the service costs of all customers. A weight $\omega_j$ is applied to the customer with the $j$-th lowest service level, i.e. with the $j$-th highest service cost. For a given solution, $x$, and its associated vector $s$ as defined above, let $\sigma$ be a permutation of the indices of $P$ such that $s_{\sigma_1}\geq\hdots\geq s_{\sigma_p}$. Then,

$$VAOM_{(Q,d,\gamma,a,\omega)}(x)= \sum_{j=1}^p\omega_{j}s_{\sigma_j}.$$

The set of feasible solutions to the problem is fully characterized by the set of feasible assignments, since an explicit representation of the open facilities is not needed. These can be obtained directly from $x$ by identifying the indices $k\in P$ with $x_{k\ell}^i=1$ for some $i\in P$, $\ell\in I$. Thus in this problem $Q$ is given by the set of feasible assignments.
For reasons that will become evident when we cast the VAOM operator as an OWA, we express the assignment vectors $x$ as one dimensional $n$ vectors with $n=p^2q$.  In particular $x$ is partitioned in $p$ blocks, each of them associated with a different customer $i\in P$. That is, {$x=({\mathbf{x}^1}^\prime | \dots | \; {\mathbf{x}^i} ^\prime \; | \dots | {\mathbf{x}^p}^\prime)'$.} In turn, each block \(\mathbf{x}^i\) consists of $p$ smaller blocks, each with $q$ components. The $k$-th block of \(\mathbf{x}^i\) contains the $q$ components $x_{k\ell}^i$ for the indices $\ell\in I$.

Now, to cast the VAOM as an OWA operator,  we define $p$ objective functions {$\overline{C}^i\mathbf{x}^i$,} one associated with each customer $i\in P$. In particular, objective {$\overline{C}^i\mathbf{x}^i$} represents the service cost of customer $i\in P$, $s^i$. With the above characterization of vectors $x\in Q$, each {$\overline{C}^i$} must be defined by a $n$ vector. Thus expressing the VAOM as an OWA becomes basically a notation issue.
For each fixed $i\in P$, again we partition the cost vector {$\overline{C}^i$} in $p$ blocks. Similarly to the partition of vectors $x\in Q$, each block corresponds to a different customer, and has {$pq$} components. We now set at value 0 all the entries except those in the block of customer $i$, which are given by the entries of the vector $\overline{C}^i$ as defined above. That is: \(C^i=(\mathbf{0}_{pq} \; | \; \dots \; | \; \overline{C}^i \; | \; \dots \; | \;\mathbf{0}_{p q} )\), where $\mathbf{0_{pq}}=(0,...,0)\in \mathds{R}^{pq}$. With this notation it becomes clear that \(C^ix=\overline{C}^i\mathbf{x}^i\). Hence,

$$VAOM_{(Q,d,\gamma,a,\omega)}(x)=OWA_{(Q, C, \omega)}(x).$$

\vspace{12pt}
%---------------------------------------
\begin{example}\label{ex3}
Consider an instance of a VAOM problem with $p=3$ customers in which $q=2$ facilities must open. Suppose all the customers have one unit of demand, i.e. $a_1=a_2=a_3=1$, and suppose the rest of the data is the following:

$$
(d_{ik})_{i,k\in P}=\left(\begin{array}{ccc}
      0 & 2 & 6 \\
      2 & 0 & 4 \\
      8 & 4 & 0 \\
    \end{array}\right),\phantom{a}
    (\gamma_{il})_{i\in P, l\in I}=\left(\begin{array}{ccc}
      0.5 & 0.5 \\
      0.5 & 0.5 \\
      1   & 0   \\
    \end{array}\right),\phantom{a}
\omega '=\left(\begin{array}{ccc}
      0 & 1 & 2 \\
    \end{array}\right).
    $$

Since $q=2$ the feasible combinations of facilities to open are $\{1, 2\}$, $\{1, 3\}$ and $\{2, 3\}$. When the distances of each customer to the potential facilities are all different, like in this example, each combination of open facilities determines a unique  feasible assignment vector $x$. For instance, when facilities 1 and 2 open, then  customer 1, has facility 1 as the closest and facility 2 as the second closest, so $x_{11}^1=x_{22}^1=1$, and $x_{12}^1=x_{21}^1=x_{31}^1=x_{32}^1=0$. The service cost of customer 1 is thus $s_1=\gamma_{11}d_{11}x_{11}^1+\gamma_{12}d_{12}x_{22}^1=0+0.5\times 2=1$. For customer 2 we have $x_{12}^2=x_{21}^2=1$, and $x_{11}^2=x_{22}^2=x_{31}^2=x_{32}^2=0$, with service cost $s_2=0+0.5\times 2=1$. With this set of open facilities, the assignment for customer 3 is $x_{12}^3=x_{21}^3=1$, and $x_{11}^3=x_{22}^3=x_{31}^3=x_{32}^3=0$ with service cost $s_3=4$. Since $s_{3}\geq s_{1}\geq s_{2}$ the objective function value for this solution is thus $0\times 3 +1\times 1 +2 \times 1=3$.\\

Proceeding similarly with the other possible combinations of open facilities we obtain the complete set of feasible solutions $Q$, which in this example is given by the set of binary vectors given in Table \ref{Tabla4}:

\vspace{12pt}

\begin{table}[h!]
\begin{tabular}{|cc|cc|cc||cc|cc|cc||cc|cc|cc|}
\hline
$x_{11}^1$ & $x_{12}^1$ &   $x_{21}^1$ & $x_{22}^1$ & $x_{31}^1$  & $x_{32}^1$ &
$x_{11}^2$ & $x_{12}^2$ &   $x_{21}^2$ & $x_{22}^2$ & $x_{31}^2$  & $x_{32}^2$ &
$x_{11}^3$ & $x_{12}^3$ &   $x_{21}^3$ & $x_{22}^3$ & $x_{31}^3$  & $x_{32}^3$ \\
\hline
 1& 0   & 0 & 1   & 0 & 0     & 0 & 1   & 1 & 0   & 0 & 0     & 0&  1   & 1& 0  & 0  & 0\\
 1& 0   & 0 & 0   & 0 & 1     & 1 & 0   & 0 & 0   & 0 & 1     & 0&  1   & 0& 0  & 1  & 0\\
 0& 0   & 1 & 0   & 0 & 1     & 0 & 0   & 1 & 0   & 0 & 1     & 0&  0   & 0& 1  & 1  & 0 \\
 \hline
 \end{tabular}
\caption[]{Complete set of feasible solution $Q$ as binary vectors for Example \ref{ex3}.}\label{Tabla4}
\end{table}

%&&\notag

\vspace{12pt}

For modeling the VAOM as an OWA we define the cost matrix $C$ as:

$$C=
\left(\begin{array}{cc|cc|cc||cc|cc|cc|| cc|cc|cc}
0 & 0  & 1 & 1   & 3 & 3  &     0&  0    & 0 & 0   & 0 & 0  &      0 &  0  & 0 & 0   & 0 & 0  \\
0 & 0  & 0 & 0   & 0 & 0  &     1&  1    & 0 & 0   & 2 & 2  &      0 &  0  & 0 & 0   & 0 & 0   \\
0 & 0  & 0 & 0   & 0 & 0  &     0&  0    & 0 & 0   & 0 & 0  &      4 &  0  & 2 & 0   & 0 & 0
\end{array}\right).$$

Table \ref{Tabla5} shows the values of $y$,  $y_\pi$  and $OWA_{(Q, C, \omega)}(x)$ for each $x\in Q$. The optimal value of the VAOM is $ \min_{x\in Q}VAOM_{(Q,d,\gamma,a,\omega)}(x)=\min\{3,3,2\}=2$.

 \begin{table}[h!]
\renewcommand{\esp}{@{\hspace{0.15cm}}}
\begin{center}
\begin{tabular}{|@{}\esp c \esp |c\esp |c\esp |c\esp  @{}|}
\hline
   $y$ &  $y_\sigma$ & $OWA_{(Q, C, \omega)}(x)=\omega ' y_\sigma$\\
\hline											
$\left(\begin{array}{ccc}1 & 1 & 4\end{array}\right)'$	& $\left(\begin{array}{ccc}4 & 1 & 1\end{array}\right)'$	& 3\\
\hline
$\left(\begin{array}{ccc}3 & 3 & 0\end{array}\right)'$	& $	 \left(\begin{array}{ccc}3 & 3 & 0\end{array}\right)'$& 3\\
\hline
$\left(\begin{array}{ccc}4 & 2 & 0\end{array}\right)'$	& $\left(\begin{array}{ccc}4 & 2 & 0\end{array}\right)' $& 2\\
\hline
\end{tabular}
\caption[]{Values of $y$,  $y_\pi$ and $OWA_{(Q, C, \omega)}(x)$ for the feasible solutions of Example \ref{ex3}.}\label{Tabla5}
\end{center}
\end{table}

\end{example}
%---------------------------------------

%\item \textbf{The Vector Assignment Ordered Median function of an abstract combinatorial optimization problem}
\subsection{The Vector Assignment Ordered Median function of an abstract combinatorial optimization problem}

In the above section we have applied the VAOM operator to the locations and allocations of a general multifacility location problem, according to the original definition by \citet{Lei2012}. Nevertheless, this operator can be also applied to the characteristic vector of a combinatorial solution of any abstract combinatorial optimization problem, as we also did with the ordered  median operator. In doing that we obtain a more general interpretation of this type of  objective function that can also be cast within the OWA operator.

Let $Q\subseteq \mathbb{R}^n$ denote the feasible domain for an optimization problem, $\omega\in \mathbb{R}^{p+}$ a given vector of nonnegative weights and $P=\{1,\ldots,p\}$. Recall that a VAOM operator considers for each objective function $s_i$, $i\in P$ different  fractions, $\gamma^{i}$, of the cost vector $d$ for the sorted elements of the decision vector $x$.

For $x\in Q$,  the evaluation of the $i$-th component of the VAOM objective {is given by $s_i=\gamma^i d_i x_i$,} for all $i \in P$.   Let $\sigma$ denote a permutation of the indices of $P$, such that $s_{\sigma_i}\ge  s_{\sigma_{i+1}}$, for $i=1,\ldots, p-1$.
The VAOM operator is $VAOM_{(Q,d,\gamma,\omega)}(x)= \sum\limits_{i\in P}\omega_i s_{\sigma_i}$. The reader may note that the original definition of VAOM can be accommodated to this general setting once we identify the combinatorial object $Q$ as the set of  location-allocations in the discrete location problem. In that case, there are $i=1,\ldots,n$ objective functions associated with each of the customers and then the fractions that apply to each customer $i$ are {non-null} only for a subset of the open facilities (servers) corresponding to the $q$-closest ones.

This can be done by defining a set of variables, one per customer $i$, with $n$ blocks. In the block $k$, $x^ i_{\centerdot k}=(x^i_{1k},...,x^i_{nk})'$ accounts for the allocation of $i$ to any facility as the $k$-th closest, therefore $x^i=(x^{i\prime}_{\centerdot 1}\; | \;x^{i\prime}_{\centerdot 2}\; | \;...\; | \;x^{i\prime}_{\centerdot n})^\prime$ for $i=1,\ldots,n$. {This way,} the cost vectors must also have the same structure by blocks, each block corresponding with the level of assignment, i.e. denoting by {$d^{i}_{\centerdot}=(d^i_{1},...,d^i_{n})'$} then {$d^{i}=(d^{i\prime}_{\centerdot}\; | \;d^{i\prime}_{\centerdot}\; | \;...\; | \;d^{i\prime}_{\centerdot})'$.} Finally, since  the fractions of costs are applied according to the level of assignment, the structure of the vector of fractions $\gamma^i$ is also by blocks. Block $k$ represents the fraction of the cost that is accounted for costumer $i$ at the $k$-th level of assignment. Denoting by  i.e. {$\gamma^{i}_\ell=(\gamma_{i\ell},...,\gamma_{i\ell})'$} then
{$\gamma^{i}=(\gamma^{i}_1\; | \;\gamma^{i\prime}_2\; | \;...\; | \;\gamma^{i\prime}_n)'$}  for $i=1,\ldots,n$.

To cast the VAOM  as an OWA operator, we only need to set $\bar{\gamma}^{i}=(\underbrace{\mathbf{0}_{np}}_{1} \; | \; \dots \; | \; \underbrace{{\gamma^{i}}'}_{i} \; | \; \dots \; | \;\underbrace{\mathbf{0}_{np}}_{p} )'$,
$\bar{d}^{i}=(\underbrace{\mathbf{0}_{np}}_{1} \; | \; \dots \; | \; \underbrace{{d^i}'}_{i} \; | \; \dots \; | \;\underbrace{\mathbf{0}_{np}}_{p} )'$,
$x=({\mathbf{x}^1}^\prime | \dots | \; {\mathbf{x}^i} ^\prime \; | \dots | {\mathbf{x}^p}^\prime)'$ and
$C^i= (\bar{\gamma}^i _j  \bar{d}^i_j)_{j=1}^{p(pn)}$.
Then, the VAOM can be written as the following OWA operator
$VAOM_{(Q,d,  \gamma ,\omega)}(x)=OWA_{(Q, C, \omega)}(x)$.\\

%\end{itemize}

As we have shown above, OWA is a very general operator. In the following, we will work in more particular settings, namely we shall restrict ourselves to assume that $Q$ is a combinatorial object which can be represented by a system of linear inequalities.
%=======================================================================================================================
%\newpage
%=======================================================================================================================
\section{Basic formulations for the OWAP, properties and reinforcements}\label{Sec:ModelingtheOWAP}
This section presents alternative Mixed Integer Programming (MIP) formulations for an OWAP, which are analyzed and compared. The starting point of our study are three basic formulations, which, broadly speaking, differ from one to another on how the permutation that defines the ordering of the cost function values is modeled.  Two of {the formulations presented} use binary variables $z$ to define the \textit{specific positions} in the ordering of the sorted cost function values, whereas the other one uses binary variables $s$ to define the \textit{relative position} in the ordering of the sorted cost function values. One of the two formulations based on the $z$ variables also uses an additional set of decision variables $y$ for modeling the specific values of the cost functions depending on their position in the ordering. All three formulations use a set of decision variables $\theta$ to compute the values of the objectives sorted at specific positions. In each case, alternative formulations are presented, which preserve the set of optimal solutions. Before addressing any concrete formulation we discuss the meaning of both sets of variables $z$ and $s$ as well as their relationships.

%------------------------------------------------------------------
\subsection{Alternative formulations for permutations\label{rel:z-s}}

The essential element in our formulations rests on the representation of ordering within a MIP model. To such end, we devote this section to describe how  a permutation can be represented with binary variables. Recall that we  have introduced $P=\{1,\dots,p\}$ as the set of the cost function indices.
Let $\pi:P \rightarrow P$  be a function representing a permutation of $P$. That is, it assigns the index $i$ of each cost function (also denoted by \emph{cost function $i$}) to a position indexed by $j$ (also denoted by \emph{position $j$}). Note that $\pi$ is a permutation if each cost function is assigned to a single position and if each position contains a single cost function index. In what follows, we use {$\pi_i=\pi(i)$} to denote the position occupied by cost function {$i\in P$} and {$\sigma_j=\pi^{-1}(j)$} to denote the index of the cost function that occupies position $j$ (we recall that the notation $\sigma$ was previously used in Section \ref{Sec:Definition}). Note that $\sigma$ also defines a permutation of the positions of $P$. In what follows we will indistinctively use $\pi$ and $\sigma$. Slightly abusing notation we refer to $\pi$ as to \emph{the cost functions permutation} and to its inverse $\sigma$ as to the \emph{positions permutation}.\\

In order to model $\pi$ as a permutation, let $z_{ij}$ be a binary decision variable defined as
\[
z_{ij}=\begin{cases}
 1 & \text{if cost function $i$ occupies position $j$, (i.e. if $\pi_i=j$)}\\
 0 & \text{otherwise.}
\end{cases}
\]
The set of variables $z$ defines a permutation if:
\begin{itemize}
\item[($i$)] each position contains a single cost function:
\formulaN{\sum_{i\in P}z_{ij}=1                                      &&  j\in P, \label{permutationz11}}
and,
\item[($ii$)] each cost function $i$ is assigned to a single position $j$:
\formulaN{\sum_{j\in P}z_{ij}=1                                      &&  i\in P. \label{permutationz12}}
\end{itemize}
In addition, we observe that since system \eqref{permutationz11}-\eqref{permutationz12} contains exactly $2p-1$ linearly independent equations, the above permutation can also be  represented without variables $z_{i1}$, for all $i\in P$, that can be replaced by $1-\sum_{j\in P: j>1} z_{ij}$. In this way, system \eqref{permutationz11}-\eqref{permutationz12} can also be rewritten as

\formulaN{\sum_{i\in P}z_{ij}=1                                      &&  j\in P: j>1, \label{permutationz21}}
\formulaN{\sum_{j\in P}z_{ij}\leq1                                      &&  i\in P. \label{permutationz22}}

\begin{example}
Let $\pi$ be a permutation defined by $\pi=\left(\begin{array}{cccc}3 & 2 & 4 & 1\end{array}\right)$ or equivalently by $\sigma=\left(\begin{array}{cccc}4 & 2 & 1 & 3\end{array}\right)$. Then, $\pi$ can be represented by using variables $z$ as follows:
$$ (z_{i,j})_{i,j\in P}=\left(\begin{array}{cccc}
      0 & 0 & 1 & 0 \\
      0 & 1 & 0 & 0 \\
      0 & 0 & 0 & 1 \\
      1 & 0 & 0 & 0 \\
    \end{array}\right), \, \textrm{or} \,
(z_{i,j})_{i,j\in P: j>1}=\left(\begin{array}{ccc}
      0 & 1 & 0 \\
      1 & 0 & 0 \\
      0 & 0 & 1 \\
      0 & 0 & 0 \\
    \end{array}\right).$$
\end{example}
\hfill$\square$\\

An alternative representation of a permutation, which we have also found useful is based on a different set of variables defined as:

\[
s_{ij}=\begin{cases}
 1 & \text{if cost function $i$ is placed before position $j$ in the ordering,}\\
 0 & \text{otherwise.}
\end{cases}
\]
The set of variables $s$ defines a permutation if:

\begin{itemize}
\item[($i$)] for all $j\in P$ there are $j-1$ cost functions placed before position $j$:
\formulaN{\sum_{i\in P}s_{ij}= j-1                                       &&  j\in P, \label{permutations11}}
and
\item[($ii$)] cost function $i$ cannot be placed  {in position $j$ unless it is also placed in position $j+1$, i.e.,}

\formulaN{s_{ij+1}-s_{ij}\geq 0                                          &&  i,j\in P:j<p. \label{permutations12}}
\end{itemize}
Again we can reduce the number of decision variables, now by eliminating $s_{i1}$ for all $i\in P$.  Indeed,  since there is no cost function placed before position 1 in any ordering, all the $s_{i1}, i\in P$ can be fixed to zero. In this way, permutation \eqref{permutationz21}-\eqref{permutationz22} can be also represented by means of the following reduced set of constraints:
\formulaN{\sum_{i\in P}s_{ij}= j-1                                       &&  j\in P: j>1, \label{permutations21}}
\formulaN{s_{ij+1}-s_{ij}\geq 0                                          &&  i,j\in P:1<j<p. \label{permutations22}}

\begin{example}
Let $\pi$ be a permutation defined by $\pi=\left(\begin{array}{cccc}3 & 2 & 4 & 1\end{array}\right)$ or equivalently by {$\sigma=\left(\begin{array}{cccc}4 & 2 & 1 & 3\end{array}\right)$}. Then, $\pi$ can be represented by using variables $s$ as follows:
$$(s_{i,j})_{i,j\in P}=\left(\begin{array}{cccc}
      0 & 0 & 0 & 1 \\
      0 & 0 & 1 & 1 \\
      0 & 0 & 0 & 0 \\
      0 & 1 & 1 & 1 \\
    \end{array}\right), \, \textrm{or} \,
(s_{i,j})_{i,j\in P: j>1}=\left(\begin{array}{ccc}
      0 & 0 & 1 \\
      0 & 1 & 1 \\
      0 & 0 & 0 \\
      1 & 1 & 1 \\
    \end{array}\right). $$
\hfill$\square$
\end{example}

With the above considerations, variables $z$ and $s$ are related by means of
\formulaN{z_{ij}=\begin{cases} \label{for:s-z-1}
 s_{ij+1}-s_{ij} &  i\in P, j=1,...,p-1 \\
 1- s_{ij} &  i\in P, j=p
\end{cases}}
and equivalently,
%\formulaN{s_{ij}=1-\sum_{j'\geq j}z_{ij'},\, i,j\in P.}
\formulaN{s_{ij}=1-\sum_{k\geq j}z_{ik},\, i,j\in P. \label{for:s-z-2}}

%---------------------------------

%=======================================================================================================================
%\newpage
%=======================================================================================================================
\subsection{$OWAP$ formulations with variables for the positions of sorted cost function values}\label{subsec:1}
%\subsection{Formulations with variables for the specific positions of sorted cost function values}\label{subsec:1}

For a given feasible set $Q \subseteq \mathbb{R}^n$, consider the binary decision variables $z$ as defined in Section \ref{rel:z-s} to represent the permutation $\pi$ associated with the sorted cost function values $C^ix$, $i\in P$. Let also $\theta_j$ be a real decision variable equal to the value of the cost function sorted in position $j$.  Next, we give an integer linear programming description of the $OWAP$  where we use $M$ to denote a non-negative upper bound of the value of all the cost functions. (We refer the interested reader to \cite{Boland2006} or \cite{Nickel2005} for similar sets of decision variables and formulations for the discrete ordered median location problem.)

%\begin{subequations}
%\vspace{-0.5cm}
%\begin{align}
%&\Hp \hspace{8cm}&&\notag\\
%&\hhp F^z_0: \hspace{0.5cm} V=\min\sum_{j\in P}\omega_j\theta_j  &&                  \tag{\rm 1a}    \label{owaa}\\ s.t.
%&\Hp \sum_{i\in P}z_{ij}=1                                      &&  j\in P          \tag{\rm 1b}    \label{owab}\\
%&\Hp \sum_{j\in P}z_{ij}=1                                      &&  i\in P          \tag{\rm 1c}    \label{owac}\\
%&\Hp C^i x\leq\theta_j+{M}(1-z_{ij})                            &&  i,j\in P        \tag{\rm 1d$^0$}\label{owad0}\\
%&\Hp \theta_j \geq \theta_{j+1}                                 &&  j\in P:j<p      \tag{\rm 1e}    \label{owae}\\
%&\Hp x\in Q, z\in\{0,1\}^{p\times p}                            &&                  \tag{\rm 1f}    \label{owaf}
%\end{align}
%\end{subequations}

\newcounter{counterequation}
\addtocounter{counterequation}{1}
\addtocounter{counterequation}{\value{equation}}
\begin{subequations}
\vspace{-0.5cm}
\begin{align}
&\Hp \hspace{8cm}&&\notag\\
&\hhp F^z_0: \hspace{0.5cm} V=\min\sum_{j\in P}\omega_j\theta_j  &&                 \label{owaa}\\ s.t.
&\Hp \sum_{i\in P}z_{ij}=1                                      &&  j\in P          \label{owab}\\
&\Hp \sum_{j\in P}z_{ij}=1                                      &&  i\in P          \label{owac}\\
&\Hp C^i x\leq\theta_j+{M}(1-z_{ij})                            &&  i,j\in P         \tag{\arabic{counterequation}d$^0$}\label{owad0}\\\addtocounter{equation}{1}
&\Hp \theta_j \geq \theta_{j+1}                                 &&  j\in P:j<p      \label{owae}\\
&\Hp x\in Q, z\in\{0,1\}^{p\times p}                            &&                  \label{owaf}
\end{align}
\end{subequations}

The objective function \eqref{owaa} minimizes the  weighted average of sorted objective function values, provided that $\theta_j$, $j\in P$, are enforced to take on the appropriate values.
As seen, constraints \eqref{owab}-\eqref{owac} define a cost functions permutation by placing at each position of {$\pi$} a single cost function and each cost function at a single position of {$\pi$.}
Constraints \eqref{owad0} relate cost function values with the values placed  in a sorted sequence.
Constraint \eqref{owae} imposes that the sorted values are ordered non-increasingly.

%------------------------------------------------------------------

In the following we denote by  $\Omega^z_0$ the domain of feasible solutions to formulation $F^z_0$. That is, $$\Omega^z_0=\left\{(x,z,\theta) \text{ satisfying constraints } \text{\text{\eqref{owab}, \eqref{owac}, \eqref{owad0}, \eqref{owae}, \eqref{owaf}}}\right\}.$$ Consider now the family of inequalities
%\formulaN{C^ix\leq \theta_j+M_{i}(1-\sum_{j'\geq j}z_{ij'})                 &&  i,j\in P,   \tag{\rm 1d} \label{owad'}}
\formulaN{C^ix\leq \theta_j+M(1-\sum_{k\geq j}z_{ik})                 &&  i,j\in P,   \tag{\arabic{counterequation}d} \label{owad}}

\noindent and note that, for $z$ satisfying \eqref{owac}, inequalities \eqref{owad} can be rewritten as

%\formulaN{C^ix\leq \theta_j+M_{i}\sum_{j'< j}z_{ij'}                &&  i,j\in P,    \tag{\ref{owad'}'} \label{owad''}}%\nonumber
\formulaN{C^ix\leq \theta_j+M\sum_{k< j}z_{ik}                &&  i,j\in P,    \tag{\arabic{counterequation}d'} \label{owad'}}%\nonumber

%\noindent since for all $i, j\in P$, $1-\sum_{j'\geq j}z_{ij'}=\sum_{j'<j}z_{ij'}$.\\
\noindent since for all $i, j\in P$, $1-\sum_{k\geq j}z_{ik}=\sum_{k<j}z_{ik}$.\\

\begin{remark}
Observe that when variables $z_{i1}, \, i\in P$ are not defined and the permutation is described by means of inequalities \eqref{permutationz21} and \eqref{permutationz22}, then constraints  \eqref{owad0}, \eqref{owad} and \eqref{owad'} must consider separately the case $j=1$ from the case $j\in P, j>1$. In particular, the case $j=1$ reduces to
\formulaN{C^ix\leq \theta_1                &&  i\in P,   }
since the first position has always a value greater than or equal to any cost function.\newline
\end{remark}

Let $\Omega^{z}=\{(x,z,\theta) \text{ satisfying constraints } \text{\eqref{owab}, \eqref{owac}, \eqref{owad}, \eqref{owae}, \eqref{owaf}}\}$ denote the domain obtained from $\Omega^{z}_0$ when constraints \eqref{owad0} are replaced by constraints \eqref{owad}.

\begin{property}\label{propertyowa1sumas}
$\Omega^{z}_0=\Omega^{z}$.
\end{property}

\texttt{Proof}.\newline
It is clear that $\Omega^{z}_0\supseteq\Omega^{z}$, since for $i,j\in P$ given, the right hand side of the associated constraint \eqref{owad} is smaller than or equal to that of constraint \eqref{owad0}.

To prove that $\Omega^{z}_0\subseteq\Omega^{z}$ also holds let $(x,z,\theta)\in \Omega^{z}_0$  and we show that $(x,z,\theta)$ satisfies constraints  \eqref{owad}. For $i,j\in P$ given, we distinguish two cases:
\begin{itemize}
  \item If $z_{ij}=1$ then \eqref{owad} holds for this pair of indices.
  \item If $z_{ij}=0$ then by \eqref{owac}, there must exist $ j'\in P$, $j'\neq j$, such that $z_{ij'}=1$. If $j'< j$, then $\sum_{k\ge j}z_{ik}=z_{ij}=0$, and \eqref{owad} holds for the pair of indices $i, j$. Otherwise, if $j'> j$, then $\sum_{k\ge j}z_{ik}=z_{ij'}=1$ so the right hand side of constraint \eqref{owad} for the pair $i,j$ takes the value $\theta_{j}$. Now constraint \eqref{owad0} for the pair of indices $i, j'$ implies that $C^ix\le \theta_{j'}$. By constraints \eqref{owae}, we also have $\theta_{j}\ge \theta_{j'}$ and thus \eqref{owad} also holds for the pair of indices $i, j$. \hfill$\square$
\end{itemize}
%\hfill$\square$\\

\begin{remark}
Since $\Omega^{z}_0=\Omega^{z}$, an equivalent formulation for the OWAP is
\begin{subequations}
\vspace{-0.5cm}
\begin{align}
&\Hp \hspace{8cm}&&\notag\\
&\hhp F^{z}: \hspace{0.5cm} V=\min\sum_{j\in P}\omega_j\theta_j  && \notag \\ s.t.
&\Hp (x,z,\theta)\in \Omega^{z}. \notag &&
\end{align}
\end{subequations}
Formulation $F^{z}$ can be preferred to formulation $F^{z}_0$ for solving an OWAP, since it may provide tighter linear programming bounds, given that, for fractional vectors $z$ satisfying constraints \eqref{owab}-\eqref{owac}, constraints \eqref{owad0} are dominated by constraints \eqref{owad}. \\
\end{remark}

In the search for optimal solutions to the OWAP any formulation whose optimal solution set coincides with that of the OWAP can be of interest. Such formulations could be preferred because they use fewer variables or constraints, or because their feasible domain has a structure which is easier to explore. Next we present three such formulations. All of them can be seen as relaxations of formulation $F^z$ in the sense that their feasible domains contain $\Omega^z$. However, all of them are valid formulations for the OWAP since they preserve the set of optimal solutions of $F^z$, i.e. their set of optimal solutions coincides with that of $F^z$. First we prove a property of optimal solutions.\\

\begin{lemma}
Let $(x^*, z^*,\theta^*)\in \Omega^z$ be an optimal solution to $F^z$. Then for each $j\in P$ there exists $i\in P$ with $\theta^*_j=C^ix^*$.
\end{lemma}

\texttt{Proof}.\newline
Let $\tilde x$ be a feasible solution in $Q$. Then, there exists a positions permutation $\sigma$ that sorts the cost functions values in non-increasing order. That is, $C^{\sigma_j}\tilde x\geq C^{\sigma_{j+1}}\tilde x, \forall j\in P\setminus \{p\}$. Therefore, we can set $\tilde  z=(z_{\sigma_j,j})_{j\in P}$ and $\theta=(C^{\sigma_j}\tilde  x)_{j\in P}$. Since this is true for each $x\in Q$, it is true in particular for $x^*$.
\hfill$\square$\newline

From above lemma, we observe that $\Omega^z$ is always non empty, provided that $Q$ is non empty.
%---------------------------------
\bigskip

Let $\Omega^{z}_{R1}=\{(x,z,\theta) \text{ satisfying constraints } \text{\eqref{owab}, \eqref{owac}, \eqref{owad}, \eqref{owaf}}\}$, i.e,  $\Omega^{z}_{R1}$ is the relaxation of the domain $\Omega^{z}$ once the set of constraints \eqref{owae} is removed. Next,  consider the formulation

\begin{subequations}
\vspace{-0.5cm}
\begin{align}
&\Hp \hspace{8cm}&&\notag\\
&\hhp F^{z}_{R1}: \hspace{0.5cm} V=\min\sum_{j\in P}\omega_j\theta_j  && \notag \\ s.t.
&\Hp  (x, z,\theta)\in \Omega^{z}_{R1}. \notag &&
\end{align}
\end{subequations}

\begin{lemma}\label{lemaFz_R1}
Every feasible solution to $F^{z}_{R1}$, $(x,z,\theta)\in \Omega^{z}_{R1}$, satisfies $ \theta_{i}\geq \max\{C^{\sigma_i}x,C^{\sigma_{i+1}}x\}, \; i=1,\ldots,p-1$ and
 $\theta_{p}\geq C^{\sigma_p}x$.
\end{lemma}

\texttt{Proof}.\newline
Let $(x,z,\theta)\in \Omega^{z}_{R1}$ be a feasible solution to $F^{z}_{R1}$ and $\sigma$ a permutation that sorts the cost function values of $x$.
For $i\in P$ given, $z_{\sigma_i,i}=1 $. Then, by \eqref{owad} we have that $\theta_{j}\geq C^{\sigma_i}x$, for $j\le i$ and, in particular,
\begin{equation}
\theta_{i}\geq C^{\sigma_i}x.  \label{proof11}
\end{equation}

When $i\le p-1$, the same argument can be applied to $z_{\sigma_{i+1},i+1}=1$, getting $\theta_{j}\geq C^{\sigma_{i+1}}x$, for $j\le i+1$ and, in particular,
\begin{equation}
\theta_{i}\geq C^{\sigma_{i+1}}x \text{ and } \theta_{i+1}\geq C^{\sigma_{i+1}}x. \label{proof12}
\end{equation}

Using \eqref{proof11} and \eqref{proof12} we obtain the result.
\hfill$\square$\newline

\begin{property}\label{propertyowa12}
Every optimal solution to  $F^{z}_{R1}$ is also optimal to $F^z$.
\end{property}

\texttt{Proof}.\newline
Since $\Omega^{z}\subseteq \Omega^{z}_{R1}$it is enough to prove that any optimal solution to $F^{z}_{R1}$ is feasible to $F^z$.
Let $(x,z,\theta)\in \Omega^{z}_{R1}$ be an optimal solution to $F^{z}_{R1}$ and $\sigma$ a permutation that sorts the cost function values of $x$. Let us see that $\theta$ verifies constraint \eqref{owae}.

By Lemma \ref{lemaFz_R1} we have that $ \theta_{i}\geq \max\{C^{\sigma_i}x,C^{\sigma_{i+1}}x\}, \; i=1,\ldots,p-1$ and  $\theta_{p}\geq C^{\sigma_p}x$.

Since we are minimizing a function  which is a linear combination with non-negative weights of the $\theta$ variables, it follows that in any optimal solution $\theta_j \geq \theta_{j+1}, \, j\in P\backslash\{p\}$ since, otherwise, the value of $\theta_{j+1}$ could be decreased to $\theta_{j}$, while keeping all other variables values unchanged, improving the objective function value. That is, \eqref{owae} holds.
\hfill$\square$\newline

%------------------------------------------------------------------

Consider now $\Omega^{z}_{R2}=\{(x,z,\theta) \text{ satisfying constraints } \text{\eqref{owab}, \eqref{owad}, \eqref{owaf}}\}$, i.e,  $\Omega^{z}_{R2}$ is the relaxation of the domain $\Omega^{z}_{R1}$ once  the set of constraints \eqref{owac} is removed. Next,  consider the formulation

\begin{subequations}
\vspace{-0.5cm}
\begin{align}
&\Hp \hspace{8cm}&&\notag\\
&\hhp F^{z}_{R2}: \hspace{0.5cm} V=\min\sum_{j\in P}\omega_j\theta_j  && \notag \\ s.t.
&\Hp  (x,z,\theta)\in \Omega^{z}_{R2}. \notag &&
\end{align}
\end{subequations}

\begin{property}\label{propertyowa13}
Every optimal solution to $F^{z}_{R2}$ is also optimal to $F^{z}$.
\end{property}

\texttt{Proof}.\newline
Since $\Omega^{z}\subseteq \Omega^{z}_{R2}$ it is enough to prove that any optimal solution to $F^{z}_{R2}$ is feasible to $F^z$.
Let $(x,z,\theta)$ be an optimal solution to $F^{z}_{R2}$.  If $(x,z,\theta)$ is optimal to $F^{z}_{R1}$ then, by using Property \ref{propertyowa12},
$(x,z,\theta)$ is also optimal to $F^{z}$. Thus, to prove that $(x,z,\theta)$ is optimal to $F^{z}$, it suffices to prove that $(x,z,\theta)$ satisfies inequalities \eqref{owac}. \\

%\rv{Miguel: Creo que la demostración podría ser más sencilla analizando \eqref{owac} %por casos.}

We prove first that $\sum_{j\in P}z_{ij}\leq 1$ for all $i\in P$. Using the notation $r_{ij}=\sum_{k\ge j}z_{ik}$, for all $i, j\in P$, constraints  \eqref{owad} can be rewritten as

$$C^{i}x\leq\theta_{j} +M(1-r_{ij}) \Leftrightarrow \theta_{j} \geq C^{i}x + M(r_{ij}-1).$$

Therefore, for all $j\in P$, $$ \theta_{j} = \max_{i\in P}\{C^ix+M(r_{ij}-1)\}.$$

Suppose there exists $i'\in P$ with $\sum_{j\in P}z_{i'j}=r>1$, and let $j'=\arg\max\{r_{i'j}=2\mid j\in P\}$. If several indices exist with $\sum_{j\in P}z_{ij}>1$ we select $i'$ as the one with maximum associated $j'$.

The criterion for the selection of $i'$ and the definition of $j'$ imply that $r_{i'j'}=2$ and $r_{ij'}\le 1$ for all $i\ne i'$.

Therefore, since $M$ is a strict upper bound on the value of any cost function, the actual value of $\theta_{j'}$ is determined by cost function $i'$, and we have

$$ \theta_{j'}=C^{i'}x + M(r_{i'j'}-1)=C^{i'}x + M.$$

Also, $r_{i'j}\ge2$ for all $j<j'$. Thus, $\theta_j\ge C^{i'}x + M$ for all $j<j'$. Furthermore, $r_{ij}\le 1$ for all $i\in P$, $j>j'$, implying that $\theta_j < M$ for all $j>j'$.\\

Observe, on the other hand, that  $\sum_{j\in P}z_{i'j}>1$ implies that there exists some $i''\in P$, $i''\ne i'$ with $\sum_{j\in P}z_{i''j}=0$. (Otherwise, adding  up all constraints \eqref{owab} we get a
contradiction.)

Let us now define the solution $(x,\overline{z},\overline{\theta}) \in \Omega^{z}_{R2}$ with the same $x$ components as above, where

\[\overline{z}_{ij}=\begin{cases}
0   &   \text{if } i=i', \text{ and } j=j'\\
1   &   \text{if } i=i'', \text{ and } j=j'\\
z_{ij} &        \text{otherwise.}
\end{cases}\]

It is clear that $\sum_{j\in P}\overline{z}_{i'j}=r-1$, and, $\sum_{k\ge j}\overline{z}_{i'k}=r_{i'j}-1$, for all $j\in P$.
It is also clear that $\sum_{j\in P}\overline{z}_{i''j}=1$, and, $\sum_{k\ge j}\overline{z}_{i''k}=1$, for all $j\le j'$, and 0 for $j>j'$. For all other $i\ne i', i''$, it holds that $\sum_{j\in P}\overline{z}_{i'j}=\sum_{j\in P}z_{i'j}$.
Since $\sum_{k\ge j'}\overline{z}_{ik}\le 1$, for all $i\in P$ we now have

$$ \overline{\theta}_{j'}=\max_{i\in P}\{C^{i}x\} < M \le C^{i'}x + M = \theta_{j'},$$

and, $\overline{\theta}_j\le \theta_j$, for all $j\ne j'$.

Therefore, since we are minimizing a linear function with non-negative weights of the $\theta$ variables, the objective function value of $(x,\overline{z},\overline{\theta})$ is smaller than that of $(x,z,\theta)$,
contradicting the optimality of $(x,z,\theta)$. Hence,  $\sum_{j\in P}z_{ij}\leq 1$ for all $i\in P$.

Let us, finally, see that  $\sum_{j\in P}z_{ij}\ne 0$ for all $i\in P$. Assume on the contrary that $\sum_{j\in P}z_{i'j}=0$ for some $i'\in P$. Then, by adding  up all constraints \eqref{owab} we get $p=\sum_{j\in P}\left(\sum_{i\in P}z_{ij}\right)=\sum_{i\in P}\left(\sum_{j\in P}z_{ij}\right)=\sum_{i\in P,\substack{i\ne i'}}\left(\sum_{j\in P}z_{ij}\right)\le p-1$, which is impossible.
\hfill$\square$\\

%---------------------------------

%---------------------------------

We now consider the inequality version of constraints \eqref{owab}
\formula{\sum_{i\in P}z_{ij}\le 1                              \qquad j\in P.         \tag{\ref{owab}$_{\le}$} \label{owab<}}

\begin{remark}\label{remark:d-d'}
Observe that when inequalities \eqref{owab<} hold, constraints \eqref{owad} are no longer equivalent to \eqref{owad'}. \\
\end{remark}

Let us define the domain $\Omega^{z}_{R3}=\{(x,z,\theta) \text{ satisfying constraints } \text{\eqref{owab<}, \eqref{owad'}, \eqref{owaf}}\}$.

 It is clear that $\Omega^{z}\subseteq \Omega^{z}_{R3}$. However, as we next see, both sets are equivalent for the minimization of the objective \eqref{owaa} in the sense that they define the same set of optimal solutions. Consider the problem

\begin{subequations}
\vspace{-0.5cm}
\begin{align}
&\Hp \hspace{8cm}&&\notag\\
&\hhp F^{z}_{R3} \hspace{0.5cm} V=\min\sum_{j\in P}\omega_j\theta_j  && \notag \\ s.t.
&\Hp  (x,z,\theta)\in \Omega^{z}_{R3}. \notag &&
\end{align}
\end{subequations}

%------------------------------------------------------------------

%------------------------------------------------------------------

\begin{lemma}
$\Omega^{z}_{R2}\subseteq \Omega^{z}_{R3}$.
\end{lemma}
\texttt{Proof}.\newline
We prove that any feasible solution $(x,z,\theta)\in \Omega^{z}_{R2}$ verifies that $(x,z,\theta)\in \Omega^{z}_{R3}$. To prove this, it is only necessary to prove that $(x,z,\theta)$ verifies \eqref{owad'}. From \eqref{owad} we have that $(x,z,\theta)$ verifies
\formulaN{\theta_j\geq \max_i\{C^ix-M(1-\sum_{k\geq j}z_{ik})\},\, j\in P\label{lem:lema6JP}}
 and for \eqref{owad'}, we have to prove that $(x,z,\theta)$ also verifies
\formulaN{\theta_j\geq \max_i\{C^ix-M(\sum_{k< j}z_{ik})\},\, j\in P. \label{lemar2r32}}
We distinguish the following cases:
\begin{itemize}
  \item If $\sum_{k\geq j}z_{i'k}=r>1$ for some $i'$ then
  \formulaN{\theta_j\geq C^{i'}x+(r-1)M \geq \max_{i}\{C^ix-M(\sum_{k< j}z_{ik})\},\label{lemar2r33}}
  and the result holds.
  \item If $\sum_{k\geq j}z_{ik}=1$ for all $i\in P$ then $\theta_j\geq \max_i\{C^ix\} \geq \max_i\{C^ix-M(\sum_{k< j}z_{ik})\}$ and the results is also proven.
  \item If $\sum_{k\geq j}z_{i'k}=0$ for some $i'$ then we distinguish to subcases. If $\sum_{k< j}z_{i'k}\geq 1$ then from \eqref{lem:lema6JP} we easily get that \eqref{lemar2r32} holds. Otherwise, $\sum_{k\in P}z_{i'k}=0$ and by \eqref{owab} it does exist an $i''$ such that $\sum_{k\geq j}z_{i''k}=r>1$. Thus, by using \eqref{lemar2r33}, equation \eqref{lemar2r32} also holds.
\end{itemize}
\hfill$\square$\\

%------------------------------------------------------------------

%------------------------------------------------------------------
\begin{property}\label{propertyowa14JP}
$F^{z}$ and $F^{z}_{R3}$ have the same set of optimal solutions.
\end{property}

\texttt{Proof}.\newline

Since $\Omega^z \subset \Omega_{R2}^z$ and $\Omega_{R2}^z\subset \Omega_{R3}^z$ then $\Omega^z \subset \Omega_{R3}^z$ and it is enough to prove that any optimal solution to $F^{z}_{R3}$ is feasible to $F^{z}$. Since the set of optimal solutions of $F^{z}$ and $F^{z}_{R2}$ coincide, we only need to prove that any optimal solution of $F^{z}_{R3}$ is feasible for $F^{z}_{R2}$.

To see that any optimal solution $(x,z,\theta)$ to $F^{z}_{R3}$ is feasible to $F^{z}_{R2}$, it is enough to see that $(x,z,\theta)\in \Omega^{z}_{R2}$, i.e. it satisfies inequalities \eqref{owab} and \eqref{owad}.\\

By a similar argument to the one applied in Property \ref{propertyowa13},  any optimal solution   $(x,z,\theta)$ of $F^z_{R3}$  satisfies {$\sum_{j\in P} z_{ij}= 1$.} Therefore, satisfying inequality \eqref{owad'} implies inequality \eqref{owad}.\\

To see that {$(x,z,\theta)$} also satisfies \eqref{owab}, let us suppose w.l.o.g. that there exists exactly one $j'\in P$ such that $\sum_{i\in P}z_{ij'}=0$. Then, by adding up all constraints \eqref{owab<} we have {$p-1\geq \sum_{j\in P}\sum_{i\in P}z_{ij}=\sum_{i\in P}\sum_{j\in P}z_{ij}$.} Therefore, there must exist $i'\in P$ such that $\sum_{j\in P}z_{i'j}=0$. Thus, we observe that we can construct $(x,\bar z,\bar \theta)$, another optimal solution to $F^z_{R3}$, setting $\bar{z}_{ij}=z_{ij}$, if $i\neq i'$ and $\bar{z}_{i'k}=1$ for any $k$. Clearly, $(x,\bar z,\bar \theta)$ is a feasible solution to $F^z_{R3}$ for some suitable $\bar \theta$, satisfying in addition
$$ C^{i'}x\le \bar \theta_k+M\sum_{\ell<j} z_{i'\ell}, \qquad \forall k\in P.$$
Therefore, this inequality allows for any $k\in P$ that $\bar \theta_k$  assumes a value smaller than or equal to $\theta_k$, the  one associated with the solution $(x, z, \theta)$, and therefore its objective value is at least as good as the previous one. Hence, $(x,\bar z,\bar \theta)$ is also optimal. In addition, values $\bar{z}$ satisfy by construction that  $\sum_{i\in P}\bar{z}_{ij'}=\sum_{i\neq i'} z_{ij'}+\bar{z}_{i'j'}=0+1=1$. Therefore \eqref{owab} holds. \\

%\JP{Hence, by constraints \eqref{owad''} it holds that $C^{i'}x\le\theta_j$, for all $j\in P$, which means that, among all cost functions, cost function $i'$ has the the smallest value. Thus, it could be placed at position $p$ or at some previous position $k$ provided that any other cost function with equal value as $C^{i'}x$ occupies the position $p$ (there are ties in the smallest cost function value). In any case, we have that $\sigma_k=i'$ and thus  $\sum_{j\in P}z_{ij}=1$ for all $i\neq i'$ and all position values are well defined by using Property \ref{lemaFz_R1}.}
\hfill$\square$\\

%---------------------------------
%------------------------------------------------------------------

We can now relate the domains of the formulations considered so far.

\begin{proposition} The following relationships hold  \label{propertyowa15JP}
$$\Omega_0^z\equiv \Omega^z \subsetneq \Omega_{R1}^z\subsetneq \Omega_{R2}^z \subsetneq  \Omega_{R3}^z$$
\end{proposition}
\texttt{Proof}.
\begin{itemize}
\item $\Omega^z \subsetneq \Omega_{R1}^z$: Every feasible solution in $\Omega^z$ verifies inequalities of $\Omega_{R1}^z$. However, a feasible solution in $\Omega_{R1}^z$ with $\theta_j\leq \theta_{j+1}$ for some $j\in P$ is not feasible in $\Omega^z$.
\item $\Omega_{R1}^z\subsetneq \Omega_{R2}^z$: Every feasible solution in $\Omega_{R1}^z$ verifies the inequalities of $\Omega_{R2}^z$. However, a feasible solution in $\Omega_{R2}^z$ where for some $j\in P$, $z_{ij}=1$, for all $i\in P$ is not feasible in $\Omega_{R1}^z$.
\item $\Omega_{R2}^z\subsetneq \Omega_{R3}^z$: Every feasible solution in $\Omega^z$ verifies the inequalities of $\Omega_{R3}^z$. However, a feasible solution in $\Omega_{R3}^z$ with $z_{ij}=0, \, i,j\in P$ is not feasible in $\Omega_{R2}^z$.
\end{itemize}
\hfill$\square$\\
%---------------------------------

%------------------------------------------------------------------

\begin{proposition}
The dimension of {$\Omega_0^z$} is $p^2-p+1+dim(Q)$.
\end{proposition}
\texttt{Proof}.\newline

Suppose $Q\subseteq  \mathbb{R}^n$. Then, $\Omega_0^z$ is embedded in a space of dimension $p^2 + p +n$. Furthermore, since there are $2p-1$ linearly independent equations in \eqref{owab} and \eqref{owac} and the dimension of $Q$ does not depend on relations \eqref{owab}-\eqref{owae}, then the dimension of \eqref{owab}-\eqref{owaf} is at most $p^2-p+1+dim(Q)$. Denote by $q=dim (Q)$ and by $\rho=p^2-2p+1$. Next, we show that there exist $q+\rho+p+1$ ( equal to $p^2-p+2+dim(Q)$) affinely independent points in $\Omega_0^z$ and consequently, the dimension of $\Omega_0^z$ is $p^2-p+1+dim(Q)$.\\

Let  $v=(v_j)_{j\in P}$ where $v_j=M+p-j+1$ for $M>0$ and sufficiently large.
Denoting by $\mathbf{e^j}\in \mathbb{R}^p$ the $j$-th vector of the canonical basis in $\mathds{R}^p$ and $0<\varepsilon<1$, let $\theta^j=\{v+\varepsilon \mathbf{e^j},j\in P\}$. Moreover, let {$\theta^{p+1}=(M,\ldots,M)'$.}  We observe that the vectors $\theta^j$, $j=1,\ldots,p+1$  are affinely independent and each one of them satisfies  inequalities (\ref{owae}).

Next, since $dim(Q)=q$, we take $q+1$ arbitrary affinely independent vectors {$x^i\in Q$, $i=1,\ldots,q+1$.} Furthermore, let $z^k\in \{0,1\}^{p^2}$ $k=1,\ldots,\rho+1$, be $\rho+1$ affinely independent vectors satisfying (\ref{owab}) and (\ref{owac}). Note that the latter is always possible since there are $p^2$ degrees of freedom for the coordinates of $z$ variables and only $2p$ equations  being one of them  linearly dependent of the others.\\

Now,  we prove that any point of the form {$((x^{i})^\prime,(z^{k})^\prime,(\theta^{l})^\prime)^\prime$} $i=1,\ldots,q+1$, $k=1,\ldots,\rho+1,\; l=1,\ldots,p+1$ satisfies \eqref{owab}-\eqref{owae}. Indeed, by construction the first block of coordinates defines a point in $Q$, the second block satisfies  (\ref{owab}) and (\ref{owac}) and the third one \eqref{owae}. Thus, it remains to prove that such a generic point  also satisfies \eqref{owad} as follows:
$$ C^ix^{i}\le M \le M+p-j+1 \le \theta^l_j + M(1-z^{k}_{ij}), \quad \forall\; i,j.$$

Consider the $q+\rho+p$ points defined as the column vectors of the  matrix $A=(A^1|A^2|A^3)$ where
$$ A^1=\left(\begin{array}{cccc} x^1 & x^2 & \ldots & x^q \\ z^2 & z^1 & \ldots & z^1 \\ \theta^2 & \theta^1 & \ldots & \theta^1 \end{array} \right), \; A^2=\left(\begin{array}{cccc} x^1 & x^1 & \ldots & x^1 \\ z^1 & z^2 & \ldots & z^{\rho} \\ \theta^2 & \theta^1 & \ldots & \theta^1 \end{array} \right), \; A^3=\left(\begin{array}{ccccc} x^1 & x^1 & x^1& \ldots & x^1 \\ z^1 & z^3 & z^1 &  \ldots & z^1 \\ \theta^1 & \theta^2 & \theta^3 &  \ldots & \theta^p \end{array} \right). $$

By construction, each submatrix $A^i$ has its column vectors linearly independent from one  another since the $i$-th block is formed by linearly independent vectors. Next, clearly each column vector of $A^1$ is linearly independent from those of $A^2$ and $A^3$ and each column vector of $A^2$ is linearly independent from those of $A^3$. Therefore, the rank of $A$ is $q+\rho+p=q+p^2-p+1$.\\

Finally, the column vectors of $A$  are linearly independent and feasible points of  \eqref{owab}-\eqref{owae}. In addition, we can easily construct another feasible point, different from those considered previously and affinely independent from all of them, namely $((x^{q+1})^\prime,(z^{\rho+1})^\prime,(\theta^{p+1})^\prime)'$ .  Hence the dimension of $\Omega^z$ is $q+\rho+p=q+p^2-p+1$.

\hfill$\square$

%------------------------------------------------------------------

\begin{proposition}
The following inequalities define facets in {$\Omega_0^z$:}
\formulaN{C^i x\leq\theta_p+M(1-z_{ip})                            &&  i\in P \label{faceta1}}
\formulaN{\theta_j \geq \theta_{j+1}                                 &&  j\in P:j<p\label{faceta2}}
\end{proposition}
\texttt{Proof}.\newline
\texttt{\eqref{faceta1} is a facet defining inequality:}\\
We prove that for each $i'\in P$ there exist {$dim(\Omega_0^z)=p^2-p+dim(Q)+1$} affinely   independent points of {$\Omega_0^z$}  that verify $C^{i'} x=\theta_p+M(1-z_{i'p})$.\\

As in the proof of the above proposition, we take $q+1$ arbitrary affinely  independent points $x^i$, $i=1,\ldots,q+1$ in $Q$. Furthermore, let $z^k\in \{0,1\}^{p^2}$ $k=1,\ldots,\rho$, be $\rho$ affinely independent points (recall that $\rho:=p^2-2p+1$) satisfying (\ref{owab}), (\ref{owac}) and $z_{i'p}=1$. Note that the latter is always possible since there are $p^2$ degrees of freedom for the coordinates of $z$ variables and  $2p$ non redundant equations ($2p-1$ as in the case above and $z_{i'p}=1$).

Let  $v^l=(v_j^l)_{j\in P}$ where $v_j^l=C^{i'}x^l+M+p-j$ if $j<p$ and $v_p^l=C^{i'}x^l$ for $M>0$ and sufficiently large.
Denoting by $\mathbf{e^j}\in \mathbb{R}^p$ the $j$-th vector of the canonical basis in $\mathds{R}^p$ and $0<\varepsilon<1$, let $\bar \theta^{lj}=\{v^l+\varepsilon \mathbf{e^j},j\in P\}$ if $j<p$ and $\bar \theta^{lp}=v^l$,  $\bar \theta^{l,p+1}=(C^{i'}x^l+M,\ldots,C^{i'}x^l+M,C^{i'}x^l)'$. We observe that for each $l$ fixed, the vectors $\bar \theta^{lj}$ $j=1,\dots,p+1$ are affinely independent and each one of them satisfies  inequalities (\ref{owae}).

Now,  we prove that any point of the form $((x^{l})^\prime,(z^{k})^\prime,(\theta^{lj})^\prime)^\prime$  $k=1,\ldots,\rho,\; j=1,\ldots,p+1$ satisfies \eqref{owab}-\eqref{owae} and $z^k_{i'p}=1$. Indeed, by construction the first block of coordinates defines a point in $Q$, the second block satisfies  (\ref{owab}), (\ref{owac})  and $z_{i'p}=1$, and the third one \eqref{owae}. Thus, it remains to prove that such a generic point  also satisfies \eqref{owad}. We distinguish two cases:
\begin{itemize}
\item If $j<p$ then
$$ C^ix^l\le C^{i'}x^l+M+p-j+1+M=C^{i'}x^l+M+p-j+M(1-z^k_{ij})=\bar \theta^{lj}+M(1-z^k_{ij}), \quad \forall i.$$
\item If $j=p$ we have that
\begin{eqnarray*}
 C^ix^{l}\le C^{i'}x^l+M=C^{i'}x^l+M(1-z^k_{ip}), & \forall\; i\neq i',\\
C^{i'}x^{l}\le C^{i'}x^l=C^{i'}x^l+M(1-z^k_{i'p}), &  \mbox{otherwise. (Recall that } z^k_{i'p}=1.)
\end{eqnarray*}

\end{itemize}

Consider the $q+\rho-1+p$ points defined as the column vectors of the  matrix $\bar A=(\bar A^1|\bar A^2|\bar A^3)$ where
$$\bar  A^1=\left(\begin{array}{cccc} x^1 & x^2 & \ldots & x^q \\ z^2 & z^1 & \ldots & z^1 \\ \bar \theta^{11} &\bar  \theta^{21} & \ldots & \bar \theta^{q1} \end{array} \right), \; \bar A^2=\left(\begin{array}{cccc} x^1 & x^1 & \ldots & x^1 \\ z^1 & z^2 & \ldots & z^{\rho-1} \\ \bar \theta^{12} & \bar \theta^{11} & \ldots & \bar \theta^{11} \end{array} \right), \; \bar A^3=\left(\begin{array}{ccccc} x^1 & x^1 & x^1& \ldots & x^1 \\ z^1 & z^3 & z^1 &  \ldots & z^1 \\ \bar \theta^{11} & \bar \theta^{12} & \bar \theta^{13} &  \ldots & \bar \theta^{1p} \end{array} \right). $$

By construction, each submatrix $\bar A^i$ has its column vectors linearly independent from one  another since the $i$-th block is formed by linearly independent vectors.   Next, clearly each column vector of $\bar A^1$ is linearly independent from those of $\bar A^2$ and $\bar A^3$ and each column vector of $\bar A^2$ is linearly independent from those of $\bar A^3$. Therefore, the rank of $A$ is $q+\rho-1+p=q+p^2-p$.

Finally, the column vectors of $A$  together with the point $((x^{q+1})^\prime,(z^{\rho+1})^\prime,(\theta^{q+1,j})^\prime)'$ are feasible points of  \eqref{owab}-\eqref{owae} that satisfy $C^{i'} x=\theta_p+M(1-z_{i'p})$; and this last vector is clearly affinely independent from the those in $\bar A$, therefore  (\ref{faceta1}) is a facet defining inequality for $\Omega^z$.\\

\texttt{\eqref{faceta2} is a facet defining inequality:}\\
In order to prove that for each $j'\in P\setminus\{p\}$ there exist {$dim(\Omega_0^z)=p^2-p+dim(Q)+1$} affinely independent points of {$\Omega_0^z$}  that verify $\theta_{j'} = \theta_{j'+1}$, we can proceed analogously as before
considering $v=(v_j)_{j=1}^p$, where $v_j=M+p-j+1$ if $j\neq j'+1$ and {$v_{j'+1}=M+p-j'+2$} and the points $\hat \theta^{j}=\{v+\varepsilon \mathbf{(e^j+e^{j'+1})},j\in P\setminus\{p\}\}$. In addition, we take $\hat \theta^{p}=(M,\ldots,M)'$. We observe that  the vectors $\hat \theta^{j}$ $j=1,\dots,p$ are affinely independent and each one of them satisfies  $\hat \theta^j_{j'}=\hat \theta^{j}_{j'+1}$.

Any point of the form {$((x^{i})^\prime,(z^{k})^\prime,(\hat \theta^{l})^\prime)^\prime$ $i=1,\ldots,q+1$,} $k=1,\ldots,\rho+1,\; l=1,\ldots,p$ satisfies \eqref{owab}-\eqref{owae} and $\hat \theta^l_{j'}=\hat \theta^{l}_{j'+1}$.

Consider the $q+\rho+p-1$ points defined as the column vectors of the  matrix $\hat A=(\hat A^1|\hat A^2|\hat A^3)$ where
$$\hat  A^1=\left(\begin{array}{cccc} x^1 & x^2 & \ldots & x^q \\ z^2 & z^1 & \ldots & z^1 \\ \hat \theta^2 & \hat \theta^1 & \ldots & \hat \theta^1 \end{array} \right), \; \hat A^2=\left(\begin{array}{cccc} x^1 & x^1 & \ldots & x^1 \\ z^1 & z^2 & \ldots & z^{\rho} \\ \hat \theta^2 & \hat \theta^1 & \ldots & \hat \theta^1 \end{array} \right), \;
\hat A^3=\left(\begin{array}{ccccc} x^1 & x^1 & x^1& \ldots & x^1 \\ z^1 & z^3 & z^1 &  \ldots & z^1 \\ \hat \theta^1 & \hat \theta^2 & \hat \theta^3 &  \ldots & \hat \theta^{p-1} \end{array} \right). $$

By construction, each submatrix $\hat A^i$ has its column vectors linearly independent from one  another since the $i$-th block is formed by linearly independent vectors.  Next, clearly each column vector of $\hat  A^1$ is linearly independent from those of $\hat  A^2$ and $\hat A^3$ and each column vector of $\hat  A^2$ is linearly independent from those of $\hat  A^3$. Therefore, the rank of $\hat A$ is $q+\rho+p-1=q+p^2-p$.

Finally, the column vectors of $\hat A$  are linearly independent and are also feasible points of  \eqref{owab}-\eqref{owae} that satisfy $\theta_{j'}=\theta_{j'+1}$. Next, we can easily add a new  feasible point, for instance $((x^{q+1})^\prime,(z^{\rho+1})^\prime,(\hat \theta^{p})^\prime)'$ that also satisfies $\theta_{j'}=\theta_{j'+1}$ and that is clearly affinely independent from the those in $\hat A$. Hence,  \eqref{faceta2} is a facet defining inequality for $\Omega^z$.

\hfill$\square$

The following table summarizes the previous proposed formulations. Formulas included on each formulation have been checked (\checkmark) whereas those not appearing are marked with a dot (.).

\renewcommand{\esp}{@{\hspace{0.2cm}}}
\begin{small}
\begin{table}[h!]
%\begin{sidewaystable}
\begin{center}
\begin{tabular}{|@{}\esp l  | c \esp c \esp c\esp c \esp c \esp @{}|}
\hline
 & $F^z_0$ & $F^z$ & $F^z_{R1}$ & $F^z_{R2}$ & $F^z_{R3}$  \\
\hline
$\displaystyle{\min\sum_{j\in P}\omega_j\theta_j }$             & \checkmark & \checkmark & \checkmark & \checkmark & \checkmark  \\
\hline
$\displaystyle{\sum_{i\in P}z_{ij}=1, \, j\in P }$              & \checkmark & \checkmark & \checkmark & \checkmark & . \\
$\displaystyle{\sum_{j\in P}z_{ij}=1, \, i\in P }$              & \checkmark & \checkmark & \checkmark & . & .  \\
$\displaystyle{\sum_{i\in P}z_{ij}\leq1, \, j\in P }$           & . & . & . & . & \checkmark \\
\hline
$\displaystyle{C^i x\leq\theta_j+M(1-z_{ij}),\,i,j\in P }$                              & \checkmark & . & . & . & . \\
$\displaystyle{C^ix\leq \theta_j+M(1-\sum_{k\geq j}z_{ik}),\,i,j\in P }$              & . & \checkmark & \checkmark & \checkmark & . \\
$\displaystyle{C^ix\leq \theta_j+M\sum_{k< j}z_{ik},\,i,j\in P }$                     & . & . & . & . & \checkmark   \\
\hline
$\displaystyle{\theta_j \geq \theta_{j+1},\,j\in P:j<p}$                    & \checkmark & \checkmark & . & . & . \\
\hline
$\displaystyle{x\in Q, z\in\{0,1\}^{p\times p}     }$                    & \checkmark & \checkmark & \checkmark & \checkmark& \checkmark \\
\hline
\end{tabular}
\end{center}
\caption[]{Summary of the proposed formulations for the OWAP.}
\label{label}
%\end{sidewaystable}
\end{table}
\end{small}

%------------------------------------------------------------------
%\newpage
%-----------------------------------------------------------------------------------------------------------------------
%\subsection{Formulation OWAP2}
\subsection{OWAP formulations with variables for the values of cost functions occupying specific sorted positions}\label{subsec:2}

%\subsection{Formulations with variables defining the values of the cost functions when occupying specific positions in the ordering}\label{subsec:2}
Another OWAP formulation can be obtained by defining an additional set of continuous variables $y=(y_{ij})_{i,j\in P}\in\mathbb{R}^{p\times p}$, where $y_{ij}$ denotes the value of cost function $i$ if it occupies the $j$-th position in the ordering. The formulation is as follows:

%\newcounter{counterequation}
%\addtocounter{counterequation}{1}
%\addtocounter{counterequation}{\value{equation}}
%\begin{subequations}
%\vspace{-0.5cm}
%\begin{align}
%&\Hp \hspace{8cm}&&\notag\\
%&\hhp F^z_0: \hspace{0.5cm} V=\min\sum_{j\in P}\omega_j\theta_j  &&                 \label{owaa}\\ s.t.
%&\Hp \sum_{i\in P}z_{ij}=1                                      &&  j\in P          \label{owab}\\
%&\Hp \sum_{j\in P}z_{ij}=1                                      &&  i\in P          \label{owac}\\
%&\Hp C^i x\leq\theta_j+{M}(1-z_{ij})                            &&  i,j\in P         \tag{\arabic{counterequation}d$^0$}\label{owad0}\\\addtocounter{equation}{1}
%&\Hp \theta_j \geq \theta_{j+1}                                 &&  j\in P:j<p      \label{owae}\\
%&\Hp x\in Q, z\in\{0,1\}^{p\times p}                            &&                  \label{owaf}
%\end{align}
%\end{subequations}

\setcounter{counterequation}{\value{equation}}
\addtocounter{counterequation}{1}
\begin{subequations}
\vspace{-0.5cm}
\begin{align}
&\Hp \hspace{8cm}&&\notag\\
&\hhp F^{zy}_0: \hspace{0.5cm} V=\min\sum_{j\in P}\omega_j\sum_{i\in P}y_{ij}   &&                        \label{owa2a}\\ s.t.
&\Hp \sum_{i\in P}z_{ij}=1                                                      &&  j\in P                \label{owa2b}\\
&\Hp \sum_{j\in P}z_{ij}=1                                                      &&  i\in P                \label{owa2c}\\
&\Hp C^ix\leq\sum_{i'\in P}y_{i'j}+M(1-z_{ij})                                  &&  i,j\in P        \tag{\arabic{counterequation}d$^0$}\label{owa2d0}\\\addtocounter{equation}{1}
&\Hp \sum_{i\in P}y_{ij} \geq\sum_{i\in P}y_{ij+1}                              &&  j\in P:j<p            \label{owa2e}\\
&\Hp x\in Q, z\in\{0,1\}^{p\times p}                                            &&                        \label{owa2f}
\end{align}

%\begin{subequations}
%\vspace{-0.5cm}
%\begin{align}
%&\Hp \hspace{8cm}&&\notag\\
%&\hhp F^{zy}_0: \hspace{0.5cm} V=\min\sum_{j\in P}\omega_j\sum_{i\in P}y_{ij}   &&                  \tag{\rm 2a}        \label{owa2a}\\ s.t.
%&\Hp \sum_{i\in P}z_{ij}=1                                                      &&  j\in P          \tag{\rm 2b}        \label{owa2b}\\
%&\Hp \sum_{j\in P}z_{ij}=1                                                      &&  i\in P          \tag{\rm 2c}        \label{owa2c}\\
%&\Hp C^ix\leq\sum_{i'\in P}y_{i'j}+M(1-z_{ij})                                  &&  i,j\in P        \tag{\rm 2d}        \label{owa2d}\\
%&\Hp \sum_{i\in P}y_{ij} \geq\sum_{i\in P}y_{ij+1}                              &&  j\in P:j<p      \tag{\rm 2e}        \label{owa2e}\\
%&\Hp x\in Q, z\in\{0,1\}^{p\times p}                                            &&                  \tag{\rm 2f}        \label{owa2f}
%\end{align}
%\end{subequations}

Next we study some properties of formulation $F^{zy}_0$ and relate it to the OWAP formulations presented above. Denote by $\Omega^{zy}_0$ the domain of Problem $F_0^{zy}$.  Consider first, for any $M>0$ sufficiently large, the following set of inequalities

\formulaN{y_{ij}\le M z_{ij}, \qquad i, j\in P. && \label{owa2g}}
\end{subequations}

\begin{property}\label{dominance}
There is an optimal solution to  $F^{zy}_0$ for which  constraints \eqref{owa2g} hold.\\
\end{property}

\texttt{Proof}.\newline
Observe that constraints  \eqref{owa2d0} imply that $\sum_{k\in P}y_{kj}\ge C^ix $ for all $i,j\in P$ with $z_{ij}=1$. Since constraints \eqref{owa2b} indicate that for $j\in P$ fixed there exists a unique index, say $i(j)$ with $z_{i(j),j}=1$,
the above condition reduces to  $\sum_{k\in P}y_{kj}\ge C^{i(j)}x $, for all $j\in P$. Because of the non-negativity of the cost coefficients, we can thus deduce that an optimal solution exists to  $F^{zy}_0$ in which
\begin{equation}\sum_{k\in P}y_{kj}= C^{i(j)}x , \text{ for all } j\in P\label{condition}.\end{equation}

Let now $(x,y,z)\in \Omega^{zy}_0$ be such an optimal solution, and suppose it violates some constraint  \eqref{owa2g}. That is, there exist $i', j'\in P$ with $y_{i'j'}> M z_{i'j'}$.
Hence, $\sum_{i\in P}y_{ij'}>M z_{i'j'}$, contradicting \eqref{condition} unless $z_{i'j'}=0$. In other words, $i(j')\ne j'$.

Consider now the solution $(x,\overline{y},z)$, with the same $x$ and $z$ values as before where $\overline{y}$ is defined as follows:

\[\overline{y}_{ij}=\begin{cases}
0   &   \text{if } i=i', \text{ and } j=j'\\
y_{i(j'),j'}+y_{i'j'}   &   \text{if } i=i', \text{ and } j=i(j')\\
y_{ij} &        \text{otherwise.}
\end{cases}\]

Indeed  $(x,\overline{y},z)\in \Omega^{zy}_0$, as it is immediate to check that it  satisfies constraints  \eqref{owa2b}-\eqref{owa2f}. Furthermore, by construction, it satisfies the constraint \eqref{owa2g} associated with $i', j'$.
Finally, note that it is optimal to $F^{zy}_0$, since $\sum_{i\in P}\overline{y}_{ij}=\sum_{i\in P}y_{ij}$, for all $j\in P$. \hfill$\square$\\

Note that if there is $j\in P$ with $\omega_j=0$ then it is possible to have optimal solutions to $F^{zy}_0$  that do not satisfy constraints \eqref{owa2g}. However, because of Property \ref{dominance}, constraints \eqref{owa2g} can be useful to restrict the domain where optimal solutions are sought. Let $$\Omega^{GS'}=\{(x,y,z,\theta) \text{ satisfying constraints } \text{\eqref{owa2b}, \eqref{owa2c}, \eqref{owa2d0}, \eqref{owa2e}, \eqref{owa2f}, \eqref{owa2g}}\}.$$
Then, a different  formulation that also ensures to obtain an optimal solutions to $F^{zy}_0$ is:

\begin{subequations}
\vspace{-0.5cm}
\begin{align}
&\Hp \hspace{8cm}&&\notag\\
&\hhp F^{GS'} \hspace{0.5cm} V=\min\sum_{j\in P}\omega_j\theta_j  && \notag \\ s.t.
&\Hp (x,y,z,\theta)\in \Omega^{GS'}. \notag &&
\end{align}
\end{subequations}

Formulation $F^{GS'}$ is closely related to the formulation used in \citet{Galand2012} for modeling the  minimum cost spanning tree OWAP. In their formulation they operate on a domain which is like $\Omega^{GS'}$ except that constraints \eqref{owa2d0} have been substituted by constraints

\formula{\sum_{j\in P}y_{ij}=C^ix \qquad i\in P. \tag{\arabic{counterequation}h} \label{owa2h}}

Let $\Omega^{GS}=\{(x,y,z,\theta) \text{ satisfying constraints \eqref{owa2b}, \eqref{owa2c}, \eqref{owa2e}, \eqref{owa2f}, \eqref{owa2g}, \eqref{owa2h}}\}$, denote the domain  used in \cite{Galand2012}. Then, it is straightforward to conclude the following.
\begin{property} \label{F2:Gallant}
The domains $\Omega^{GS}$ and $\Omega^{GS'}$ satisfy $\Omega^{GS}\subseteq \Omega^{GS'}$. Moreover, if $(x^*,y^*,z^*)$ is an optimal solution of $F^{GS'}$ then it is also optimal for $F^{GS}$ and conversely.
\end{property}

%------------------------------------------------------------------

We can also relate $F^{zy}_0$ with $F^{z}_0$ and its variations. In particular, because of the relationship \formulaN{\theta_{j}=\sum_{i \in P}y_{ij}, \, j\in P. \label{relOWAP1OWAP2}} we have:

\begin{property}\label{proprelOWAP1OWAP2}
For each  optimal solution to $F^{zy}_0$, $(x^*,z^*,\theta^*)$, there exists $(x^*,y^*,z^*,\theta^*)$ optimal solution for $F^{z}_0$ and conversely. Moreover, $\sum_{j \in P} w_j \sum_{i\in P} y_{ij}^*=\sum_{j\in P} w_j \theta_j^*$.
\end{property}

%---------------------------------

By above result,  we can derive variations of $F^{zy}$ similar to the ones obtained for $F^z$ with similar properties. These constructions are straightforward and therefore are left for the interested readers.

%------------------------------------------------------------------

The following table summarizes the proposed formulations of this subsection that can be derived from those of Subsection \ref{subsec:1}.

\renewcommand{\esp}{@{\hspace{0.2cm}}}
\begin{small}
\begin{table}[h!]
%\begin{sidewaystable}
\begin{center}
\begin{tabular}{|@{}\esp l  | c \esp c \esp c\esp c \esp c \esp @{}|}
\hline
 & $F^{zy}_0$ & $F^{zy}$ & $F^{zy}_{R1}$ & $F^{zy}_{R2}$ & $F^{zy}_{R3}$  \\
\hline
$\displaystyle{\min\sum_{j\in P}\omega_j\theta_j }$             & \checkmark & \checkmark & \checkmark & \checkmark & \checkmark  \\
\hline
$\displaystyle{\sum_{i\in P}z_{ij}=1, \, j\in P }$              & \checkmark & \checkmark & \checkmark & \checkmark & . \\
$\displaystyle{\sum_{j\in P}z_{ij}=1, \, i\in P }$              & \checkmark & \checkmark & \checkmark & . & .  \\
$\displaystyle{\sum_{i\in P}z_{ij}\leq1, \, j\in P }$           & . & . & . & . & \checkmark \\
\hline
$\displaystyle{C^i x\leq\sum_{i'\in P}y_{i'j}+M(1-z_{ij}),\,i,j\in P }$                               & \checkmark & . & . & . & . \\
$\displaystyle{C^ix\leq \sum_{i'\in P}y_{i'j}+M(1-\sum_{k\geq j}z_{ik}),\,i,j\in P }$              & . & \checkmark & \checkmark & \checkmark & . \\
$\displaystyle{C^ix\leq \sum_{i'\in P}y_{i'j}+M\sum_{k< j}z_{ik},\,i,j\in P }$                    & . & . & . & . & \checkmark   \\
\hline
$\displaystyle{\sum_{i\in P}y_{ij} \geq\sum_{i\in P}y_{ij+1},\,j\in P:j<p}$           & \checkmark & \checkmark & . & . & . \\
\hline
$\displaystyle{x\in Q, z\in\{0,1\}^{p\times p}     }$                    & \checkmark & \checkmark & \checkmark & \checkmark& \checkmark \\
\hline
\end{tabular}
\end{center}
\caption[]{Summary of the proposed formulations for the OWAP.}
\label{label}
%\end{sidewaystable}
\end{table}
\end{small}

%=======================================================================================================================
\subsection{Using variables defining relative positions of sorted cost function values}\label{subsec:3}

We close this section with another formulation which uses decision variables defining the relative positions of the sorted cost function values. As we have seen in Section \ref{rel:z-s} it is possible to describe permutations with variables representing the relative positions of the sorted values. Next we use such variables to obtain formulations for the OWAP.

For $i,j\in P$, consider the binary variable $s_{ij},\; i,j\in P$ as
\[
s_{ij}=\begin{cases}
 1 & \text{if cost function $i$ is placed after position $j$ in the ordering,}\\
 0 & \text{otherwise.}
\end{cases}
\]

As we have seen in Section \ref{rel:z-s}, for all $i, j\in P$, $s_{ij}=1-\sum_{k\geq j}z_{ik},\, i,j\in P $. Therefore, variables $z$ and $s$ are related by means of

\formulaN{z_{ij}=\begin{cases}
 s_{ij+1}-s_{ij} &  i\in P, j=1,...,p-1 \\
 1- s_{ij} &  i\in P, j=p
\end{cases}}

Thus, we can reformulate the OWAP in the new space of the $s$ variables as

\begin{subequations}
\vspace{-0.5cm}
\begin{align}
&\Hp \hspace{8cm}&&\notag\\
&\hhp F^s: \hspace{0.5cm} V=\min\sum_{j\in P}\omega_j\theta_j     &&                    \label{owa3a}\\ s.t.
&\Hp \sum_{i\in P}s_{ij}= j-1                                       &&  j\in P          \label{owa3b}\\
&\Hp s_{ij+1}-s_{ij}\geq 0                                          &&  i,j\in P:j<p    \label{owa3c}\\
&\Hp C^i x\leq\theta_j+M s_{ij}                                     &&  i,j\in P        \label{owa3d}\\
&\Hp \theta_j \geq \theta_{j+1}                                     &&  j\in P:j<p      \label{owa3e}\\
&\Hp x\in Q, s\in\{0,1\}^{p\times p}                                &&                  \label{owa3f}
\end{align}
\end{subequations}

%\begin{subequations}
%\vspace{-0.5cm}
%\begin{align}
%&\Hp \hspace{8cm}&&\notag\\
%&\hhp F^s: \hspace{0.5cm} V=\min\sum_{j\in P}\omega_j\theta_j     &&                  \tag{\rm 3a} \label{owa3a}\\ s.t.
%&\Hp \sum_{i\in P}s_{ij}= j-1                                       &&  j\in P        \tag{\rm 3b}  \label{owa3b}\\
%&\Hp s_{ij+1}-s_{ij}\geq 0                                          &&  i,j\in P:j<p  \tag{\rm 3c}  \label{owa3c}\\
%&\Hp C^i x\leq\theta_j+M s_{ij}                                     &&  i,j\in P      \tag{\rm 3d}  \label{owa3d}\\
%&\Hp \theta_j \geq \theta_{j+1}                                     &&  j\in P:j<p    \tag{\rm 3e}  \label{owa3e}\\
%&\Hp x\in Q                                                         &&                \tag{\rm 3f}  \label{owa3f}
%\end{align}
%\end{subequations}

Since $F^s$ is obtained from $F^z$ by a change of variable and there is a one to one correspondence between feasible solutions, we can state the following result. Let $\Omega^s$ be the feasible region of Problem $F^s$.

%------------------------------------------------------------------
\begin{property}\label{proprelOWAP1OWAP3}
For each solution $(x,s,\theta)\in \Omega^s$ there exists $(x,z,\theta)\in \Omega^z$ with equal objective value and conversely.
\end{property}
%---------------------------------
%------------------------------------------------------------------

By analogy with the notation used in Section \ref{subsec:1} let us define the following domains and problems related to $F^s$:

\begin{subequations}
\vspace{-0.5cm}
\begin{align}
&\Hp \hspace{8cm}&&\notag\\
&\hhp F^{s}_{R1} \hspace{0.5cm} V=\min\sum_{j\in P}\omega_j\theta_j  && \notag \\ s.t.
&\Hp  (x, s,\theta)\in \Omega^{s}_{R1}. \notag &&
\end{align}
\end{subequations}

\noindent with $\Omega^{s}_{R1}=\{(x,s,\theta) \text{ satisfying constraints } \text{\eqref{owa3b}, \eqref{owa3c}, \eqref{owa3d}, \eqref{owa3f}}\}$.

\begin{subequations}
\vspace{-0.5cm}
\begin{align}
&\Hp \hspace{8cm}&&\notag\\
&\hhp F^{s}_{R2} \hspace{0.5cm} V=\min\sum_{j\in P}\omega_j\theta_j  && \notag \\ s.t.
&\Hp  (x,s,\theta)\in \Omega^{s}_{R2}. \notag &&
\end{align}

\end{subequations}
\noindent with $\Omega^{s}_{R2}=\{(x,s,\theta) \text{ satisfying constraints } \text{\eqref{owa3b}, \eqref{owa3d}, \eqref{owa3f}}\}$.

\begin{subequations}
\vspace{-0.5cm}
\begin{align}
&\Hp \hspace{8cm}&&\notag\\
&\hhp F^{z}_{R3} \hspace{0.5cm} V=\min\sum_{j\in P}\omega_j\theta_j  && \notag \\ s.t.
&\Hp  (x,z,\theta)\in \Omega^{z}_{R3}. \notag &&
\end{align}
\end{subequations}

\noindent with
$\Omega^{s}_{R3}=\{(x,s,\theta) \text{ satisfying constraints } \text{\eqref{owa3b<}, \eqref{owa3d}, \eqref{owa3f}}\}$, where \eqref{owa3b<} are the inequality version of constraints \eqref{owa3b}. That is,
\formula{\sum_{i\in P}s_{ij}\le j-1    \qquad  j\in P.        \tag{\ref{owa3b}$_\le$}  \label{owa3b<}}

\begin{property} The following relationships hold.
\begin{enumerate}
  \item Every optimal solution to $F^{s}_{R1}$ is optimal to $F^{s}$ and conversely.
  \item Every optimal solution to $F^{s}_{R2}$ is optimal to $F^{s}$ and conversely.
  \item Every optimal solution to $F^{s}_{R3}$ is optimal to $F^{s}$ and conversely.
  \item $\Omega^s\subsetneq \Omega^s_{R1}\subsetneq \Omega^s_{R2} \subsetneq  \Omega^s_{R3}$.
\end{enumerate}

\end{property}

\texttt{Proof}.\newline The proofs of the above statements follow directly from the relationship that links variables $z$ and $s$, namely \eqref{for:s-z-1} and \eqref{for:s-z-2}. Specifically, statement 1 follows from Property \ref{propertyowa12}, statement 2 from Property \ref{propertyowa13}, statement 3 from Property \ref{propertyowa14JP} and statement 4 from Property \ref{propertyowa15JP}.

%---------------------------------

\hfill$\square$\\
%---------------------------------
%%------------------------------------------------------------------

%=======================================================================================================================
%\newpage
%=======================================================================================================================

\subsection{Formulations summary}
The following table summarizes the previous proposed formulations.

\renewcommand{\esp}{@{\hspace{0.2cm}}}
\begin{small}
\begin{table}[h!]
%\begin{sidewaystable}
\begin{center}
\begin{tabular}{|@{}\esp l  | c \esp c \esp c\esp c \esp @{}|}
\hline
 & $F^{s}$ & $F^{s}_{R1}$ & $F^{s}_{R2}$ & $F^{s}_{R3}$  \\
\hline
$\displaystyle{\min\sum_{j\in P}\omega_j\theta_j }$             & \checkmark & \checkmark & \checkmark & \checkmark   \\
\hline
$\displaystyle{\sum_{i\in P}z_{ij}=1, \, j\in P }$              & \checkmark  & \checkmark & \checkmark & . \\
$\displaystyle{\sum_{j\in P}z_{ij}=1, \, i\in P }$              & \checkmark  & \checkmark & . & .  \\
$\displaystyle{\sum_{i\in P}z_{ij}\leq1, \, j\in P }$           & . & . & . & \checkmark \\
\hline
$\displaystyle{C^i x\leq\theta_j+M s_{ij},\,i,j\in P}$                              & \checkmark & \checkmark & \checkmark & \checkmark \\
\hline
$\displaystyle{\theta_j \geq \theta_{j+1},\,j\in P:j<p}$                    & \checkmark & . & . & . \\
\hline
$\displaystyle{x\in Q, z\in\{0,1\}^{p\times p}     }$                    & \checkmark & \checkmark & \checkmark & \checkmark \\
\hline
\end{tabular}
\end{center}
\caption[]{Summary of the proposed formulations for the OWAP.}
\label{label}
%\end{sidewaystable}
\end{table}
\end{small}

\section{Valid inequalities and reinforcements for the OWAP formulation}\label{Sec:ValIneq}
\subsection{Valid inequalities for the (OWAP) formulation}\label{SecOwa1Vi}
In this section we derive different valid inequalities  for all the formulations presented in previous sections. For the sake of simplicity, we present all inequalities for the formulations developed in Subsection \ref{subsec:1}. However, all these inequalities can be easily adapted to the remaining formulations just by means of the substitutions explained by Equations \eqref{for:s-z-2} and \eqref{relOWAP1OWAP2}.

\begin{itemize}
\item %%\TP{cota\_inf\_i, cota\_sup\_i}
\emph{Constraints related to bounds of cost function values.} Let $l_i$ ($u_i$) denote the minimum (maximum) objective value relative to cost function $i\in P$, respectively. It is clear that $l_i$ ($u_i$) are valid lower (upper) bounds on the value of objective $i$, independently of the position of cost function $i$ in the ordering. Therefore we obtain the following two sets of constraints which are valid for the OWAP:
\formulaN{l_{i} \le C^i x \le u_{i} && i\in P \label{cotai_inf}}

\item %\TP{cota\_inf\_j, cota\_sup\_j}
\emph{Constraints related to bounds of values in specific positions.} Let $l^\pi_j$ ({$u^\pi_j$}) denote the $j$-th lowest (largest) value of $l_i$ ($u_i$). Then, {$l^\pi_j$} ({$u^\pi_j$}) is a valid lower (upper) bound of the objective function sorted in position $j$, that is
{\formulaN{l^{\pi}_j \le \theta_j \le u^{\pi}_j && j\in P \label{cotai_inf_ord}}}

\item % \TP{cota\_inf\_i\_uij, cota\_sup\_i\_uij, cota\_inf\_j\_uij, cota\_sup\_j\_uij}
\emph{Constraints related to bounds of cost function values in specific positions.} Let $l_{ij}$ and $u_{ij}$ denote valid lower and upper bounds on the value of objective $i$ if it occupies position $j$, respectively. Then, lower and upper bounds on the value of objective $i$  are
\formulaN{\min_{j\in P}l_{ij} \le C^i x \le \max_{j\in P}u_{ij} && i\in P \label{cotai_uij}}
Analogously to \eqref{cotai_inf_ord}, we can sort the $j$-th lowest (largest) value of $\min_{j\in P}l_{ij}$ obtaining the following inequality
\formulaN{\min_{i\in P}l_{ij} \le \theta_j \le \max_{i\in P}u_{ij} && j\in P \label{cotaj_uij}}
%Note that, in general
%\formula{u_{ij} \leq u_{i} && i,j\in P}
%\formula{u_{ij} \geq l_{i} && i,j\in P}
%\formula{u_{ij} \leq u_{\pi_j} && i,j\in P}
%\formula{u_{ij} \geq l_{\pi_j} && i,j\in P}

\item %\TP{cota\_y\_z\_inf, cota\_y\_z\_sup}
There are also different bounds on the value of the cost function $i$ and the value of the cost function sorted in position $j$:
\formulaN{\sum_{j\in P} \max\{l_{i},l^{\pi}_j\} z_{ij}\leq C^ix \leq \sum_{j\in P} \min\{u_{i}, u^{\pi}_j\} z_{ij}                              && i\in P \label{cotai_uij_max}}

%\item %\TP{cota\_H\_z\_inf, cota\_H\_z\_sup}
%There are also different bounds of the cost function sorted in position $j$:
\formulaN{\sum_{i\in P} \max\{l_{i},l^{\pi}_j\} z_{ij}\leq \theta_j \leq \sum_{i\in P} \min\{u_{i}, u^{\pi}_j\} z_{ij}                              && j\in P \label{cotaj_uij_max}}

%\item \TP{cota\_y\_z\_3\_uij}
%\emph{Constraints establishing relations between variables $y$, $\theta$ and $z$} Constraint \eqref{owad'} can be strengthened replacing $M$ by a smaller bound as:
%\formulaN{C^i x\leq \theta_j+u_{ij}(1-\sum_{j'\geq j}z_{ij'})                 &&  i,j\in P   }

\item%\TP{cota\_z\_y\_3\_uij\_1, cota\_z\_y\_3\_uij\_2}\\
The inclusion of the following constraint also allows to consider, in the original formulations in Section \ref{Sec:ModelingtheOWAP}, weights $\omega \in \mathds{R}$ that, consequently, could take both negative and positive values.
\formulaN{\theta_j\leq \max_{i\in P}\{u_{ij},C^i x+M(1-z_{ij})\}  && i,j\in P\label{cotazy}}

\item %\TP{orden\_H}
\emph{Constraints related to positions in the ordering.}
Constraints \eqref{validordering} impose that the position values are ordered in non-increasing order.
\formulaN{\theta_j \geq \theta_{j+1} && j\in P\backslash\{p\} \label{validordering}}

\item %\TP{subsetY\_1, subsetY\_2, subsetY\_3, subsetY\_4}
\emph{Constraints related to subsets  of cost functions.}
Next, we observe that for any subset $I\subseteq P$, of size $k=1,...,p$
\formulaN{\sum_{i\in I}y_i\leq\sum_{j=1}^{k}\theta_j                            &&  I\subseteq P   \label{validsubsets}}
In particular, we consider the cases when $I=\{i\}$, $I=\{i,i'\in P\}$, $I=P\setminus\{i\}$ and $I=P$.

\end{itemize}

%=======================================================================================================================

%=======================================================================================================================
\subsection{Valid inequalities for the (OWAP2) formulation}\label{SecOwa2Vi}
Note first that all previous inequalities from Section \ref{SecOwa1Vi} can be applied to the two-index formulation of the OWAP substituting $\theta_j=\sum_{i\in P}y_{ij}$. Additionally, the following inequalities provide a reinforcement to the {formulations using $y$ variables:}
\begin{itemize}
\item %\TP{valor\_yy\_eq}
The following inequality combined with \eqref{owa2e} improves considerably the LP relaxation of the OWAP
\formulaN{\sum_{k\in P}y_{ik}=C^ix        &&  i\in P\label{vi:owa2eq}}
%\formulaN{\sum_{j'\in P}y_{ij'}=C^ix        &&  i\in P}

\item %\TP{cota\_yy\_z\_3\_dis1}
Constraint (\ref{owa2e}) can be disaggregated by $j\in P$ as:
%\formulaN{y_{ij}\leq\sum_{i'\in P}y_{i'j}+\min\{u_{i},u_{\pi_j} \}(1-\sum_{j'\geq j}z_{ij'})         &&  i,j\in P}
\formulaN{y_{ij}\leq\sum_{i'\in P}y_{i'j}+\min\{u_{i},u^\pi_j \}(1-\sum_{k\geq j}z_{ik})         &&  i,j\in P\label{vi:cotayydis1}}

\item %\TP{cota\_yy\_z\_3\_dis2}
We can also establish a lower bound on the value of cost function $i\in P$ if it is ordered in position $j\in P$ by  relating the $x$, $y$ and $z$ variables as follows:
\formulaN{C^ix\leq y_{ij}+u^\pi_j(1-z_{ij})         &&  i,j\in P\label{vi:cotayydis2}}

Observe that, for $i, j$ fixed, the above constraint imposes a lower bound on the value $y_{ij}$ only when cost function $i\in P$ is ordered in position $j\in P$, and becomes inactive otherwise.
\item %\TP{yy\_related1}
We can also relate the values of two different cost functions between them, depending on their positions. In particular,
\formulaN{\sum_{k \geq j+1}y_{ik}\leq y_{i'j}+u_{i}(1-z_{i'j}-z_{ij}) && i,i',j\in P, i\neq i', j\neq p \label{vi:yyrel1}}
For i, i', j fixed, constraint \eqref{vi:yyrel1} establishes that when cost function {$i'$} occupies position $j$, its value cannot be smaller than that of cost function {$i$}, provided that cost function {$i$} is ordered after $j$. Observe that the constraint becomes inactive when {$i$} is ordered before $j$ (since in this case $\sum\limits_{k \geq j+1} y_{ik}=0$) and when $i$ does not occupy position $j$.

\item %\TP{yy\_related2}
A better effectiveness of the previous inequalities can be obtained by means of
\formulaN{y_{ij+1}\leq y_{i'j}+(1-z_{ij+1})u_{ij+1}+(1-z_{i'j})u_{i'j} && i,i',j\in P, i\neq i', j\neq p\label{vi:yyrel2}}
which can be further reinforced to
\formulaN{y_{ij+1}\leq y_{i'j}+(1-z_{ij+1})\min\{u_{i},u^\pi_{j+1} \}+(1-z_{i'j})\min\{u_{i'},u^\pi_j \} && i,i',j\in P, i\neq i', j\neq p.\label{vi:yyrel3}}

%\item \TP{yy\_related3}
%Let $U$ be an upper bound on the value of any objective (for instance, $U=\max_{i\in P}u_i$. Then, for all $i,j\in P$
%\formulaN{\sum_{i' \in P} y_{i'j+1}\le y_{ij}+U(1-z_{ij}) && i,j\in P}
%Constraint \eqref{c:valor_j+1} is only activated when objective $i$ occupies position $j$. In this case, the value of the objective in position $j+1$ cannot exceed the value of objective $i$.
\end{itemize}
%=======================================================================================================================

%=======================================================================================================================
\subsection{Lower and upper bounds: Elimination tests}\label{Lower and upper bounds}
Several of the inequalities presented above use valid lower and upper bounds on the values of the different cost functions, $l_i$ and $u_i$, respectively. As mentioned above, the minimum and maximum objective value with respect to each cost function provide such bounds. However, tighter bounds can be very useful for obtaining tighter constraints. One possibility is to use lower and upper bounds on the value of each objective for the different positions in the ordering. In particular, if $L_{ij}$ and $U_{ij}$ denote valid lower and upper bounds on the value of objective $i$ if it occupies position $j$, respectively, then lower and upper bounds on the value of objective $i$  are $l_i=\min_{j\in P}L_{ij}$ and $u_i=\max_{j\in P}U_{ij}$, respectively.  For $i, j\in P$ given, $L_{ij}$ and $U_{ij}$ can be obtained in different ways. One alternative is to solve the linear programming (LP) relaxation of the formulation, both for the minimization and the maximization of cost function $i$, with the additional constraint that it occupies position $j$. In this case $L_{ij}$ ($U_{ij}$) is the optimal value of the minimization (maximization) OWAP problem in which we fix the ordering variable at value 1, i.e. $z_{ij}=1$.

Next we present simple tests which can help to eliminate some variables by fixing their values. Broadly speaking these tests compare the value of a lower bound associated with the decision of setting (or not setting) objective $i$ at position $j$ with the value of a known upper bound. If the value of the lower bound exceeds the value of the upper bound, the associated decision variable can be fixed. Any feasible solution yields a valid upper bound, which corresponds to its value with respect to the objective function. In the following we use $U$ to denote the value of the upper bound corresponding to the best-known solution. We also denote by $L_{ij}^0$  the optimal value of the minimization OWAP problem in which we fix the ordering variable at value 0, i.e. $z_{ij}=0$. Then for each $i\in P$, $j\in P$ we have
\begin{itemize}
\item If $L_{ij}>U$ then $z_{ij}=0$ (no optimal solution will have objective $i$ in position $j$).
\item If $L_{ij}^0>U$ then $z_{ij}=1$ (no optimal solution will not have objective $i$ in position $j$).
\end{itemize}
\section{The OWA problem on shortest paths and minimum cost perfect matchings} \label{Sec:Comb-OWAP}
This section presents the formulations of the combinatorial objects that we use in our computational experiments, namely shortest paths and minimum cost perfect matchings. In order to test our results we have chosen two of the most well-known formulations for these two problems. These formulations have to be combined with those presented in previous sections to provide valid OWAP models for the Shortest Path Problem (SPP) (see e.g. \citealp{Cherkassky1996, Ramaswamy2005}) and the Perfect Matching Problem (PMP) (see e.g. \citealp{Edmonds1965, Groetschel1985}).  All the details are given in what follows.
\subsection{The shortest path problem }\label{Sec:ShortestPath}
We consider now the $OWAP$ when $Q$ is the SPP (see e.g. \citealp{Cherkassky1996}). Let $G=(V, E)$ be an undirected graph with set of vertices $V$, $|V|=n$ and set of edges $E$, $|E|=m$. In addition to the sets of variables required to model the order of the $p$ cost functions ranked by non-increasing values, we will need additional variables used to model the structure of the combinatorial object (shortest path in this case). For modeling the shortest path between two selected vertices, $u_1, u_n\in V$ we use a flow-based formulation, in which binary design variables $x$ are related to continuous flow variables $\varphi$. In particular, for each $e=(u,v)\in E$ let

\[
x_{e}\equiv x_{uv}=\begin{cases}
 1 & \text{edge $e\equiv(u, v)$ is in the shortest path,}\\
 0 & \text{otherwise.}
\end{cases}
\]

\noindent As usual, paths between $u_1, u_n\in V$ can be obtained by identifying the arcs that are used when one unit of flow is sent from $u_1$ to $u_n$. For the flow variables we consider a directed network, with set of vertices $V$ and set of arcs $A$ which contains two arcs, one in each direction, associated with each edge of $E$. For each $(u,v)\in A$ we define the decision variables  $\varphi_{uv}$ which represents the amount of flow through arc $(u, v)$. Then a characterization of the domain of feasible solutions ($Q$) for the SPP is:

\begin{subequations}
\begin{align}
& \hspace{8cm}&&\notag\\
%&\hhp (SPP) \hspace{0.5cm} \min\sum_{e\in E}c^i_e x_e                           &&  i\in P                  \label{SPPa}\\
&\Hp \sum_{(u,v)\in A}\varphi_{u,v}- \sum_{(u,v)\in A}\varphi_{v,u}= 1          &&  u=u_1                   \label{SPPb}\\
&\Hp \sum_{(u,v)\in A}\varphi_{u,v}- \sum_{(u,v)\in A}\varphi_{v,u}= -1         &&  u=u_n                   \label{SPPc}\\
&\Hp \sum_{(u,v)\in A}\varphi_{u,v}- \sum_{(u,v)\in A}\varphi_{v,u}= 0          &&  u\in V\setminus\{u_1,u_n\}  \label{SPPd}\\
&\Hp \varphi_{u,v}+\varphi_{v,u}\le x_{uv}                                      && (u,v)\in E               \label{SPPe}\\
&\Hp \varphi_{uv}\ge 0                                                          && (u,v)\in A               \label{SPPf}\\
&\Hp x_e\in\{0, 1\}                                                             && e\in E                   \label{SPPg}
\end{align}
\end{subequations}

Constraints (\ref{SPPb})--(\ref{SPPd}) guarantee flow conservation at any vertex of the network.
Constraints (\ref{SPPe}) relate the $\varphi$ and $x$ variables, by imposing that all the edges used for sending flow in some direction are activated.

%\begin{subequations}
%\vspace{-0.5cm}
%\begin{align}
%& \hspace{8cm}&&\notag\\
%&\hhp (SP) \hspace{0.5cm} \min\sum_{e\in E}c^i_e x_e                                &&  i\in P                          \label{SPa}\\
%&\hp \sum_{e\in \delta(u)}x_e= 1                                                    &&  u={u^+,u^-}                     \label{SPb}\\
%&\Hp \sum_{e\in \delta(u)}x_e\in\{0,2\}                                             &&  u\in V\setminus\{u^+,u^-\}      \label{SPc}\\
%&\Hp x_e\in\{0, 1\}                                                                 &&  e\in E                          \label{SPd}
%\end{align}
%\end{subequations}

%=======================================================================================================================
%\newpage
%=======================================================================================================================
\subsection{The perfect matching problem }\label{Sec:PerfectMatching}
We consider now the $OWAP$ when $Q$ is the PMP (see e.g. \citealp{Edmonds1965}).
It is well known that the PMP is polynomially solvable by using the Blossom algorithm \citep{Edmonds1965}. However, to the best of our knowledge it is not known how such an algorithm could be used for solving an OWAP in which $Q$ is given by the set of perfect matchings on a given graph. Indeed, this can be done by using any of the OWAP formulations we have introduced in the previous sections.\\
Let $G=(V, E)$ be an undirected graph with set of vertices $V$, $|V|=n$ and set of edges $E$, $|E|=m$. In addition to the sets of variables required to model the order of the $p$ cost functions ranked by non-increasing values, we will need additional variables used to model the structure of the combinatorial object (perfect matching in this case). For modeling the perfect matching we use binary design variables $x$ associated with the edges of the graph. In particular, for each $e=(u,v)\in E$ let

\[
x_{e}\equiv x_{uv}=\begin{cases}
 1 & \text{edge $e\equiv(u, v)$ is in the matching,}\\
 0 & \text{otherwise.}
\end{cases}
\]

We introduce some additional notation. For $S\subset V$, $E(S)=\{e=(u,v)\in E\mid u,v\in S\}$ and $\delta(S)=\{e=(u,v)\in E \mid u\in S, v\notin S\}$ . When $S$ is a singleton, i.e. $S=\{u\}$ with $u\in V$ we simply write $\delta(\{u\})=\delta(u)$.
Then, a characterization of the domain of feasible solutions for the PMP ($Q$) is:

\begin{subequations}
\vspace{-0.5cm}
\begin{align}
&\Hp \hspace{8cm}&&\notag\\
&\Hp \sum_{e\in \delta(u)}x_e= 1          &&  u\in V                   \label{PM1}\\
&\Hp x_e\in\{0, 1\}           && e\in E                   \label{PM2}
\end{align}
\end{subequations}

Constraints \eqref{PM1} guarantee that in the solution the degree of every vertex is one.

%\rvv{We can reinforce the LP relaxation by adding the well known odd circuit constraints. In particular for every $S\subset V$, $|S|$ odd it holds that:}{}
%
%\formula{\sum_{e\in E(S)}x_e \le \frac{|S|-1}{2} && S\subset V: |S| \text{  odd}}

%\MP{¿No sería más correcto decir que:
%\formula{ \sum_{e\in \delta(S)}x_e \le \min\{|S|,|V|-|S|\} && S\subset V: |S|=2n-1, n\in \mathds{Z}^+, n\le \lfloor\frac{|V|}{2}\rfloor}
%o que
%\formula{ \sum_{e\in E(S)}x_e \le \frac{|S|-1}{2} && S\subset V: |S| odd}
%}
%=======================================================================================================================
%\newpage
\section{Computational experience}\label{Sec:Computational}
\TP{Descripción general del experimento}
In this section we report on the results of some computational experiments we have run, in order to compare empirically the proposed formulations and reinforcements. We have studied the OWAP over the two combinatorial objects proposed: Shortest Paths and  Minimum Cost Perfect Matchings. The best formulation obtained for each combinatorial object, has been later used for studying the proposed valid inequalities, including them one by one separately. Then, for each combinatorial object, we have obtained results for 16 basic formulations (i.e., without adding any valid inequality) plus 19 ``reinforced'' formulations. For the sake of readability, we display results in tables just for the three best basic formulations and graphics for both basic and reinforced formulations. %For further details, the reader may refer to \citet{Fernandez2013} in order to check all the results obtained in the computational experiments organized by tables.

\TP{Detalles sobre el objetivo, tamaño/número de instancias, PC}
In the computational experience we study a particular case of the OWAP operator, namely the Hurwicz criterion \citep{Hurwicz1951}, defined as $\alpha\max_{i\in P}y_i+(1-\alpha)\min_{i\in P}y_i$. This objective has been already considered when analyzing the behavior of OWA operators in multiobjective optimization (see e.g. \citealp{Galand2012}) and it is of special interest for being non-convex since the sorting weights, $\alpha$,  are not in non-increasing order (\citealp{Grzybowski2011}, \citealp{Puerto2005}). The considered values of $\alpha$  are $\{0.4, 0.6, 0.8\}$ and the number of objectives ranges in $|P|\in \{4, 7, 10\}$. Graphs generation is described below considering three different sizes of the graph according to $|V|\in \{100, 225, 400\}$. In addition, for each selection of the parameters $(|V|, p, \alpha)$, 10 instances were randomly generated so, in total, we have a set of 270 benchmark instances. All instances were solved with the MIP Xpress optimizer, under a Windows 7 environment in an Intel(R) Core(TM)i7 CPU 2.93 GHz processor and 8 GB RAM. Default values were used for all solver parameters. A CPU time limit of 600 seconds was set.

\TP{Descripción de los grafos utilizados}
For the benchmark instances, we generated square grid networks produced as with the SPGRID generator of \citet{Cherkassky1996} for both combinatorial objects. Nodes of these graphs correspond to points on the plane with integer coordinates $[x,y]$, $1\leq x\leq \sqrt{|V|}$, $1\leq y\leq \sqrt{|V|}$. These points are connected ``forward'' by arcs of the form $([x,y],[x+1,y]), \, 1\leq x< \sqrt{|V|}, \, 1\leq y\leq \sqrt{|V|}$; ``up'' by arcs of the form $([x,y],[x,y+1]), \, 1\leq x\leq \sqrt{|V|}, \, 1\leq y< \sqrt{|V|}$ and ``down'' by arcs of the form $([x,y],[x,y-1]), \, 1\leq x\leq \sqrt{|V|}, \, 1< y\leq \sqrt{|V|}$ and by arcs of the form $([x,y],[x+1,y-1]), \, 1\leq x\leq \sqrt{|V|}, \, 1< y\leq \sqrt{|V|}$. The components of the cost vectors are randomly drawn from a uniform distribution on $[1,100]$. Note also that shortest paths are computed between nodes 1 and $|V|$ whereas node $|V|$ is removed for the PMP when $|V|$ is odd.

\TP{Descripción de los elementos de las tablas numéricas (para su posterior análisis)}
Each of our tables reports the following items. Each row corresponds to a group of 10 instances with the same characteristics $(|P|, |V|, \alpha)$ indicated in the first three columns. Column $t (\#)$ reports firstly the average running time in seconds of the 10 instances of the row. In addition, if at least one instance reaches the CPU time limit, we indicate in brackets the number of instances that could be solved to optimality within the maximum CPU time limit and, in such a case, we compute the average running time by using the CPU time limit for those instances that could not be solved to optimality. Column $t^*/gap^*$ reports the biggest CPU time over the 10 instances of the group. Whenever the time limit is reached, the relative gap (indicated with a percentage \%) is reported instead. Column $\#nodes$ indicates the average number of nodes explored in the branch and bound tree and column $gap_{LR}$ reports the relative gap computed with the best solution found by the solver and the linear relaxation optima at the root node. All tables report analogous items for the different formulations described along the paper. The best three formulations for each combinatorial object  are $F^z_{R2}$, $F^{zy}_{R2}$, and $F^s$ for the SPP; and $F^z_{R1}$, $F^{zy}_{R1}$, and $F^s_{R1}$ for the PMP.  Entries  in bold remark best values among the 16 basic formulations.% (all tables are available at \citet{Fernandez2013}).

\begin{figure}[h!]
\centering
\includegraphics[width=15cm]{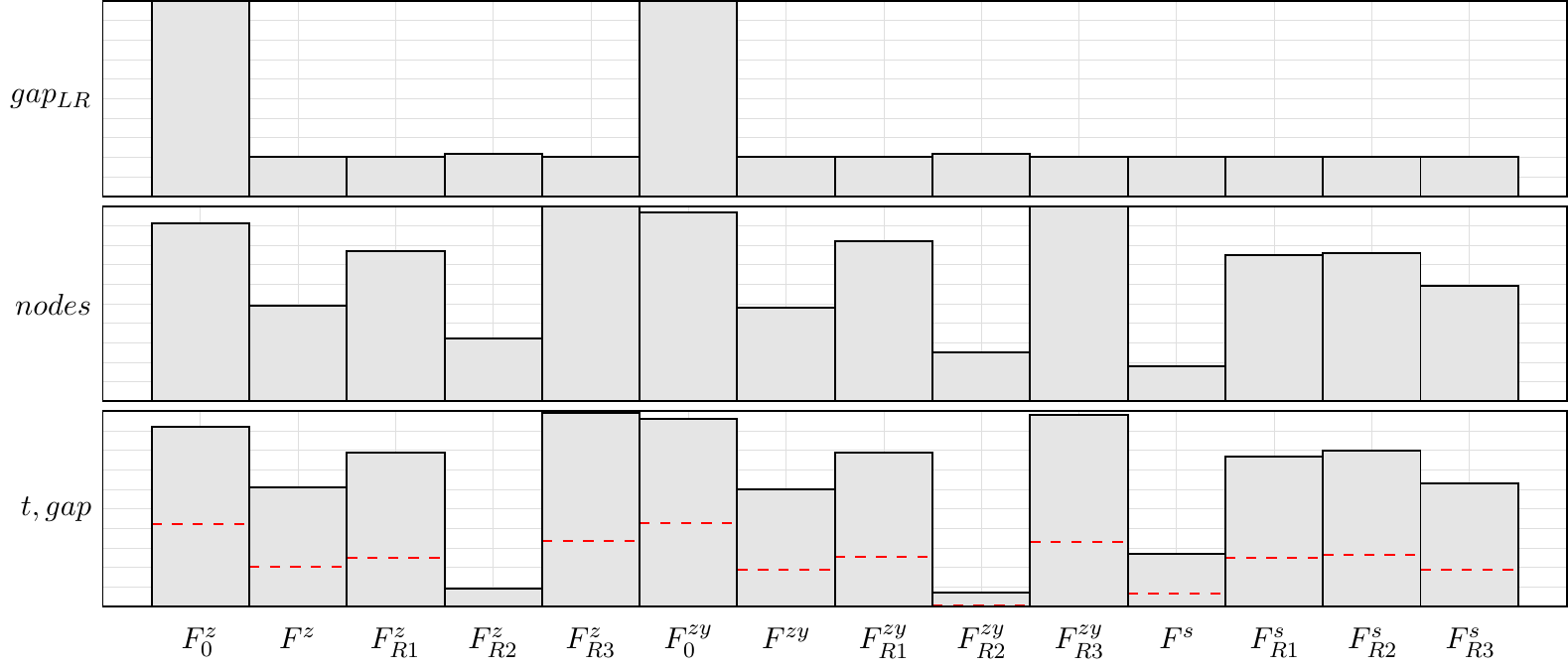}
\caption{Comparative results for the proposed OWAP basic formulations applied to the Shortest Path Problem ($p=10$, $|V|=400$)}
\label{ResFormulationsSPP}
\end{figure}

\begin{figure}[h!]
\centering
\includegraphics[width=15cm]{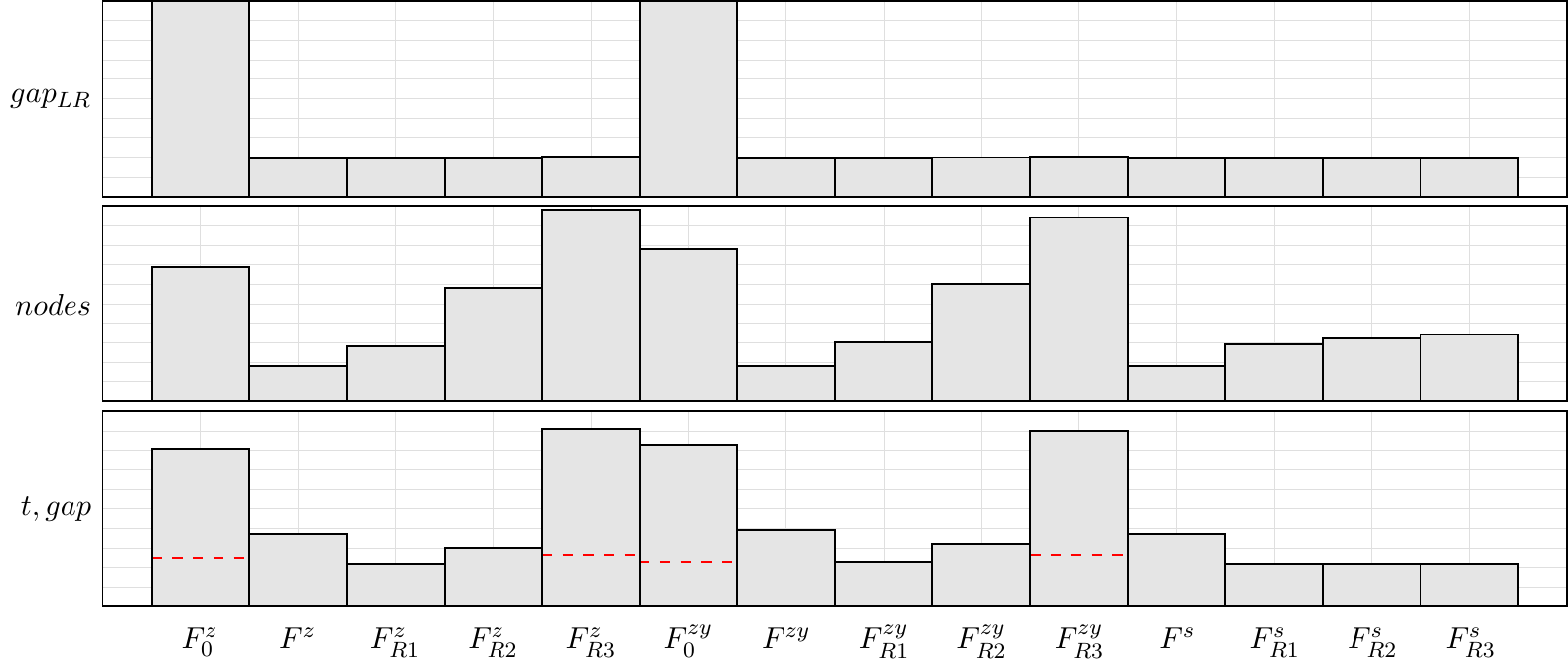}
\caption{Comparative results for the proposed OWAP basic formulations applied to the Perfect Matching Problem ($p=10$, $|V|=400$)}
\label{ResFormulationsMPM}
\end{figure}

\TP{Descripción de las gráficas (para su posterior análisis)}
Figures \ref{ResFormulationsSPP} and \ref{ResFormulationsMPM} summarize the comparative results of all proposed basic formulations applied to each combinatorial object respectively. In these graphics the $x-$axis displays the different variations of the formulations presented in Section \ref{Sec:ModelingtheOWAP} and the $y-$axis the features analyzed. All displayed bars represent percentages of mean values computed over 90 instances with $|V|=400$. These are the 90 hardest instances for the solver among the 270 we generated.

In particular the row labeled with ``$t, gap$'' shows a bar with the mean values of the running times measured in percentage over 600 seconds. For those instances reaching the time limit, we compute the mean running time taking the value of the time limit. Moreover, a dashed line indicates the percentage of worst case gap among those instances that have reached the time limit.
The columns in the row labeled with ``$nodes$'' show the percentage of nodes over $10^6$ that have been visited in the branch and bound tree.
The columns in the row labeled with ``$gap_{LR}$'' report the percentage gap relative to the best solution found by the solver and the linear relaxation optima at the root node.

%-----------------------------------------------------------------------------------------------------------------------
%-----------------------------------------------------------------------------------------------------------------------

%-----------------------------------------------------------------------------------------------------------------------
\begin{table}[h!]
\renewcommand{\esp}{@{\hspace{0.15cm}}}
\begin{center}
\begin{scriptsize}
\begin{tabular}{|@{}\esp c \esp c \esp c\esp |c\esp c\esp c\esp c\esp |c\esp c\esp c\esp c\esp |c\esp c\esp c\esp c\esp  @{}|}
\hline
\multicolumn{3}{|c|}{Inst} &
\multicolumn{4}{|c|}{$F_{R2}^z$} &
\multicolumn{4}{|c|}{$F_{R2}^{zy}$} &
\multicolumn{4}{|c|}{$F^s$}   \\
$|V|$ & $p$ & $\alpha$ &
$t(\#)$ & $t^*/gap^*$ & $\#nodes$ & $gap_{LR}$ &
$t(\#)$ & $t^*/gap^*$ & $\#nodes$ & $gap_{LR}$ &
$t(\#)$ & $t^*/gap^*$ & $\#nodes$ & $gap_{LR}$\\
\hline
100	&	4	&	0.4	&	0.5 & \textbf{0.6} & 15 & 55.79	&	 \textbf{0.4} & \textbf{0.6} & \textbf{13} & 55.79	&	8.4 & 59.9 & 16959 & \textbf{53.72}	\\
100	&	4	&	0.6	&	\textbf{0.5} & \textbf{0.6} & \textbf{41} & 40.17	&	12 & 116.5 & 56370 & 40.17	&	31.9 & 213.4 & 113492 & \textbf{37.39}	\\
100	&	4	&	0.8	&	\textbf{0.4} & \textbf{0.5} & 61 & 24.26	&	 \textbf{0.4} & \textbf{0.5} & \textbf{48} & 24.26	&	2.4 & 7.9 & 2871 & \textbf{20.79}	\\
100	&	7	&	0.4	&	\textbf{0.6} & \textbf{0.7} & 200 & 52.77	&	 \textbf{0.6} & 0.8 & \textbf{177} & 52.77	&	121 (8) & 40.66\% & 210086 & \textbf{51.11}	\\
100	&	7	&	0.6	&	\textbf{0.7} & \textbf{0.8} & \textbf{360} & 38.19	&	0.8 & 1.6 & 468 & 38.19	&	121 (8) & 22.74\% & 193167 & \textbf{36.01}	\\
100	&	7	&	0.8	&	\textbf{0.9} & 1.6 & 839 & 23.76	&	 \textbf{0.9} & \textbf{1.4} & \textbf{760} & 23.76	&	20.8 & 125.9 & 36870 & \textbf{21.08}	\\
100	&	10	&	0.4	&	\textbf{2.1} & \textbf{4.2} & \textbf{6658} & 51.98	&	2.5 & 10.3 & 9239 & 51.98	&	195.2 (7) & 43.64\% & 273035 & \textbf{49.29}	\\
100	&	10	&	0.6	&	4.1 & 13.2 & 16386 & 37.89	&	\textbf{2.8} & \textbf{11.1} & \textbf{9985} & 37.89	&	178.6 (8) & 24.44\% & 238513 & \textbf{34.4}	\\
100	&	10	&	0.8	&	5.5 & 27.9 & 23353 & 24.83	&	13.1 & 49.4 & 57599 & 24.83	&	95.3 & 500.7 & 127230 & \textbf{20.61}	\\
225	&	4	&	0.4	&	\textbf{0.8} & \textbf{1} & 48 & 55.77	&	 \textbf{0.8} & 1.1 & \textbf{45} & 55.77	&	64.4 (9) & 52.43\% & 29874 & 55	\\
225	&	4	&	0.6	&	\textbf{0.8} & \textbf{1} & \textbf{44} & 39.42	 &	\textbf{0.8} & \textbf{1} & 49 & 39.42	&	91.5 (9) & 31.77\% & 41747 & \textbf{38.31}	\\
225	&	4	&	0.8	&	\textbf{0.8} & \textbf{1.1} & 95 & 22.13	&	 \textbf{0.8} & 1.2 & \textbf{84} & 22.13	&	243.8 (6) & 14.08\% & 70842 & \textbf{20.7}	\\
225	&	7	&	0.4	&	\textbf{1.2} & \textbf{1.3} & \textbf{99} & 52.61	&	1.3 & 1.8 & 151 & 52.61	&	129 (8) & 49.52\% & 41763 & 51.29	\\
225	&	7	&	0.6	&	\textbf{3.3} & \textbf{8.8} & \textbf{1554} & 37.63	&	16.2 & 143.6 & 10871 & 37.63	&	185.6 (7) & 31.25\% & 63146 & 35.83	\\
225	&	7	&	0.8	&	4.6 & 22.1 & 3082 & 22.76	&	\textbf{2.6} & \textbf{6.1} & \textbf{1204} & 22.76	&	305 (5) & 14.19\% & 105127 & \textbf{20.44}	\\
225	&	10	&	0.4	&	9.1 & 62.7 & 6427 & 51.68	&	\textbf{5.4} & \textbf{24.9} & \textbf{4222} & 51.68	&	317.1 (5) & 49.98\% & 95076 & 50.33	\\
225	&	10	&	0.6	&	15.2 & 56.6 & 10148 & 37.07	&	\textbf{10.8} & \textbf{39.7} & \textbf{7537} & 37.07	&	319.5 (5) & 32.14\% & 96370 & \textbf{35.15}	\\
225	&	10	&	0.8	&	38.1 & 147.8 & 41223 & 23.12	&	 \textbf{29.6} & \textbf{141.5} & \textbf{32090} & 23.12	&	279.6 (6) & 15.16\% & 85419 & 20.81	\\
400	&	4	&	0.4	&	1.4 & 1.8 & 57 & 55.07	&	\textbf{1.3} & \textbf{1.6} & \textbf{55} & 55.07	&	3.3 & 16.8 & 286 & \textbf{54.44}	\\
400	&	4	&	0.6	&	1.6 & 2 & 95 & 38.71	&	\textbf{1.5} & \textbf{1.8} & \textbf{76} & 38.71	&	88.5 (9) & 35.79\% & 13806 & \textbf{37.84}	\\
400	&	4	&	0.8	&	\textbf{1.8} & \textbf{2.9} & \textbf{182} & 21.57	&	\textbf{1.8} & 3.1 & 265 & 21.57	&	255.9 (6) & 17.21\% & 49806 & \textbf{20.47}	\\
400	&	7	&	0.4	&	\textbf{6.5} & \textbf{41.1} & \textbf{1102} & 52.72	&	19.3 & 169 & 4370 & 52.72	&	76.4 (9) & 50.67\% & 9192 & \textbf{51.85}	\\
400	&	7	&	0.6	&	\textbf{9.4} & \textbf{62.6} & \textbf{2952} & 37.41	&	63.4 (9) & 33.32\% & 10416 & 37.48	&	70.9 (9) & 34.59\% & 8711 & \textbf{36.27}	\\
400	&	7	&	0.8	&	8.1 & 30.2 & \textbf{1994} & 21.87	&	 \textbf{7.2} & \textbf{24.6} & 1999 & 21.87	&	368.8 (4) & 18.29\% & 32614 & 20.41	\\
400	&	10	&	0.4	&	158.5 (9) & 1.09\% & 100979 & 51.8	&	 \textbf{116.1} & \textbf{242.7} & 83991 & 51.8	&	306.4 (5) & 48.93\% & 24184 & \textbf{50.73}	\\
400	&	10	&	0.6	&	61.8 & 121.5 & 37448 & 36.48	&	 \textbf{33.2} & \textbf{115.9} & 17308 & 36.48	&	132.4 (8) & 31.48\% & \textbf{10395} & \textbf{35.01}	\\
400	&	10	&	0.8	&	229.9 (8) & 0.61\% & 143034 & 21.82	&	155.6 (\textbf{9}) & \textbf{0.04}\% & 104042 & 21.82	&	 \textbf{142.6} (8) & 17.18\% & 12829 & 19.99	 \\

\hline
\end{tabular}
\end{scriptsize}
\end{center}
\caption{Results obtained for the three best OWAP basic formulations applied to the Shortest Path Problem}
\label{table:owap_spp_pure}
\end{table}
%-----------------------------------------------------------------------------------------------------------------------

\TP{SP: Conclusiones de las tablas y gráficas para los modelos ``puros'' aplicados al SP y PMP}
From the results displayed in Table \ref{table:owap_spp_pure} and Figure \ref{ResFormulationsSPP}, we observe first that the $gap_{LR}$ is similar for all formulations except for $F_0^{z}$ and $F_0^{zy}$, where a 100\% of gap is reached. Formulations $F_{R2}^{z}$ and $F_{R2}^{zy}$ increase slightly the $gap_{LR}$ in comparison with the remaining formulations but this does not affect negatively in the exploration as we see next. The values of $nodes$ and $t, gap$ are strongly related for each one of the formulations.  $F_{0}^{z}$, $F_{R3}^{z}$, $F_{0}^{zy}$ and $F_{R3}^{zy}$ give the worst values. In contrast, $F_{R2}^{z}$, $F_{R2}^{zy}$ and $F^{s}$ produce the best values. In addition, we observe a regular behavior among all formulations with $s$ variables, namely $F^{s}$, $F_{R1}^{s}$, $F_{R2}^{s}$ and $F_{R3}^{s}$. Regarding to the PMP, analogous conclusions can be obtained in Table \ref{table:owap_pmp_pure} and Figure \ref{ResFormulationsMPM} for the $gap_{LR}$ and the relations between $nodes$ and $t, gap$. However, in this case, formulations $F_{R1}^{z}$ and $F_{R1}^{zy}$ produce the best values together with $F_{R1}^{s}$, $F_{R2}^{s}$ and $F_{R3}^{s}$.

%-----------------------------------------------------------------------------------------------------------------------
%-----------------------------------------------------------------------------------------------------------------------

\begin{figure}[h!]
\centering
\begin{tikzpicture}[scale=1]
\tiny
\node[rectangle, scale=1, anchor=west]    () at ( 0,5) {\includegraphics[width=15cm]{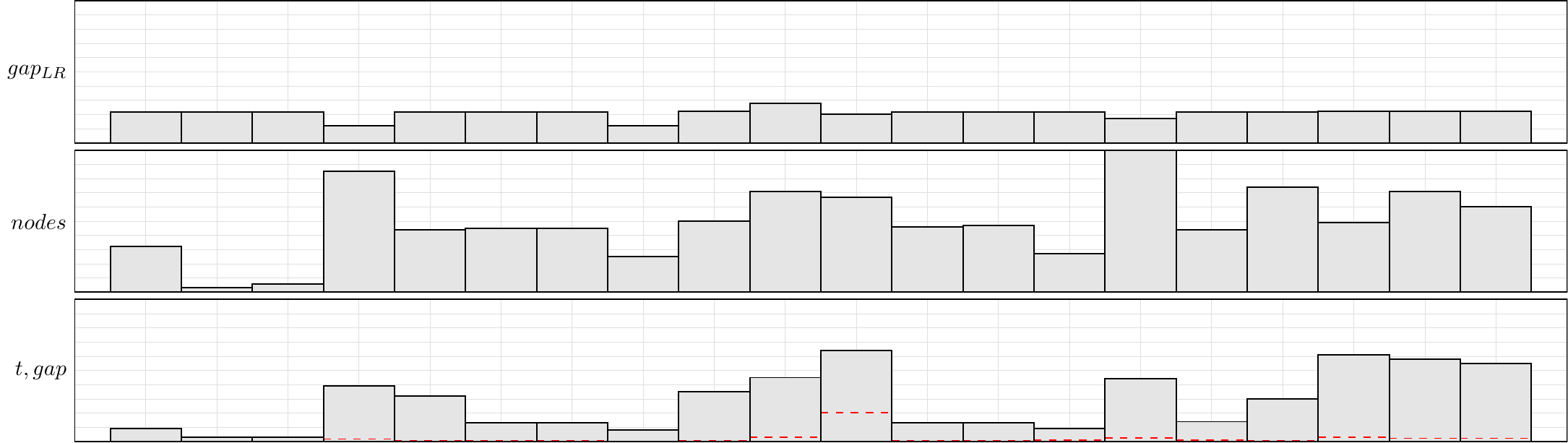}};
\newlength{\len}
\setlength\len{1.5cm}
\foreach \marca in {$F_{R2}^{zy}$, \ref{cotai_inf}.1, \ref{cotai_inf}.2, \ref{cotai_inf_ord}.1, \ref{cotai_inf_ord}.2, \ref{cotai_uij_max}.1, \ref{cotai_uij_max}.2, \ref{cotaj_uij_max}.1, \ref{cotaj_uij_max}.2, \ref{cotazy}, \ref{validsubsets}.1, \ref{validsubsets}.2, \ref{validsubsets}.3, \ref{validsubsets}.4, \ref{vi:owa2eq}, \ref{vi:cotayydis1}, \ref{vi:cotayydis2}, \ref{vi:yyrel1}, \ref{vi:yyrel2}, \ref{vi:yyrel3}}{ %Marcas de la malla grande
\node[circle, scale=1]    () at (\len,2.5) {\tiny{\marca}};
\addtolength{\len}{0.68cm}
\global\len=\len
}
\end{tikzpicture}
\caption{Comparative results for the proposed OWAP reinforced formulations applied to the Shortest Path Problem ($p=10$, $|V|=400$)}
\label{ResValIneqSPP}
\end{figure}

\begin{figure}[h!]
\centering
\begin{tikzpicture}[scale=1]
\tiny
\node[rectangle, scale=1, anchor=west]    () at ( 0,5) {\includegraphics[width=15cm]{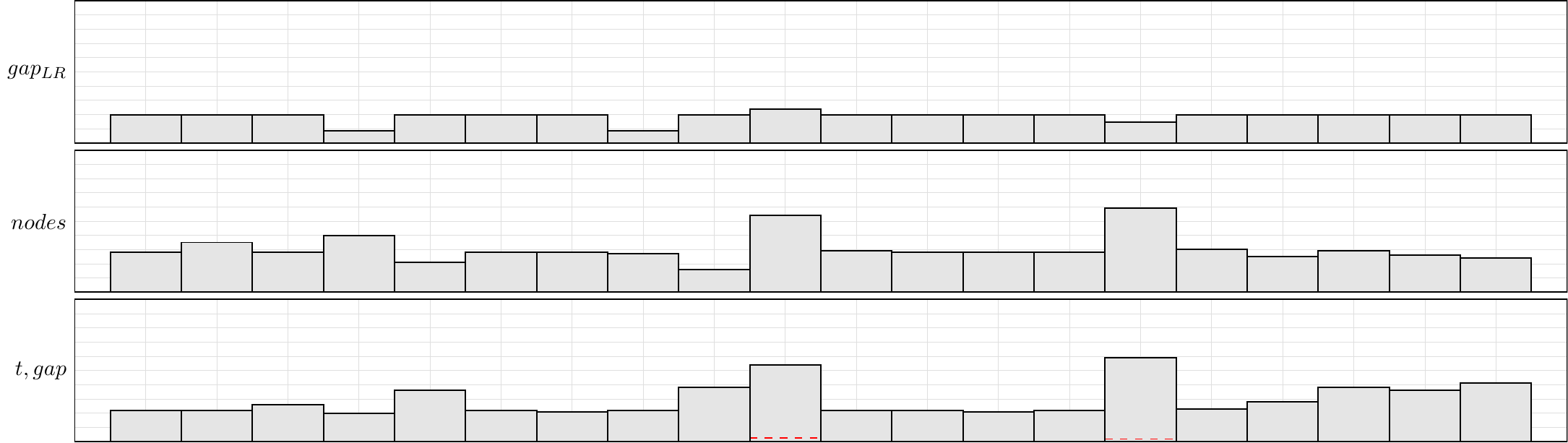}};
\setlength\len{1.5cm}
\foreach \marca in {$F_{R1}^z$, \ref{cotai_inf}.1, \ref{cotai_inf}.2, \ref{cotai_inf_ord}.1, \ref{cotai_inf_ord}.2, \ref{cotai_uij_max}.1, \ref{cotai_uij_max}.2, \ref{cotaj_uij_max}.1, \ref{cotaj_uij_max}.2, \ref{cotazy}, \ref{validsubsets}.1, \ref{validsubsets}.2, \ref{validsubsets}.3, \ref{validsubsets}.4, \ref{vi:owa2eq}, \ref{vi:cotayydis1}, \ref{vi:cotayydis2}, \ref{vi:yyrel1}, \ref{vi:yyrel2}, \ref{vi:yyrel3}}{ %Marcas de la malla grande
\node[circle, scale=1]    () at (\len,2.5) {\tiny{\marca}};
\addtolength{\len}{0.68cm}
\global\len=\len
}
\end{tikzpicture}
\caption{Comparative results for the proposed OWAP reinforced formulations applied to the Perfect Matching Problem ($p=10$, $|V|=400$)}
\label{ResValIneqMPM}
\end{figure}

%-----------------------------------------------------------------------------------------------------------------------
\begin{table}[h!]
\renewcommand{\esp}{@{\hspace{0.15cm}}}
\begin{center}
\begin{scriptsize}
\begin{tabular}{|@{}\esp c \esp c \esp c\esp |c\esp c\esp c\esp c\esp |c\esp c\esp c\esp c\esp |c\esp c\esp c\esp c\esp  @{}|}
\hline
\multicolumn{3}{|c|}{Inst} &
\multicolumn{4}{|c|}{$F_{R1}^z$} &
\multicolumn{4}{|c|}{$F_{R1}^{zy}$} &
\multicolumn{4}{|c|}{$F_{R1}^s$}   \\
$|V|$ & $p$ & $\alpha$ &
$t(\#)$ & $t^*/gap^*$ & $\#nodes$ & $gap_{LR}$ &
$t(\#)$ & $t^*/gap^*$ & $\#nodes$ & $gap_{LR}$ &
$t(\#)$ & $t^*/gap^*$ & $\#nodes$ & $gap_{LR}$\\
\hline
100	&	4	&	0.4	&	0.6 & \textbf{0.7} & 186 & \textbf{55.44}	&	 0.6 & \textbf{0.7} & 139 & \textbf{55.44}	&	\textbf{0.5} & \textbf{0.7} & 147 & \textbf{55.44}	\\
100	&	4	&	0.6	&	\textbf{0.6} & 0.7 & 152 & \textbf{38.98}	&	 \textbf{0.6} & 0.7 & 140 & \textbf{38.98}	&	\textbf{0.6} & 0.8 & 175 & \textbf{38.98}	\\
100	&	4	&	0.8	&	\textbf{0.6} & \textbf{0.7} & 302 & \textbf{21.53}	&	0.7 & 0.8 & 329 & \textbf{21.53}	&	 \textbf{0.6} & \textbf{0.7} & 157 & \textbf{21.53}	\\
100	&	7	&	0.4	&	1 & 1.3 & 236 & \textbf{52.18}	&	1 & \textbf{1.2} & 256 & \textbf{52.18}	 &	1 & \textbf{1.2} & \textbf{205} & \textbf{52.18}	\\
100	&	7	&	0.6	&	\textbf{1.1} & 1.4 & 480 & \textbf{35.97}	&	 1.2 & 1.8 & 591 & \textbf{35.97}	&	1.2 & 1.6 & 529 & \textbf{35.97}	\\
100	&	7	&	0.8	&	1.4 & 2 & 965 & \textbf{20.27}	&	1.5 & 2.3 & 1008 & \textbf{20.27}	 &	 \textbf{1.3} & \textbf{1.8} & 1075 & \textbf{20.27}	\\
100	&	10	&	0.4	&	\textbf{1.5} & 1.9 & 299 & \textbf{50.66}	&	 1.7 & 4.1 & 580 & \textbf{50.66}	&	\textbf{1.5} & 1.9 & 333 & \textbf{50.66}	\\
100	&	10	&	0.6	&	\textbf{1.9} & \textbf{2.6} & 963 & \textbf{34.85}	&	\textbf{1.9} & 2.9 & 985 & \textbf{34.85}	&	2 & 2.8 & 922 & \textbf{34.85}	\\
100	&	10	&	0.8	&	6 & 19.4 & 6329 & \textbf{20.2}	&	5.6 & 17.6 & \textbf{5364} & \textbf{20.2}	&	6.1 & 19.8 & 7018 & \textbf{20.2}	 \\
225	&	4	&	0.4	&	2.1 & 4.4 & 1188 & \textbf{55.09}	&	2 & \textbf{2.9} & 990 & \textbf{55.09}	 &	\textbf{1.9} & 4.1 & 1095 & \textbf{55.09}	\\
225	&	4	&	0.6	&	\textbf{1.7} & 2.9 & 1236 & \textbf{38.57}	&	 \textbf{1.7} & 2.5 & 1239 & \textbf{38.57}	&	\textbf{1.7} & \textbf{2.2} & \textbf{982} & \textbf{38.57}	 \\
225	&	4	&	0.8	&	\textbf{1.9} & 3.2 & 1101 & \textbf{21.09}	&	 \textbf{1.9} & 3.7 & 1240 & \textbf{21.09}	&	2 & 3.6 & 1221 & \textbf{21.09}	\\
225	&	7	&	0.4	&	\textbf{7.1} & \textbf{22.8} & 9208 & \textbf{52.34}	&	8.4 & 36 & 5617 & \textbf{52.34}	&	8.7 & 29.3 & 6308 & \textbf{52.34}	\\
225	&	7	&	0.6	&	10 & 16 & 6038 & \textbf{36.27}	&	9.7 & 18.3 & 6206 & \textbf{36.27}	&	 \textbf{8.8} & \textbf{15.9} & \textbf{5432} & \textbf{36.27}	\\
225	&	7	&	0.8	&	17.2 & 62.5 & 10491 & \textbf{20.32}	&	 17.1 & \textbf{48.7} & 10746 & \textbf{20.32}	&	\textbf{14.7} & 50.1 & \textbf{9525} & \textbf{20.32}	\\
225	&	10	&	0.4	&	7.5 & 13.2 & 2136 & \textbf{50.25}	&	7.4 & 12.4 & 2464 & \textbf{50.25}	 &	7.8 & 15.5 & 2265 & \textbf{50.25}	 \\
225	&	10	&	0.6	&	32.4 & 123.2 & 15537 & \textbf{34.56}	&	 33.9 & 90.1 & 13763 & \textbf{34.56}	&	31.5 & 70.1 & 15465 & \textbf{34.56}	\\
225	&	10	&	0.8	&	295 (8) & 0.32\% & 114029 & \textbf{19.62}	&	 338.7 (7) & 12.07\% & 130079 & 19.7	&	344.7 (8) & 0.33\% & 133025 & \textbf{19.62}	\\
400	&	4	&	0.4	&	7.3 & 22.3 & 3345 & \textbf{55.37}	&	6.3 & 15.5 & 2546 & \textbf{55.37}	 &	\textbf{6.1} & \textbf{9.6} & 2777 & \textbf{55.37}	\\
400	&	4	&	0.6	&	\textbf{6.7} & \textbf{11.9} & 4103 & \textbf{39.04}	&	7.5 & 16.7 & 4044 & \textbf{39.04}	&	8.7 & 25.3 & 6589 & \textbf{39.04}	\\
400	&	4	&	0.8	&	9 & 22.1 & 5397 & \textbf{21.03}	&	11.4 & 44.9 & 6263 & \textbf{21.03}	&	 9.2 & 19.4 & 5194 & \textbf{21.03}	 \\
400	&	7	&	0.4	&	\textbf{34.4} & 144.4 & \textbf{10464} & \textbf{52.05}	&	48.9 & 257 & 15696 & \textbf{52.05}	&	37.7 & 218.2 & 11164 & \textbf{52.05}	\\
400	&	7	&	0.6	&	83.4 & 250.9 & 27604 & \textbf{36.12}	&	 \textbf{74.5} & 209.5 & \textbf{26944} & \textbf{36.12}	&	 78.7 & \textbf{185.1} & 28692 & \textbf{36.12}	 \\
400	&	7	&	0.8	&	\textbf{84.4} & 187.6 & \textbf{28762} & \textbf{20.19}	&	98.2 & 182.5 & 35369 & \textbf{20.19}	&	92.6 & 206.4 & 34328 & \textbf{20.19}	\\
400	&	10	&	0.4	&	68.4 & \textbf{197.4} & 13777 & \textbf{50.58}	 &	86.7 & 387.2 & 17514 & \textbf{50.58}	&	91.9 & 407.6 & 19024 & \textbf{50.58}	\\
400	&	10	&	0.6	&	289.4 (9) & 0.11\% & 61886 & \textbf{34.54}	&	 335 (9) & 0.24\% & 69457 & \textbf{34.54}	&	285.7 & \textbf{563.5} & 59428 & \textbf{34.54}	\\
400	&	10	&	0.8	&	583.5 (\textbf{1}) & 0.42\% & 97022 & \textbf{19.5}	&	599 (\textbf{1}) & 0.4\% & 93171 & \textbf{19.5}	 &	\textbf{577} (\textbf{1}) & 0.43\% & 97258 & 19.52	\\

\hline
\end{tabular}
\end{scriptsize}
\end{center}
\caption{Results obtained for the three best OWAP basic formulations applied to the Perfect Matching Problem}
\label{table:owap_pmp_pure}
\end{table}
%-----------------------------------------------------------------------------------------------------------------------
\TP{Descripción de los elementos de las gráficas II (desigualdades válidas $p=10$, $|V|=400$) aplicados al SP y PMP}
Figures \ref{ResValIneqSPP} and \ref{ResValIneqMPM} report analogous items as Figures \ref{ResFormulationsSPP} and \ref{ResFormulationsMPM}, but now when the valid inequalities of Section \ref{Sec:ValIneq} are incorporated to the best basic formulations obtained for each combinatorial object. The $x-$axis displays the different variations in the formulations, starting first with the best basic formulation. Next labels refer to the valid inequality that has been added. Labels of the valid inequalities correspond with those of Section \ref{Sec:ValIneq}, where ``.1'' and ``.2'' refer to the two inequalities displayed in a single equation (for example the two valid inequalities of equation \eqref{cotai_inf} are labeled as (\ref{cotai_inf}.1) and (\ref{cotai_inf}.2)). In the following we will refer indistinctly to a valid inequality and the formulation that includes such valid inequality. All displayed bars represent percentages of mean values computed over 30 random instances with $p=10$, $|V|=400$ and $\alpha\in\{0.4,0.6,0.8\}$.

%-----------------------------------------------------------------------------------------------------------------------
%-----------------------------------------------------------------------------------------------------------------------

\TP{Conclusiones de las gráficas II (desigualdades válidas $p=10$, $|V|=400$) aplicados al SP y PMP}
From the results displayed in Figure \ref{ResValIneqSPP}, we observe first that the $gap_{LR}$ is similar for all formulations but (\ref{cotai_inf_ord}.1), (\ref{cotaj_uij_max}.1), (\ref{cotazy}), (\ref{validsubsets}.1) and (\ref{vi:owa2eq}). As compared with with $F_{R2}^{zy}$, formulation (\ref{cotaj_uij_max}.1) improves the values of $gap_{LR}$, $nodes$ and $t, gap$. However, (\ref{cotai_inf_ord}.1), (\ref{validsubsets}.1) and (\ref{vi:owa2eq}) improve $gap_{LR}$ but are not able to improve  $nodes$ or $t, gap$. We also note that (\ref{cotazy}) increases $gap_{LR}$ since this gap is computed with a (low quality) best solution found by the solver and the linear relaxation optima at the root node. In addition, formulations (\ref{cotai_inf}.1), (\ref{cotai_inf}.2) and (\ref{validsubsets}.4) provide promising results in comparison with the values of $nodes$ and $t, gap$.

%-----------------------------------------------------------------------------------------------------------------------
%-----------------------------------------------------------------------------------------------------------------------

From the results displayed in Figure \ref{ResValIneqMPM}, we observe first that the $gap_{LR}$ is similar for all formulations but (\ref{cotai_inf_ord}.1), (\ref{cotaj_uij_max}.1), (\ref{cotazy}) and (\ref{vi:owa2eq}). As compared with $F_{R1}^z$, formulations (\ref{cotai_inf_ord}.1) and (\ref{cotaj_uij_max}.1), improve $gap_{LR}$ and $nodes$ or $t, gap$. However, (\ref{vi:owa2eq}) improves $gap_{LR}$ but is not able to improve  $nodes$ or $t, gap$ in comparison with the best basic formulation for PMP, namely $F_{R1}^z$. We also note that (\ref{cotazy}) increases $gap_{LR}$ since this gap is computed with a (low quality) best solution found by the solver and the linear relaxation optima at the root node. In addition, formulations (\ref{cotai_uij_max}.2), (\ref{cotaj_uij_max}.2) and (\ref{validsubsets}.3) provide promising results in comparison with the values of $nodes$ or $t, gap$.

In summary, we observe the performance of the OWAP formulation depends on its combination with the considered combinatorial object. In particular we conclude, from our computational experience, that for the SPP, it is convenient to apply $F^{zy}_{R2}$ reinforced with (\ref{cotai_inf}.1) and (\ref{cotai_inf}.2); although rather similar results can be obtained with $F^{z}_{R2}$. The conclusion for the PMP is different, because the best basic formulation is now $F^z_{R1}$ and the reinforcements (\ref{cotai_inf_ord}.1). Once more, rather similar results are obtained for $F^{zy}_{R1}$ and $F^s_{R1}$. Therefore, we can not conclude whether there is a formulation superior to all the others regardless the domain $Q$ to be considered. For this reason it is important to have developed the catalogue of formulations and valid inequalities presented in this paper. In general, it is advisable to test them depending on the combinatorial object to be considered.

%%-----------------------------------------------------------------------------------------------------------------------
%\TP{Descripción de las gráficas III y IV (formulaciones $p=10$, $|V|=400$, Hurwicz)}
%Figure \ref{ResFormulationsCeil} reports analogous items as Figure \ref{ResFormulations} applying a particular case of the OWAP operator, namely the minimization of the $\left\lceil\frac{p}{2}\right\rceil$ greatest objectives defined by giving values to the components of $\omega$ as $\omega_j=1 \iif j\leq \left\lceil\frac{p}{2}\right\rceil$, or $0$ otherwise. This objective function is of special interest since...  All displayed bars represent percentages of mean values computed over 10 random instances with $p=8$, $|V|=225$. Figure \ref{ResFormulationsCeil2} shows analogous results for $p=10$, $|V|=400$.
%
%\TP{Conclusiones de las gráficas II (desigualdades válidas $p=10$, $|V|=400$, Hurwicz)}
%=======================================================================================================================

%=======================================================================================================================
\newpage
\section{Appendix: Complete results obtained in the experiments}
%=======================================================================================================================
%\subsubsection*{SP: Formulations $F_0^z$, $F^z$, $F_{R1}^z$, $F_{R2}^z$, $F_{R3}^z$}
%-----------------------------------------------------------------------------------------------------------------------
\begin{table}[h!]
\renewcommand{\esp}{@{\hspace{0.15cm}}}
\begin{center}
\begin{scriptsize}
\begin{tabular}{|@{}\esp c \esp c \esp c\esp |c\esp c\esp c\esp c\esp |c\esp c\esp c\esp c\esp |c\esp c\esp c\esp c\esp  @{}|}
\hline
\multicolumn{3}{|c|}{Inst} &
\multicolumn{4}{|c|}{$F_0^z$} &
\multicolumn{4}{|c|}{$F^z$} &
\multicolumn{4}{|c|}{$F^z$$^\prime$}   \\
$|V|$ & $p$ & $\alpha$ &
$t(\#)$ & $t^*/gap^*$ & $\#nodes$ & $gap_{LR}$ &
$t(\#)$ & $t^*/gap^*$ & $\#nodes$ & $gap_{LR}$ &
$t(\#)$ & $t^*/gap^*$ & $\#nodes$ & $gap_{LR}$\\
\hline
100	&	4	&	0.4	&	61.9 (9) & 41.15\% & 159186 & 100	&	2.2 & 12.3 & 2087 & \textbf{53.72}	&	1.4 & 4.5 & 1005 & \textbf{53.72}	\\		
100	&	4	&	0.6	&	1.6 & 2.5 & 2775 & 100	&	1.9 & 6.7 & 2149 & \textbf{37.39}	&	17.1 & 117.1 & 72137 & \textbf{37.39}	\\		
100	&	4	&	0.8	&	15.7 & 64 & 47308 & 100	&	8.9 & 38.9 & 15109 & \textbf{20.79}	&	4.9 & 30.1 & 4885 & \textbf{20.79}	\\		
100	&	7	&	0.4	&	539.6 (1) & 55.52\% & 670773 & 100	&	73.8 (9) & 38.38\% & 169704 & \textbf{51.11}	&	114.9 (9) & 40.83\% & 252250 & \textbf{51.11}	\\		
100	&	7	&	0.6	&	599.2 (0) & 42.4\% & 821449 & 100	&	71 (9) & 19.84\% & 177585 & \textbf{36.01}	&	62.1 (9) & 16.07\% & 143671 & \textbf{36.01}	\\		
100	&	7	&	0.8	&	599.2 (0) & 30.88\% & 977772 & 100	&	21.8 & 96.5 & 45613 & \textbf{21.08}	&	19.9 & 118 & 42932 & \textbf{21.08}	\\		
100	&	10	&	0.4	&	599.2 (0) & 59.33\% & 598177 & 100	&	181.6 (8) & 38.78\% & 299151 & \textbf{49.29}	&	186.8 (7) & 40.11\% & 204286 & \textbf{49.29}	\\		
100	&	10	&	0.6	&	599.4 (0) & 47.99\% & 721323 & 100	&	131.6 (8) & 22.48\% & 184532 & \textbf{34.4}	&	192.6 (7) & 23.86\% & 237500 & 34.43	\\		
100	&	10	&	0.8	&	599.5 (0) & 36.12\% & 816612 & 100	&	70.2 & 479.4 & 101591 & \textbf{20.61}	&	57.8 & 376.8 & 83723 & \textbf{20.61}	\\		
225	&	4	&	0.4	&	4 & 11.5 & 1814 & 100	&	129.7 (8) & 51.55\% & 48659 & \textbf{54.96}	&	138.9 (8) & 53.16\% & 52164 & 55.06	\\		
225	&	4	&	0.6	&	365 (6) & 19.25\% & 308001 & 100	&	138.3 (8) & 32.58\% & 47943 & \textbf{38.31}	&	294.5 (6) & 33.99\% & 133651 & 38.39	\\		
225	&	4	&	0.8	&	599.2 (0) & 19.46\% & 405027 & 100	&	273 (7) & 12.11\% & 130154 & \textbf{20.7}	&	350.7 (5) & 14.57\% & 130975 & \textbf{20.7}	\\		
225	&	7	&	0.4	&	599.4 (0) & 62.68\% & 216392 & 100	&	361.5 (4) & 48.61\% & 181096 & 51.74	&	151.4 (9) & 39.95\% & 58002 & \textbf{51.2}	\\		
225	&	7	&	0.6	&	599.5 (0) & 51.17\% & 253663 & 100	&	363.1 (4) & 32.18\% & 176663 & 36.05	&	361 (4) & 32.66\% & 146512 & 35.83	\\		
225	&	7	&	0.8	&	599.5 (0) & 41.18\% & 265937 & 100	&	420.8 (3) & 14.49\% & 221975 & \textbf{20.44}	&	288.8 (6) & 14.55\% & 133194 & 20.45	\\		
225	&	10	&	0.4	&	599.7 (0) & 63.54\% & 166127 & 100	&	396.8 (4) & 50.6\% & 118651 & 50.61	&	332.3 (5) & 47.82\% & 92470 & \textbf{50.3}	\\		
225	&	10	&	0.6	&	599.7 (0) & 55.79\% & 208149 & 100	&	421.8 (3) & 32.03\% & 125420 & 35.43	&	335.9 (5) & 32.85\% & 111328 & 35.36	\\		
225	&	10	&	0.8	&	599.7 (0) & 47.37\% & 227825 & 100	&	324.9 (5) & 15.12\% & 97109 & 20.84	&	298.3 (6) & 14.91\% & 109912 & 20.71	\\		
400	&	4	&	0.4	&	211.8 (7) & 46.17\% & 44001 & 100	&	200.9 (7) & 52.86\% & 32709 & 54.47	&	227.4 (7) & 51.74\% & 43043 & 54.46	\\		
400	&	4	&	0.6	&	572.4 (1) & 40.28\% & 123465 & 100	&	371.9 (4) & 37.16\% & 51494 & 38.03	&	182.5 (8) & 36.23\% & 31789 & 37.89	\\		
400	&	4	&	0.8	&	599.3 (0) & 63.06\% & 146115 & 100	&	464.5 (3) & 17.95\% & 72353 & \textbf{20.47}	&	391.2 (5) & 17.89\% & 64058 & \textbf{20.47}	\\		
400	&	7	&	0.4	&	599.2 (0) & 64.86\% & 78332 & 100	&	599.8 (0) & 53.96\% & 85672 & 53.77	&	496.5 (2) & 54.52\% & 51199 & 53.01	\\		
400	&	7	&	0.6	&	599.3 (0) & 54.05\% & 97758 & 100	&	540.5 (1) & 36.56\% & 77694 & 37.07	&	550.7 (1) & 36.37\% & 53000 & 36.78	\\		
400	&	7	&	0.8	&	599.3 (0) & 45.98\% & 106458 & 100	&	481.7 (2) & 18.53\% & 70106 & 20.52	&	599.7 (0) & 17.89\% & 59361 & 20.42	\\		
400	&	10	&	0.4	&	599.5 (0) & 65.97\% & 59878 & 100	&	314.7 (5) & 50.22\% & \textbf{22479} & 50.93	&	376.5 (4) & 52.17\% & 23753 & 51.09	\\		
400	&	10	&	0.6	&	599.3 (0) & 56.75\% & 82904 & 100	&	198.9 (7) & 34.1\% & 14858 & 35.12	&	225.6 (7) & 32.13\% & 18486 & 35.06	\\		
400	&	10	&	0.8	&	599.4 (0) & 52.8\% & 82418 & 100	&	145.1 (8) & 15.5\% & \textbf{11060} & \textbf{19.98}	&	258.4 (6) & 17.55\% & 21698 & 20.04	\\

\hline
\end{tabular}
\end{scriptsize}
\end{center}
\caption{Results obtained for OWAP formulations with the Shortest Path Problem}
%\label{tab_exact_resolution}
\end{table}
%-----------------------------------------------------------------------------------------------------------------------
%-----------------------------------------------------------------------------------------------------------------------
\begin{table}[h!]
\renewcommand{\esp}{@{\hspace{0.15cm}}}
\begin{center}
\begin{scriptsize}
\begin{tabular}{|@{}\esp c \esp c \esp c\esp |c\esp c\esp c\esp c\esp |c\esp c\esp c\esp c\esp |c\esp c\esp c\esp c\esp  @{}|}
\hline
\multicolumn{3}{|c|}{Inst} &
\multicolumn{4}{|c|}{$F_{R1}^z$} &
\multicolumn{4}{|c|}{$F_{R2}^z$} &
\multicolumn{4}{|c|}{$F_{R3}^z$}   \\
$|V|$ & $p$ & $\alpha$ &
$t(\#)$ & $t^*/gap^*$ & $\#nodes$ & $gap_{LR}$ &
$t(\#)$ & $t^*/gap^*$ & $\#nodes$ & $gap_{LR}$ &
$t(\#)$ & $t^*/gap^*$ & $\#nodes$ & $gap_{LR}$\\
\hline
100	&	4	&	0.4	&	4.5 & 17.9 & 7045 & \textbf{53.72}	&	0.5 & \textbf{0.6} & 15 & 55.79	&	360.7 (4) & 35.38\% & 1088918 & \textbf{53.72}	\\		
100	&	4	&	0.6	&	8.4 & 43.7 & 20767 & \textbf{37.39}	&	\textbf{0.5} & \textbf{0.6} & \textbf{41} & 40.17	&	398.1 (4) & 19.95\% & 940627 & \textbf{37.39}	\\		
100	&	4	&	0.8	&	6.3 & 21.4 & 9361 & \textbf{20.79}	&	\textbf{0.4} & \textbf{0.5} & 61 & 24.26	&	28.1 & 84.1 & 73283 & \textbf{20.79}	\\		
100	&	7	&	0.4	&	183.8 (7) & 39.39\% & 505058 & 51.21	&	\textbf{0.6} & \textbf{0.7} & 200 & 52.77	&	599.2 (0) & 44.67\% & 1267412 & 51.52	\\		
100	&	7	&	0.6	&	177.4 (8) & 22.51\% & 552821 & \textbf{36.01}	&	\textbf{0.7} & \textbf{0.8} & \textbf{360} & 38.19	&	599.2 (0) & 28.4\% & 1451450 & 36.12	\\		
100	&	7	&	0.8	&	31.2 & 79 & 68528 & \textbf{21.08}	&	\textbf{0.9} & 1.6 & 839 & 23.76	&	494.3 (4) & 6.67\% & 2046319 & \textbf{21.08}	\\		
100	&	10	&	0.4	&	131 (8) & 38.43\% & 211989 & \textbf{49.29}	&	\textbf{2.1} & \textbf{4.2} & \textbf{6658} & 51.98	&	560.6 (1) & 47.89\% & 1405677 & 49.98	\\		
100	&	10	&	0.6	&	132.7 (8) & 18.42\% & 360647 & 34.43	&	4.1 & 13.2 & 16386 & 37.89	&	599.4 (0) & 29.69\% & 1566146 & 34.86	\\		
100	&	10	&	0.8	&	86.2 & 388.2 & 176071 & \textbf{20.61}	&	5.5 & 27.9 & 23353 & 24.83	&	539.1 (2) & 13.52\% & 1982531 & \textbf{20.61}	\\		
225	&	4	&	0.4	&	128.9 (8) & 52.67\% & 66387 & 55.01	&	\textbf{0.8} & \textbf{1} & 48 & 55.77	&	599.4 (0) & 52.65\% & 346252 & 55.09	\\		
225	&	4	&	0.6	&	308.3 (5) & 33.54\% & 154345 & 38.39	&	\textbf{0.8} & \textbf{1} & \textbf{44} & 39.42	&	360 (4) & 34.52\% & 207131 & 38.34	\\		
225	&	4	&	0.8	&	361.9 (5) & 15.51\% & 140540 & \textbf{20.7}	&	\textbf{0.8} & \textbf{1.1} & 95 & 22.13	&	439.3 (4) & 15.58\% & 206648 & \textbf{20.7}	\\		
225	&	7	&	0.4	&	314.6 (5) & 49.58\% & 160573 & 51.65	&	\textbf{1.2} & \textbf{1.3} & \textbf{99} & 52.61	&	599.6 (0) & 50.61\% & 312524 & 52.1	\\		
225	&	7	&	0.6	&	441.4 (3) & 32.37\% & 233819 & 35.92	&	\textbf{3.3} & \textbf{8.8} & \textbf{1554} & 37.63	&	599.7 (0) & 32.98\% & 336958 & 36.43	\\		
225	&	7	&	0.8	&	526.7 (2) & 13.95\% & 308827 & \textbf{20.44}	&	4.6 & 22.1 & 3082 & 22.76	&	599.7 (0) & 14.79\% & 307722 & \textbf{20.44}	\\		
225	&	10	&	0.4	&	554.5 (1) & 50.01\% & 270862 & 51.07	&	9.1 & 62.7 & 6427 & 51.68	&	599.7 (0) & 51.52\% & 331748 & 52.53	\\		
225	&	10	&	0.6	&	599.8 (0) & 31.67\% & 337820 & 35.75	&	15.2 & 56.6 & 10148 & 37.07	&	599.8 (0) & 33.34\% & 334555 & 36.52	\\		
225	&	10	&	0.8	&	599.7 (0) & 14.51\% & 279970 & \textbf{20.69}	&	38.1 & 147.8 & 41223 & 23.12	&	599.7 (0) & 16.15\% & 345339 & 20.85	\\		
400	&	4	&	0.4	&	195.8 (7) & 52.37\% & 33251 & \textbf{54.44}	&	1.4 & 1.8 & 57 & 55.07	&	599.7 (0) & 54.69\% & 118795 & 54.6	\\		
400	&	4	&	0.6	&	325 (5) & 35.5\% & 71726 & 38.01	&	1.6 & 2 & 95 & 38.71	&	539.9 (1) & 36.64\% & 111458 & 38.04	\\		
400	&	4	&	0.8	&	367.4 (4) & 17.8\% & 43381 & \textbf{20.47}	&	\textbf{1.8} & \textbf{2.9} & \textbf{182} & 21.57	&	599.8 (0) & 18.03\% & 102960 & \textbf{20.47}	\\		
400	&	7	&	0.4	&	541.2 (1) & 54.39\% & 74342 & 53.16	&	\textbf{6.5} & \textbf{41.1} & \textbf{1102} & 52.72	&	599.7 (0) & 53.57\% & 121813 & 52.76	\\		
400	&	7	&	0.6	&	541.6 (1) & 35.17\% & 92740 & 36.53	&	\textbf{9.4} & \textbf{62.6} & \textbf{2952} & 37.41	&	599.7 (0) & 36.68\% & 132839 & 36.89	\\		
400	&	7	&	0.8	&	544.3 (1) & 17.86\% & 106851 & 20.4	&	8.1 & 30.2 & \textbf{1994} & 21.87	&	599.7 (0) & 18.53\% & 110457 & 20.57	\\		
400	&	10	&	0.4	&	531.6 (2) & 53.98\% & 89609 & 52.43	&	158.5 (9) & 1.09\% & 100979 & 51.8	&	599.7 (0) & 54.33\% & 101358 & 53.88	\\		
400	&	10	&	0.6	&	599.7 (0) & 35.53\% & 84013 & 36.4	&	61.8 & 121.5 & 37448 & 36.48	&	599.7 (0) & 36.05\% & 129485 & 36.6	\\		
400	&	10	&	0.8	&	599.7 (0) & 17.75\% & 96917 & 20.25	&	229.9 (8) & 0.61\% & 143034 & 21.82	&	599.7 (0) & 17.86\% & 111753 & 20.2	\\

\hline
\end{tabular}
\end{scriptsize}
\end{center}
\caption{Results obtained for the OWAP formulations with the Shortest Path Problem}
%\label{tab_exact_resolution}
\end{table}
%-----------------------------------------------------------------------------------------------------------------------
%=======================================================================================================================
\newpage
%=======================================================================================================================
%\subsubsection*{SP: Formulations $F_0^{zy}$, $F^{zy}$, $F_{R1}^{zy}$, $F_{R2}^{zy}$, $F_{R3}^{zy}$}
%-----------------------------------------------------------------------------------------------------------------------
\begin{table}[h!]
\renewcommand{\esp}{@{\hspace{0.15cm}}}
\begin{center}
\begin{scriptsize}
\begin{tabular}{|@{}\esp c \esp c \esp c\esp |c\esp c\esp c\esp c\esp |c\esp c\esp c\esp c\esp |c\esp c\esp c\esp c\esp  @{}|}
\hline
\multicolumn{3}{|c|}{Inst} &
\multicolumn{4}{|c|}{$F_{0}^{zy}$} &
\multicolumn{4}{|c|}{$F_{}^{zy}$} &
\multicolumn{4}{|c|}{$F_{}^{zy}$$^\prime$}   \\
$|V|$ & $p$ & $\alpha$ &
$t(\#)$ & $t^*/gap^*$ & $\#nodes$ & $gap_{LR}$ &
$t(\#)$ & $t^*/gap^*$ & $\#nodes$ & $gap_{LR}$ &
$t(\#)$ & $t^*/gap^*$ & $\#nodes$ & $gap_{LR}$\\
\hline
100	&	4	&	0.4	&	1.2 & 3.1 & 1192 & 100	&	121.2 (8) & 31.45\% & 451392 & \textbf{53.72}	&	85 (9) & 33.02\% & 416097 & \textbf{53.72}	\\		
100	&	4	&	0.6	&	2.9 & 9.1 & 5786 & 100	&	15.1 & 104.1 & 29720 & \textbf{37.39}	&	4.4 & 25.5 & 7032 & \textbf{37.39}	\\		
100	&	4	&	0.8	&	14.9 & 29.7 & 43225 & 100	&	5.7 & 22.2 & 8568 & \textbf{20.79}	&	4.4 & 25 & 4359 & \textbf{20.79}	\\		
100	&	7	&	0.4	&	599.1 (0) & 56.76\% & 777811 & 100	&	82.3 (9) & 35.05\% & 196817 & \textbf{51.11}	&	121.8 (8) & 33.93\% & 232458 & 51.19	\\		
100	&	7	&	0.6	&	599.1 (0) & 44.1\% & 815725 & 100	&	63.6 (9) & 8.86\% & 165583 & \textbf{36.01}	&	67.9 (9) & 18.4\% & 97103 & \textbf{36.01}	\\		
100	&	7	&	0.8	&	599.2 (0) & 26.56\% & 992683 & 100	&	20.2 & 107.4 & 42676 & \textbf{21.08}	&	11.8 & 53.3 & 18601 & \textbf{21.08}	\\		
100	&	10	&	0.4	&	599.3 (0) & 59.5\% & 584890 & 100	&	132.3 (8) & 42.34\% & 220685 & \textbf{49.29}	&	187.6 (7) & 38.16\% & 228058 & \textbf{49.29}	\\		
100	&	10	&	0.6	&	599.4 (0) & 47.08\% & 698012 & 100	&	176.2 (8) & 23.12\% & 325255 & \textbf{34.4}	&	184.5 (7) & 24.44\% & 158769 & \textbf{34.4}	\\		
100	&	10	&	0.8	&	599.5 (0) & 37.21\% & 832031 & 100	&	83.1 & 260.9 & 118290 & \textbf{20.61}	&	122.9 & 499.5 & 148775 & \textbf{20.61}	\\		
225	&	4	&	0.4	&	42.7 & 293.1 & 26879 & 100	&	129.8 (8) & 49.74\% & 57642 & 55.01	&	285.9 (6) & 52.02\% & 123411 & 54.97	\\		
225	&	4	&	0.6	&	404.3 (5) & 17.48\% & 335937 & 100	&	185.7 (7) & 31.19\% & 99753 & 38.36	&	362.6 (4) & 34.45\% & 160383 & 38.39	\\		
225	&	4	&	0.8	&	599.7 (0) & 14.1\% & 436690 & 100	&	327.3 (5) & 14.76\% & 103352 & \textbf{20.7}	&	540.4 (1) & 14.31\% & 223399 & \textbf{20.7}	\\		
225	&	7	&	0.4	&	599.9 (0) & 61.69\% & 216243 & 100	&	480.2 (2) & 50.43\% & 227345 & 51.59	&	374.2 (4) & 51.31\% & 135242 & 51.69	\\		
225	&	7	&	0.6	&	599.5 (0) & 57.25\% & 248810 & 100	&	420.6 (3) & 32.51\% & 210813 & 35.87	&	448.8 (3) & 33.08\% & 152915 & 35.89	\\		
225	&	7	&	0.8	&	599.9 (0) & 51.27\% & 275751 & 100	&	307.1 (5) & 14.84\% & 151530 & 20.45	&	529.6 (3) & 15.35\% & 205445 & 20.56	\\		
225	&	10	&	0.4	&	599.8 (0) & 64.4\% & 165068 & 100	&	421.5 (3) & 49.78\% & 121243 & 50.58	&	311.9 (5) & 47.81\% & 87477 & 50.36	\\		
225	&	10	&	0.6	&	599.9 (0) & 64.91\% & 211343 & 100	&	311.6 (5) & 31.66\% & 93967 & 35.35	&	303.6 (5) & 32.22\% & 89976 & 35.68	\\		
225	&	10	&	0.8	&	600 (0) & 47.48\% & 215323 & 100	&	247.5 (6) & 15.04\% & 73111 & 20.71	&	321.3 (5) & 15.62\% & 89902 & 20.73	\\		
400	&	4	&	0.4	&	361.2 (4) & 49.35\% & 82543 & 100	&	153 (8) & 52.96\% & 27764 & \textbf{54.44}	&	100.9 (9) & 53.09\% & 19609 & \textbf{54.44}	\\		
400	&	4	&	0.6	&	598 (1) & 38.99\% & 145270 & 100	&	298.6 (6) & 36.15\% & 56477 & 37.91	&	371.1 (4) & 36.31\% & 51542 & 38.01	\\		
400	&	4	&	0.8	&	599.8 (0) & 34.48\% & 131878 & 100	&	425.9 (3) & 17.52\% & 63344 & \textbf{20.47}	&	432.6 (3) & 18.51\% & 49572 & \textbf{20.47}	\\		
400	&	7	&	0.4	&	599.8 (0) & 64.49\% & 78184 & 100	&	481 (2) & 52.37\% & 68194 & 53.06	&	563.9 (1) & 53.36\% & 55424 & 53.02	\\		
400	&	7	&	0.6	&	599.8 (0) & 53.68\% & 96075 & 100	&	540.4 (1) & 36.55\% & 76295 & 37.16	&	540.3 (1) & 36.44\% & 44656 & 36.82	\\		
400	&	7	&	0.8	&	599.8 (0) & 42.76\% & 108088 & 100	&	481.8 (2) & 18.09\% & 71226 & 20.37	&	450.8 (3) & 18.24\% & 68513 & 20.4	\\		
400	&	10	&	0.4	&	599.7 (0) & 66.58\% & 58295 & 100	&	371.3 (4) & 50.51\% & 29467 & 50.99	&	371.6 (5) & 50.67\% & 24697 & 50.84	\\		
400	&	10	&	0.6	&	599.7 (0) & 56.48\% & 83222 & 100	&	252.2 (6) & 33.31\% & 20435 & 35.04	&	278.1 (6) & 33.63\% & 23601 & 35.2	\\		
400	&	10	&	0.8	&	599.7 (0) & 49.98\% & 84953 & 100	&	259.4 (6) & 16.84\% & 21451 & 19.99	&	148.3 (8) & 17.27\% & 11315 & 20.03	\\		
\hline
\end{tabular}
\end{scriptsize}
\end{center}
\caption{Results obtained for the OWAP formulations with the Shortest Path Problem}
%\label{tab_exact_resolution}
\end{table}
%-----------------------------------------------------------------------------------------------------------------------

%-----------------------------------------------------------------------------------------------------------------------
\begin{table}[h!]
\renewcommand{\esp}{@{\hspace{0.15cm}}}
\begin{center}
\begin{scriptsize}
\begin{tabular}{|@{}\esp c \esp c \esp c\esp |c\esp c\esp c\esp c\esp |c\esp c\esp c\esp c\esp |c\esp c\esp c\esp c\esp  @{}|}
\hline
\multicolumn{3}{|c|}{Inst} &
\multicolumn{4}{|c|}{$F_{R1}^{zy}$} &
\multicolumn{4}{|c|}{$F_{R2}^{zy}$} &
\multicolumn{4}{|c|}{$F_{R3}^{zy}$}   \\
$|V|$ & $p$ & $\alpha$ &
$t(\#)$ & $t^*/gap^*$ & $\#nodes$ & $gap_{LR}$ &
$t(\#)$ & $t^*/gap^*$ & $\#nodes$ & $gap_{LR}$ &
$t(\#)$ & $t^*/gap^*$ & $\#nodes$ & $gap_{LR}$\\
\hline
100	&	4	&	0.4	&	65.6 (9) & 23.26\% & 241924 & \textbf{53.72}	&	\textbf{0.4} & \textbf{0.6} & \textbf{13} & 55.79	&	322 (5) & 42.47\% & 845867 & \textbf{53.72}	\\		
100	&	4	&	0.6	&	7.5 & 40.8 & 15505 & \textbf{37.39}	&	12 & 116.5 & 56370 & 40.17	&	241.6 (8) & 19.94\% & 672083 & \textbf{37.39}	\\		
100	&	4	&	0.8	&	8.3 & 25.4 & 14166 & \textbf{20.79}	&	\textbf{0.4} & \textbf{0.5} & \textbf{48} & 24.26	&	42.5 & 246.4 & 113968 & \textbf{20.79}	\\		
100	&	7	&	0.4	&	315.7 (5) & 43.06\% & 656760 & \textbf{51.11}	&	\textbf{0.6} & 0.8 & \textbf{177} & 52.77	&	599.5 (0) & 44.43\% & 1266772 & \textbf{51.11}	\\		
100	&	7	&	0.6	&	184 (7) & 21.71\% & 362904 & 36.02	&	0.8 & 1.6 & 468 & 38.19	&	599.3 (0) & 28.49\% & 1362716 & 36.19	\\		
100	&	7	&	0.8	&	28.7 & 61.1 & 64042 & \textbf{21.08}	&	\textbf{0.9} & \textbf{1.4} & \textbf{760} & 23.76	&	491.4 (3) & 6.46\% & 2116465 & \textbf{21.08}	\\		
100	&	10	&	0.4	&	165.3 (8) & 39.18\% & 309014 & \textbf{49.29}	&	2.5 & 10.3 & 9239 & 51.98	&	599.4 (0) & 44.7\% & 1410323 & 49.62	\\		
100	&	10	&	0.6	&	106.9 (9) & 15.18\% & 257461 & \textbf{34.4}	&	\textbf{2.8} & \textbf{11.1} & \textbf{9985} & 37.89	&	599.1 (0) & 29.21\% & 1489995 & 34.7	\\		
100	&	10	&	0.8	&	59.3 & 290.2 & 117576 & \textbf{20.61}	&	13.1 & 49.4 & 57599 & 24.83	&	558.7 (1) & 11.48\% & 2253482 & 20.64	\\		
225	&	4	&	0.4	&	86 (9) & 48.2\% & 39407 & \textbf{54.96}	&	\textbf{0.8} & 1.1 & \textbf{45} & 55.77	&	599.4 (0) & 51.97\% & 300625 & 55.09	\\		
225	&	4	&	0.6	&	203.6 (7) & 33.71\% & 107278 & 38.34	&	\textbf{0.8} & \textbf{1} & 49 & 39.42	&	479.5 (2) & 32.79\% & 244773 & \textbf{38.31}	\\		
225	&	4	&	0.8	&	309.9 (5) & 14.89\% & 99639 & \textbf{20.7}	&	\textbf{0.8} & 1.2 & \textbf{84} & 22.13	&	471.1 (3) & 15.78\% & 194755 & \textbf{20.7}	\\		
225	&	7	&	0.4	&	427.2 (3) & 49.21\% & 226670 & 51.59	&	1.3 & 1.8 & 151 & 52.61	&	599.6 (0) & 50.54\% & 316386 & 52.35	\\		
225	&	7	&	0.6	&	484.1 (2) & 32.39\% & 250074 & 35.92	&	16.2 & 143.6 & 10871 & 37.63	&	599.6 (0) & 32.92\% & 332670 & 36.32	\\		
225	&	7	&	0.8	&	551.3 (1) & 13.9\% & 272995 & \textbf{20.44}	&	\textbf{2.6} & \textbf{6.1} & \textbf{1204} & 22.76	&	539.7 (1) & 15.18\% & 303691 & 20.45	\\		
225	&	10	&	0.4	&	599.6 (0) & 48.93\% & 308090 & 50.64	&	\textbf{5.4} & \textbf{24.9} & \textbf{4222} & 51.68	&	599.6 (0) & 51.96\% & 325027 & 52.86	\\		
225	&	10	&	0.6	&	580.6 (1) & 31.75\% & 307517 & 35.6	&	\textbf{10.8} & \textbf{39.7} & \textbf{7537} & 37.07	&	599.9 (0) & 33.36\% & 385794 & 36.12	\\		
225	&	10	&	0.8	&	600 (0) & 14.4\% & 298384 & 20.73	&	\textbf{29.6} & \textbf{141.5} & \textbf{32090} & 23.12	&	599.7 (0) & 16.66\% & 346034 & 20.84	\\		
400	&	4	&	0.4	&	208.6 (7) & 51.54\% & 44423 & \textbf{54.44}	&	\textbf{1.3} & \textbf{1.6} & \textbf{55} & 55.07	&	541 (1) & 52.68\% & 102735 & \textbf{54.44}	\\		
400	&	4	&	0.6	&	372.1 (4) & 37.31\% & 66668 & 37.97	&	\textbf{1.5} & \textbf{1.8} & \textbf{76} & 38.71	&	541 (1) & 35.41\% & 112968 & 38.03	\\		
400	&	4	&	0.8	&	288.8 (6) & 18.01\% & 44492 & \textbf{20.47}	&	\textbf{1.8} & 3.1 & 265 & 21.57	&	599.6 (0) & 17.94\% & 125836 & \textbf{20.47}	\\		
400	&	7	&	0.4	&	599.7 (0) & 54.35\% & 106106 & 53.05	&	19.3 & 169 & 4370 & 52.72	&	599.8 (0) & 55.43\% & 108756 & 54.22	\\		
400	&	7	&	0.6	&	509.6 (2) & 35.93\% & 81877 & 37.05	&	63.4 (9) & 33.32\% & 10416 & 37.48	&	599.6 (0) & 35.83\% & 116993 & 37.01	\\		
400	&	7	&	0.8	&	544.1 (1) & 18.44\% & 105707 & 20.6	&	\textbf{7.2} & \textbf{24.6} & 1999 & 21.87	&	599.7 (0) & 18\% & 122164 & 20.61	\\		
400	&	10	&	0.4	&	599.8 (0) & 53.58\% & 102856 & 52.38	&	\textbf{116.1} & \textbf{242.7} & 83991 & 51.8	&	599.7 (0) & 54.19\% & 122568 & 53.59	\\		
400	&	10	&	0.6	&	529.3 (2) & 36.02\% & 72585 & 36.14	&	\textbf{33.2} & \textbf{115.9} & 17308 & 36.48	&	599.8 (0) & 35.96\% & 117557 & 36.53	\\		
400	&	10	&	0.8	&	599.7 (0) & 17.61\% & 115177 & 20.27	&	155.6 (\textbf{9}) & \textbf{0.04}\% & 104042 & 21.82	&	599.7 (0) & 18.17\% & 79055 & 20.26	\\		
\hline
\end{tabular}
\end{scriptsize}
\end{center}
\caption{Results obtained for the OWAP formulations with the Shortest Path Problem}
%\label{tab_exact_resolution}
\end{table}
%-----------------------------------------------------------------------------------------------------------------------
\newpage
%=======================================================================================================================
%\subsubsection*{SP: Formulations $F_0^{s}$, $F^{s}$, $F_{R1}^{s}$, $F_{R2}^{s}$, $F_{R3}^{s}$}
%-----------------------------------------------------------------------------------------------------------------------
\begin{table}[h!]
\renewcommand{\esp}{@{\hspace{0.15cm}}}
\begin{center}
\begin{scriptsize}
\begin{tabular}{|@{}\esp c \esp c \esp c\esp |c\esp c\esp c\esp c\esp |c\esp c\esp c\esp c\esp |c\esp c\esp c\esp c\esp  @{}|}
\hline
\multicolumn{3}{|c|}{Inst} &
\multicolumn{4}{|c|}{$F_{}^{s}$} &
\multicolumn{4}{|c|}{$F_{R1}^{s}$} &
\multicolumn{4}{|c|}{$F_{R2}^{s}$}   \\
$|V|$ & $p$ & $\alpha$ &
$t(\#)$ & $t^*/gap^*$ & $\#nodes$ & $gap_{LR}$ &
$t(\#)$ & $t^*/gap^*$ & $\#nodes$ & $gap_{LR}$ &
$t(\#)$ & $t^*/gap^*$ & $\#nodes$ & $gap_{LR}$\\
\hline
100	&	4	&	0.4	&	8.4 & 59.9 & 16959 & \textbf{53.72}	&	66.4 (9) & 27.51\% & 204809 & \textbf{53.72}	&	16.9 & 114.5 & 40710 & \textbf{53.72}	\\		
100	&	4	&	0.6	&	31.9 & 213.4 & 113492 & \textbf{37.39}	&	72.3 & 460.2 & 275423 & \textbf{37.39}	&	4.3 & 17 & 6899 & \textbf{37.39}	\\		
100	&	4	&	0.8	&	2.4 & 7.9 & 2871 & \textbf{20.79}	&	8.3 & 23.2 & 11449 & \textbf{20.79}	&	8.5 & 20.9 & 14436 & \textbf{20.79}	\\		
100	&	7	&	0.4	&	121 (8) & 40.66\% & 210086 & \textbf{51.11}	&	124.2 (8) & 11.44\% & 492697 & \textbf{51.11}	&	63.5 (9) & 43.15\% & 102947 & \textbf{51.11}	\\		
100	&	7	&	0.6	&	121 (8) & 22.74\% & 193167 & \textbf{36.01}	&	185.7 (7) & 23.17\% & 477873 & \textbf{36.01}	&	11.2 & 50.4 & 28267 & \textbf{36.01}	\\		
100	&	7	&	0.8	&	20.8 & 125.9 & 36870 & \textbf{21.08}	&	38.8 & 105.4 & 94229 & \textbf{21.08}	&	26.5 & 65.1 & 62967 & \textbf{21.08}	\\		
100	&	10	&	0.4	&	195.2 (7) & 43.64\% & 273035 & \textbf{49.29}	&	75.8 (9) & 22\% & 172053 & \textbf{49.29}	&	20 & 121.2 & 50657 & \textbf{49.29}	\\		
100	&	10	&	0.6	&	178.6 (8) & 24.44\% & 238513 & \textbf{34.4}	&	108.7 (9) & 7.87\% & 322084 & \textbf{34.4}	&	151.2 (8) & 23.02\% & 303828 & \textbf{34.4}	\\		
100	&	10	&	0.8	&	95.3 & 500.7 & 127230 & \textbf{20.61}	&	63.7 & 298.7 & 132170 & \textbf{20.61}	&	57.4 & 201.2 & 132973 & \textbf{20.61}	\\		
225	&	4	&	0.4	&	64.4 (9) & 52.43\% & 29874 & 55	&	72.6 (9) & 49.94\% & 34891 & 55.04	&	198.3 (7) & 52.09\% & 104058 & 55.08	\\		
225	&	4	&	0.6	&	91.5 (9) & 31.77\% & 41747 & \textbf{38.31}	&	140.4 (9) & 30.08\% & 75346 & \textbf{38.31}	&	304 (5) & 31.95\% & 149300 & \textbf{38.31}	\\		
225	&	4	&	0.8	&	243.8 (6) & 14.08\% & 70842 & \textbf{20.7}	&	257.8 (7) & 13.65\% & 123343 & \textbf{20.7}	&	437 (4) & 14.76\% & 178778 & \textbf{20.7}	\\		
225	&	7	&	0.4	&	129 (8) & 49.52\% & 41763 & 51.29	&	313.3 (5) & 50.41\% & 159171 & 51.38	&	489.8 (2) & 49.98\% & 232281 & 51.73	\\		
225	&	7	&	0.6	&	185.6 (7) & 31.25\% & 63146 & 35.83	&	265.6 (6) & 31.59\% & 154086 & \textbf{35.76}	&	472.1 (3) & 32.27\% & 269321 & 36.02	\\		
225	&	7	&	0.8	&	305 (5) & 14.19\% & 105127 & \textbf{20.44}	&	562.7 (3) & 14.08\% & 385137 & \textbf{20.44}	&	577.7 (1) & 13.89\% & 353657 & \textbf{20.44}	\\		
225	&	10	&	0.4	&	317.1 (5) & 49.98\% & 95076 & 50.33	&	599.5 (0) & 49.77\% & 270100 & 51.61	&	599.4 (0) & 51.24\% & 286845 & 51.72	\\		
225	&	10	&	0.6	&	319.5 (5) & 32.14\% & 96370 & \textbf{35.15}	&	565.7 (1) & 31.87\% & 268250 & 35.64	&	599.7 (0) & 31.65\% & 270293 & 35.8	\\		
225	&	10	&	0.8	&	279.6 (6) & 15.16\% & 85419 & 20.81	&	599.3 (0) & 14.71\% & 305531 & 20.9	&	599.5 (0) & 14.01\% & 292841 & 20.79	\\		
400	&	4	&	0.4	&	3.3 & 16.8 & 286 & \textbf{54.44}	&	253.7 (6) & 55.16\% & 36622 & 54.6	&	311 (5) & 54.17\% & 51331 & 54.76	\\		
400	&	4	&	0.6	&	88.5 (9) & 35.79\% & 13806 & \textbf{37.84}	&	355 (5) & 36.67\% & 68752 & 37.93	&	258.4 (6) & 35.91\% & 48491 & 37.92	\\		
400	&	4	&	0.8	&	255.9 (6) & 17.21\% & 49806 & \textbf{20.47}	&	226.7 (7) & 18.11\% & 25218 & \textbf{20.47}	&	374.7 (4) & 17.85\% & 41044 & \textbf{20.47}	\\		
400	&	7	&	0.4	&	76.4 (9) & 50.67\% & 9192 & \textbf{51.85}	&	546 (1) & 54.01\% & 100473 & 52.97	&	497.6 (2) & 52.93\% & 92358 & 52.81	\\		
400	&	7	&	0.6	&	70.9 (9) & 34.59\% & 8711 & \textbf{36.27}	&	599.6 (0) & 36.56\% & 97419 & 36.79	&	490.8 (2) & 36.45\% & 76146 & 36.92	\\		
400	&	7	&	0.8	&	368.8 (4) & 18.29\% & 32614 & 20.41	&	485.1 (2) & 17.99\% & 69436 & 20.71	&	599.3 (0) & 18.38\% & 81188 & 20.57	\\		
400	&	10	&	0.4	&	306.4 (5) & 48.93\% & 24184 & \textbf{50.73}	&	508.1 (2) & 52.44\% & 79421 & 52.03	&	599.6 (0) & 54.16\% & 89357 & 53.69	\\		
400	&	10	&	0.6	&	132.4 (8) & 31.48\% & \textbf{10395} & \textbf{35.01}	&	599.4 (0) & 33.61\% & 117393 & 35.56	&	599.6 (0) & 34.98\% & 98170 & 36.4	\\		
400	&	10	&	0.8	&	\textbf{142.6} (8) & 17.18\% & 12829 & 19.99	&	599.4 (0) & 17.64\% & 80443 & 20.29	&	599.4 (0) & 17.58\% & 109094 & 20.18	\\		
\hline
\end{tabular}
\end{scriptsize}
\end{center}
\caption{Results obtained for the OWAP formulations with the Shortest Path Problem}
%\label{tab_exact_resolution}
\end{table}
\begin{table}[h!]
\renewcommand{\esp}{@{\hspace{0.15cm}}}
\begin{center}
\begin{scriptsize}
\begin{tabular}{|@{}\esp c \esp c \esp c\esp |c\esp c\esp c\esp c\esp   @{}|}
\hline
\multicolumn{3}{|c|}{Inst} &
\multicolumn{4}{|c|}{$F_{R3}^{s}$}  \\
$|V|$ & $p$ & $\alpha$ &
$t(\#)$ & $t^*/gap^*$ & $\#nodes$ & $gap_{LR}$ \\
\hline
100	&	4	&	0.4	&	91.2 & 402.2 & 384235 & \textbf{53.72}	\\
100	&	4	&	0.6	&	53.7 & 296.4 & 220361 & \textbf{37.39}	\\
100	&	4	&	0.8	&	3.3 & 13.2 & 3412 & \textbf{20.79}	\\
100	&	7	&	0.4	&	66.8 (9) & 32.13\% & 184653 & \textbf{51.11}	\\
100	&	7	&	0.6	&	62.4 (9) & 16.05\% & 177939 & \textbf{36.01}	\\
100	&	7	&	0.8	&	4.8 & 25 & 8846 & \textbf{21.08}	\\
100	&	10	&	0.4	&	123.6 (8) & 39.94\% & 219511 & \textbf{49.29}	\\
100	&	10	&	0.6	&	121.9 (8) & 15.68\% & 371819 & \textbf{34.4}	\\
100	&	10	&	0.8	&	\textbf{4.6} & \textbf{9.6} & \textbf{7553} & \textbf{20.61}	\\
225	&	4	&	0.4	&	190.2 (7) & 51.96\% & 91128 & 54.97	\\
225	&	4	&	0.6	&	127.4 (8) & 33.79\% & 71841 & 38.41	\\
225	&	4	&	0.8	&	245.5 (6) & 14.16\% & 78074 & \textbf{20.7}	\\
225	&	7	&	0.4	&	313.7 (5) & 50.43\% & 151042 & 51.36	\\
225	&	7	&	0.6	&	563.6 (1) & 32.75\% & 279849 & 36.07	\\
225	&	7	&	0.8	&	449.7 (3) & 14.28\% & 223068 & \textbf{20.44}	\\
225	&	10	&	0.4	&	293.5 (6) & 49.76\% & 134368 & 50.53	\\
225	&	10	&	0.6	&	489.2 (2) & 32.65\% & 219326 & 35.92	\\
225	&	10	&	0.8	&	338 (5) & 14.78\% & 162308 & 20.73	\\
400	&	4	&	0.4	&	132.2 (8) & 52.61\% & 29642 & 54.47	\\
400	&	4	&	0.6	&	277.2 (6) & 35.79\% & 57596 & 37.98	\\
400	&	4	&	0.8	&	423.8 (3) & 18.1\% & 69373 & \textbf{20.47}	\\
400	&	7	&	0.4	&	357 (5) & 52.98\% & 58912 & 52.73	\\
400	&	7	&	0.6	&	357.2 (5) & 35.86\% & 56379 & 36.64	\\
400	&	7	&	0.8	&	435.6 (3) & 18.38\% & 64149 & \textbf{20.33}	\\
400	&	10	&	0.4	&	559.4 (1) & 53.23\% & 88799 & 52.39	\\
400	&	10	&	0.6	&	446.1 (3) & 34.22\% & 53445 & 35.83	\\
400	&	10	&	0.8	&	386.7 (4) & 17.59\% & 48957 & 20.16	\\

\hline
\end{tabular}
\end{scriptsize}
\end{center}
\caption{Results obtained for the OWAP formulations with the Shortest Path Problem}
%\label{tab_exact_resolution}
\end{table}
%-----------------------------------------------------------------------------------------------------------------------

\newpage
%=======================================================================================================================
%\subsubsection*{PMP: Formulations $F_0^z$, $F^z$, $F_{R1}^z$, $F_{R2}^z$, $F_{R3}^z$}
%-----------------------------------------------------------------------------------------------------------------------
\begin{table}[h!]
\renewcommand{\esp}{@{\hspace{0.15cm}}}
\begin{center}
\begin{scriptsize}
\begin{tabular}{|@{}\esp c \esp c \esp c\esp |c\esp c\esp c\esp c\esp |c\esp c\esp c\esp c\esp |c\esp c\esp c\esp c\esp  @{}|}
\hline
\multicolumn{3}{|c|}{Inst} &
\multicolumn{4}{|c|}{$F_0^z$} &
\multicolumn{4}{|c|}{$F^z$} &
\multicolumn{4}{|c|}{$F^z$$^\prime$}   \\
$|V|$ & $p$ & $\alpha$ &
$t(\#)$ & $t^*/gap^*$ & $\#nodes$ & $gap_{LR}$ &
$t(\#)$ & $t^*/gap^*$ & $\#nodes$ & $gap_{LR}$ &
$t(\#)$ & $t^*/gap^*$ & $\#nodes$ & $gap_{LR}$\\
\hline
100	&	4	&	0.4	&	1.3 & 1.6 & 694 & 100	&	0.7 & 1 & \textbf{109} & \textbf{55.44}	&	0.8 & 1 & 137 & \textbf{55.44}	\\
100	&	4	&	0.6	&	132.1 & 570.7 & 395386 & 100	&	0.7 & 1 & 196 & \textbf{38.98}	&	0.7 & 0.9 & 136 & \textbf{38.98}	\\
100	&	4	&	0.8	&	65.2 & 363.1 & 174806 & 100	&	0.8 & 1 & 235 & \textbf{21.53}	&	0.8 & 1.3 & 244 & \textbf{21.53}	\\
100	&	7	&	0.4	&	164.1 (8) & 15.52\% & 237915 & 100	&	2.4 & 3.5 & 405 & \textbf{52.18}	&	1.6 & 2.1 & 218 & \textbf{52.18}	\\
100	&	7	&	0.6	&	445.1 (5) & 24.15\% & 438889 & 100	&	2.3 & 3 & 452 & \textbf{35.97}	&	1.7 & 2.2 & \textbf{383} & \textbf{35.97}	\\
100	&	7	&	0.8	&	599.3 (0) & 38.43\% & 473919 & 100	&	3.3 & 5.3 & 1302 & \textbf{20.27}	&	2.2 & 3 & 925 & \textbf{20.27}	\\
100	&	10	&	0.4	&	569.9 (1) & 24.29\% & 267816 & 100	&	5.3 & 6.9 & 485 & \textbf{50.66}	&	2.9 & 4.3 & \textbf{284} & \textbf{50.66}	\\
100	&	10	&	0.6	&	599.3 (0) & 13\% & 248655 & 100	&	7.5 & 11.7 & 1162 & \textbf{34.85}	&	4.1 & 5.7 & 1001 & \textbf{34.85}	\\
100	&	10	&	0.8	&	599.4 (0) & 40.14\% & 229928 & 100	&	21.2 & 61.7 & 5799 & \textbf{20.2}	&	14.4 & 48.2 & 6302 & \textbf{20.2}	\\
225	&	4	&	0.4	&	133.1 (8) & 22.67\% & 196711 & 100	&	2.6 & 4.4 & 977 & \textbf{55.09}	&	2.7 & 5.3 & 1339 & \textbf{55.09}	\\
225	&	4	&	0.6	&	169.1 (8) & 24.65\% & 159699 & 100	&	2.4 & 3.7 & 1156 & \textbf{38.57}	&	2.4 & 3.1 & 1148 & \textbf{38.57}	\\
225	&	4	&	0.8	&	151.7 (8) & 12.9\% & 115889 & 100	&	2.6 & 3.8 & 1432 & \textbf{21.09}	&	2.6 & 5.7 & 1064 & \textbf{21.09}	\\
225	&	7	&	0.4	&	599.4 (0) & 26.96\% & 164416 & 100	&	21 & 62.4 & 9718 & \textbf{52.34}	&	15.4 & 54.9 & \textbf{5102} & \textbf{52.34}	\\
225	&	7	&	0.6	&	599.3 (0) & 12.26\% & 148547 & 100	&	34.1 & 126.8 & 12967 & \textbf{36.27}	&	19.3 & 28.5 & 5883 & \textbf{36.27}	\\
225	&	7	&	0.8	&	599.3 (0) & 21.17\% & 126129 & 100	&	47 & 169.8 & 12076 & \textbf{20.32}	&	30 & 117.5 & 10155 & \textbf{20.32}	\\
225	&	10	&	0.4	&	599.3 (0) & 35.34\% & 97866 & 100	&	26.8 & 51.9 & 2100 & \textbf{50.25}	&	19.3 & 33.7 & 2101 & \textbf{50.25}	\\
225	&	10	&	0.6	&	599.7 (0) & 43.06\% & 72292 & 100	&	179.4 (9) & 20.95\% & 16700 & 34.68	&	98.6 & 339.8 & 15018 & \textbf{34.56}	\\
225	&	10	&	0.8	&	599.5 (0) & 48.17\% & 63656 & 100	&	596.9 (1) & 0.61\% & 58769 & 19.64	&	516 (3) & 0.45\% & 70215 & 19.63	\\
400	&	4	&	0.4	&	149.6 (9) & 2.99\% & 70585 & 100	&	10.7 & 30.5 & 2933 & \textbf{55.37}	&	8.7 & 16.9 & 2966 & \textbf{55.37}	\\
400	&	4	&	0.6	&	361.2 (6) & 23.07\% & 144951 & 100	&	10.9 & 15.7 & 4019 & \textbf{39.04}	&	10.8 & 24.4 & 5798 & \textbf{39.04}	\\
400	&	4	&	0.8	&	274.6 (8) & 15.16\% & 96947 & 100	&	12.2 & 19.6 & \textbf{4720} & \textbf{21.03}	&	14.4 & 29.4 & 7329 & \textbf{21.03}	\\
400	&	7	&	0.4	&	599.8 (0) & 10.76\% & 65411 & 100	&	125.3 (9) & 0.1\% & 14707 & \textbf{52.05}	&	84.3 & 409.6 & 12412 & \textbf{52.05}	\\
400	&	7	&	0.6	&	599.9 (0) & 39.05\% & 66520 & 100	&	224.5 & 397.9 & 31091 & \textbf{36.12}	&	190.2 & 467.7 & 34097 & \textbf{36.12}	\\
400	&	7	&	0.8	&	599.9 (0) & 47.19\% & 60028 & 100	&	223.7 & 575.4 & 35131 & \textbf{20.19}	&	205 & 502.3 & 38941 & \textbf{20.19}	\\
400	&	10	&	0.4	&	599.8 (0) & 67.35\% & 44020 & 100	&	301.3 (9) & 0.53\% & 15116 & 50.59	&	168 (9) & 0.3\% & \textbf{12093} & 50.59	\\
400	&	10	&	0.6	&	599.8 (0) & 63.84\% & 36930 & 100	&	511 (2) & 0.55\% & 22643 & 34.58	&	471.8 (5) & 0.26\% & 31918 & \textbf{34.54}	\\
400	&	10	&	0.8	&	599.8 (0) & 67.82\% & 32663 & 100	&	599.9 (0) & 0.63\% & 36012 & 19.59	&	599.8 (0) & 0.42\% & 33206 & 19.53	\\

\hline
\end{tabular}
\end{scriptsize}
\end{center}
\caption{Results obtained for the OWAP formulations with the Perfect Matching Problem}
%\label{tab_exact_resolution}
\end{table}
%-----------------------------------------------------------------------------------------------------------------------

%-----------------------------------------------------------------------------------------------------------------------
\begin{table}[h!]
\renewcommand{\esp}{@{\hspace{0.15cm}}}
\begin{center}
\begin{scriptsize}
\begin{tabular}{|@{}\esp c \esp c \esp c\esp |c\esp c\esp c\esp c\esp |c\esp c\esp c\esp c\esp |c\esp c\esp c\esp c\esp  @{}|}
\hline
\multicolumn{3}{|c|}{Inst} &
\multicolumn{4}{|c|}{$F_{R1}^z$} &
\multicolumn{4}{|c|}{$F_{R2}^z$} &
\multicolumn{4}{|c|}{$F_{R3}^z$}   \\
$|V|$ & $p$ & $\alpha$ &
$t(\#)$ & $t^*/gap^*$ & $\#nodes$ & $gap_{LR}$ &
$t(\#)$ & $t^*/gap^*$ & $\#nodes$ & $gap_{LR}$ &
$t(\#)$ & $t^*/gap^*$ & $\#nodes$ & $gap_{LR}$\\
\hline
100	&	4	&	0.4	&	0.6 & \textbf{0.7} & 186 & \textbf{55.44}	&	0.7 & 0.8 & 138 & 55.98	&	212.2 (7) & 19.55\% & 628046 & \textbf{55.44}	\\
100	&	4	&	0.6	&	\textbf{0.6} & 0.7 & 152 & \textbf{38.98}	&	\textbf{0.6} & \textbf{0.6} & 149 & 39.71	&	130.2 (8) & 13.38\% & 354878 & \textbf{38.98}	\\
100	&	4	&	0.8	&	\textbf{0.6} & \textbf{0.7} & 302 & \textbf{21.53}	&	\textbf{0.6} & 0.8 & 234 & 22.47	&	164.1 (8) & 5.69\% & 487138 & \textbf{21.53}	\\
100	&	7	&	0.4	&	1 & 1.3 & 236 & \textbf{52.18}	&	1.2 & 1.5 & 1020 & 53.23	&	599.3 (0) & 46.18\% & 590657 & 52.19	\\
100	&	7	&	0.6	&	\textbf{1.1} & 1.4 & 480 & \textbf{35.97}	&	60.9 (9) & 13.62\% & 99261 & 37.43	&	599.3 (0) & 31.27\% & 662540 & \textbf{35.97}	\\
100	&	7	&	0.8	&	1.4 & 2 & 965 & \textbf{20.27}	&	1.6 & 3.1 & 1419 & 22.02	&	599.3 (0) & 16.8\% & 880878 & 20.3	\\
100	&	10	&	0.4	&	\textbf{1.5} & 1.9 & 299 & \textbf{50.66}	&	7.4 & 28.9 & 22875 & 51.75	&	599.3 (0) & 46.74\% & 339601 & 50.83	\\
100	&	10	&	0.6	&	\textbf{1.9} & \textbf{2.6} & 963 & \textbf{34.85}	&	2.6 & 4.5 & 2032 & 36.29	&	599.3 (0) & 32.35\% & 392747 & 35.35	\\
100	&	10	&	0.8	&	6 & 19.4 & 6329 & \textbf{20.2}	&	11.3 & 36.2 & 13433 & 21.97	&	599.3 (0) & 18.05\% & 730420 & 20.56	\\
225	&	4	&	0.4	&	2.1 & 4.4 & 1188 & \textbf{55.09}	&	\textbf{1.9} & 3 & 1003 & 55.48	&	286.9 (6) & 44.05\% & 212277 & 55.3	\\
225	&	4	&	0.6	&	\textbf{1.7} & 2.9 & 1236 & \textbf{38.57}	&	\textbf{1.7} & 2.5 & 1262 & 39.1	&	193.8 (8) & 22.64\% & 157652 & \textbf{38.57}	\\
225	&	4	&	0.8	&	\textbf{1.9} & 3.2 & 1101 & \textbf{21.09}	&	2.1 & 3.8 & 1504 & 21.77	&	436.6 (3) & 13.99\% & 357325 & \textbf{21.09}	\\
225	&	7	&	0.4	&	\textbf{7.1} & \textbf{22.8} & 9208 & \textbf{52.34}	&	17.9 & 96.4 & 21252 & 52.73	&	599.3 (0) & 48.08\% & 176693 & 52.59	\\
225	&	7	&	0.6	&	10 & 16 & 6038 & \textbf{36.27}	&	12.3 & 33.5 & 8888 & 36.79	&	599.3 (0) & 32.95\% & 191613 & 36.45	\\
225	&	7	&	0.8	&	17.2 & 62.5 & 10491 & \textbf{20.32}	&	17.6 & 67.6 & 11652 & 20.97	&	599.3 (0) & 17.7\% & 232228 & 20.43	\\
225	&	10	&	0.4	&	7.5 & 13.2 & 2136 & \textbf{50.25}	&	18.7 & 59.6 & 15457 & 50.65	&	599.3 (0) & 48.46\% & 102767 & 50.7	\\
225	&	10	&	0.6	&	32.4 & 123.2 & 15537 & \textbf{34.56}	&	45.3 & 101.4 & 27800 & 35.09	&	599.3 (0) & 33.15\% & 139120 & 35.01	\\
225	&	10	&	0.8	&	295 (8) & 0.32\% & 114029 & \textbf{19.62}	&	402.6 (6) & 0.31\% & 169450 & 20.26	&	599.3 (0) & 18.82\% & 164745 & 19.95	\\
400	&	4	&	0.4	&	7.3 & 22.3 & 3345 & \textbf{55.37}	&	6.5 & 17 & 3311 & 55.5	&	455.6 (5) & 45.64\% & 119265 & 55.4	\\
400	&	4	&	0.6	&	\textbf{6.7} & \textbf{11.9} & 4103 & \textbf{39.04}	&	8.6 & 27.3 & 6816 & 39.2	&	360.4 (7) & 31.89\% & 99906 & 39.06	\\
400	&	4	&	0.8	&	9 & 22.1 & 5397 & \textbf{21.03}	&	12.1 & 25.5 & 8260 & 21.24	&	512.7 (2) & 16.34\% & 187051 & \textbf{21.03}	\\
400	&	7	&	0.4	&	\textbf{34.4} & 144.4 & \textbf{10464} & \textbf{52.05}	&	53.4 & \textbf{111.5} & 23341 & 52.26	&	599.9 (0) & 48.8\% & 85058 & 52.56	\\
400	&	7	&	0.6	&	83.4 & 250.9 & 27604 & \textbf{36.12}	&	112.2 & 240.7 & 42546 & 36.4	&	599.9 (0) & 33.56\% & 93893 & 36.6	\\
400	&	7	&	0.8	&	\textbf{84.4} & 187.6 & \textbf{28762} & \textbf{20.19}	&	224.3 (9) & 0.44\% & 146178 & 20.55	&	599.9 (0) & 18.1\% & 118886 & 20.38	\\
400	&	10	&	0.4	&	68.4 & \textbf{197.4} & 13777 & \textbf{50.58}	&	278.2 (7) & 0.45\% & 104776 & 50.79	&	599.8 (0) & 49.31\% & 53902 & 51.49	\\
400	&	10	&	0.6	&	289.4 (9) & 0.11\% & 61886 & \textbf{34.54}	&	338.6 (8) & 0.19\% & 81662 & 34.8	&	599.9 (0) & 33.99\% & 59623 & 35.36	\\
400	&	10	&	0.8	&	583.5 (\textbf{1}) & 0.42\% & 97022 & \textbf{19.5}	&	599.9 (0) & 0.48\% & 108421 & 19.85	&	599.9 (0) & 18.46\% & 66410 & 19.96	\\

\hline
\end{tabular}
\end{scriptsize}
\end{center}
\caption{Results obtained for the OWAP formulations with the Perfect Matching Problem}
%\label{tab_exact_resolution}
\end{table}
%-----------------------------------------------------------------------------------------------------------------------
%=======================================================================================================================
\newpage
%=======================================================================================================================
%\subsubsection*{PMP: Formulations $F_0^{zy}$, $F^{zy}$, $F_{R1}^{zy}$, $F_{R2}^{zy}$, $F_{R3}^{zy}$}
%-----------------------------------------------------------------------------------------------------------------------
\begin{table}[h!]
\renewcommand{\esp}{@{\hspace{0.15cm}}}
\begin{center}
\begin{scriptsize}
\begin{tabular}{|@{}\esp c \esp c \esp c\esp |c\esp c\esp c\esp c\esp |c\esp c\esp c\esp c\esp |c\esp c\esp c\esp c\esp  @{}|}
\hline
\multicolumn{3}{|c|}{Inst} &
\multicolumn{4}{|c|}{$F_{0}^{zy}$} &
\multicolumn{4}{|c|}{$F_{}^{zy}$} &
\multicolumn{4}{|c|}{$F_{}^{zy}$$^\prime$}   \\
$|V|$ & $p$ & $\alpha$ &
$t(\#)$ & $t^*/gap^*$ & $\#nodes$ & $gap_{LR}$ &
$t(\#)$ & $t^*/gap^*$ & $\#nodes$ & $gap_{LR}$ &
$t(\#)$ & $t^*/gap^*$ & $\#nodes$ & $gap_{LR}$\\
\hline
100	&	4	&	0.4	&	1.5 & 2.7 & 1346 & 100	&	0.8 & 1 & 162 & \textbf{55.44}	&	0.7 & 0.9 & \textbf{109} & \textbf{55.44}	\\
100	&	4	&	0.6	&	19.2 & 93.3 & 48335 & 100	&	0.7 & 1 & 121 & \textbf{38.98}	&	0.7 & 0.8 & 151 & \textbf{38.98}	\\
100	&	4	&	0.8	&	52.6 & 440.5 & 133230 & 100	&	0.8 & 1 & 229 & \textbf{21.53}	&	0.7 & 1 & 197 & \textbf{21.53}	\\
100	&	7	&	0.4	&	70.2 & 271.4 & 74398 & 100	&	2.4 & 3.5 & 326 & \textbf{52.18}	&	1.9 & 2.7 & 230 & \textbf{52.18}	\\
100	&	7	&	0.6	&	323.9 (7) & 2.62\% & 311151 & 100	&	2.5 & 3.4 & 892 & \textbf{35.97}	&	1.9 & 2.5 & 403 & \textbf{35.97}	\\
100	&	7	&	0.8	&	557.8 (1) & 35.79\% & 458556 & 100	&	3.8 & 6.7 & 1540 & \textbf{20.27}	&	2.7 & 4 & 892 & \textbf{20.27}	\\
100	&	10	&	0.4	&	557 (1) & 24.67\% & 284261 & 100	&	4.8 & 6 & 448 & \textbf{50.66}	&	4.6 & 6.9 & 417 & \textbf{50.66}	\\
100	&	10	&	0.6	&	599.3 (0) & 11.6\% & 240287 & 100	&	6.7 & 8.9 & 1162 & \textbf{34.85}	&	5.3 & 7.6 & 837 & \textbf{34.85}	\\
100	&	10	&	0.8	&	599.4 (0) & 41\% & 223345 & 100	&	25.3 & 85.7 & 7395 & \textbf{20.2}	&	18.1 & 45.9 & 6194 & \textbf{20.2}	\\
225	&	4	&	0.4	&	154.4 (8) & 12.71\% & 121716 & 100	&	2.7 & 6.6 & 1576 & \textbf{55.09}	&	2.4 & 3.4 & 965 & \textbf{55.09}	\\
225	&	4	&	0.6	&	250.2 (6) & 12.41\% & 208907 & 100	&	2.3 & 3.4 & 1187 & \textbf{38.57}	&	2.3 & 3.9 & 1118 & \textbf{38.57}	\\
225	&	4	&	0.8	&	120 (9) & 9.23\% & 94574 & 100	&	2.5 & 4.6 & 1317 & \textbf{21.09}	&	2.3 & 4.2 & 1056 & \textbf{21.09}	\\
225	&	7	&	0.4	&	532.5 (2) & 3.39\% & 148297 & 100	&	19.1 & 58.9 & 7577 & \textbf{52.34}	&	22.6 & 96.9 & 6535 & \textbf{52.34}	\\
225	&	7	&	0.6	&	599.3 (0) & 40.02\% & 169446 & 100	&	23.7 & 43.5 & 6293 & \textbf{36.27}	&	26.5 & 57.7 & 6779 & \textbf{36.27}	\\
225	&	7	&	0.8	&	599.3 (0) & 46.4\% & 128394 & 100	&	43.4 & 167.1 & 12348 & \textbf{20.32}	&	38 & 151 & 9632 & \textbf{20.32}	\\
225	&	10	&	0.4	&	599.5 (0) & 34.2\% & 91387 & 100	&	26.5 & 61.8 & 2006 & \textbf{50.25}	&	24 & 37.9 & 1866 & \textbf{50.25}	\\
225	&	10	&	0.6	&	599.6 (0) & 39.98\% & 77386 & 100	&	130 & 362.9 & 18614 & \textbf{34.56}	&	123.3 & 415 & 13898 & \textbf{34.56}	\\
225	&	10	&	0.8	&	599.8 (0) & 57.06\% & 63465 & 100	&	588.8 (1) & 0.51\% & \textbf{53502} & 19.63	&	576.1 (2) & 0.52\% & 56654 & 19.63	\\
400	&	4	&	0.4	&	285.4 (6) & 5.64\% & 161874 & 100	&	12.2 & 37.6 & 3617 & \textbf{55.37}	&	9.1 & 15.4 & 3329 & \textbf{55.37}	\\
400	&	4	&	0.6	&	231.2 (8) & 11.41\% & 94521 & 100	&	14.4 & 29.8 & 5798 & \textbf{39.04}	&	10.3 & 23 & 4340 & \textbf{39.04}	\\
400	&	4	&	0.8	&	384.6 (7) & 16.27\% & 129285 & 100	&	13.9 & 33.8 & 5320 & \textbf{21.03}	&	12.5 & 25.2 & 5468 & \textbf{21.03}	\\
400	&	7	&	0.4	&	599.9 (0) & 31.44\% & 64070 & 100	&	128.6 (9) & 0.23\% & 14895 & \textbf{52.05}	&	104.5 & 383 & 13188 & \textbf{52.05}	\\
400	&	7	&	0.6	&	599.9 (0) & 36.91\% & 63436 & 100	&	258.9 & 459.9 & 37140 & \textbf{36.12}	&	205 (9) & 0.05\% & 33027 & \textbf{36.12}	\\
400	&	7	&	0.8	&	599.8 (0) & 48.62\% & 67620 & 100	&	267.4 & 526.6 & 39724 & \textbf{20.19}	&	221.6 & 447.7 & 33662 & \textbf{20.19}	\\
400	&	10	&	0.4	&	599.9 (0) & 37.68\% & 40006 & 100	&	253.7 (9) & 0.26\% & 13493 & 50.59	&	252.2 (9) & 0.31\% & 13052 & 50.59	\\
400	&	10	&	0.6	&	599.9 (0) & 42.48\% & 42371 & 100	&	558 (1) & 0.39\% & 22812 & 34.57	&	488.3 (4) & 0.35\% & \textbf{22428} & 34.56	\\
400	&	10	&	0.8	&	599.8 (0) & 63\% & 37416 & 100	&	599.8 (0) & 0.59\% & \textbf{21413} & 19.58	&	599.9 (0) & 0.49\% & 24123 & 19.53	\\

\hline
\end{tabular}
\end{scriptsize}
\end{center}
\caption{Results obtained for the OWAP formulations with the Perfect Matching Problem}
%\label{tab_exact_resolution}
\end{table}
%-----------------------------------------------------------------------------------------------------------------------

%-----------------------------------------------------------------------------------------------------------------------
\begin{table}[h!]
\renewcommand{\esp}{@{\hspace{0.15cm}}}
\begin{center}
\begin{scriptsize}
\begin{tabular}{|@{}\esp c \esp c \esp c\esp |c\esp c\esp c\esp c\esp |c\esp c\esp c\esp c\esp |c\esp c\esp c\esp c\esp  @{}|}
\hline
\multicolumn{3}{|c|}{Inst} &
\multicolumn{4}{|c|}{$F_{R1}^{zy}$} &
\multicolumn{4}{|c|}{$F_{R2}^{zy}$} &
\multicolumn{4}{|c|}{$F_{R3}^{zy}$}   \\
$|V|$ & $p$ & $\alpha$ &
$t(\#)$ & $t^*/gap^*$ & $\#nodes$ & $gap_{LR}$ &
$t(\#)$ & $t^*/gap^*$ & $\#nodes$ & $gap_{LR}$ &
$t(\#)$ & $t^*/gap^*$ & $\#nodes$ & $gap_{LR}$\\
\hline
100	&	4	&	0.4	&	0.6 & \textbf{0.7} & 139 & \textbf{55.44}	&	0.7 & 0.8 & 133 & 55.98	&	189.7 (7) & 22.3\% & 505255 & 55.47	\\
100	&	4	&	0.6	&	\textbf{0.6} & 0.7 & 140 & \textbf{38.98}	&	\textbf{0.6} & 0.7 & 147 & 39.71	&	29.5 & 160.6 & 91753 & \textbf{38.98}	\\
100	&	4	&	0.8	&	0.7 & 0.8 & 329 & \textbf{21.53}	&	\textbf{0.6} & 0.8 & 285 & 22.47	&	238.3 (7) & 5.76\% & 703988 & \textbf{21.53}	\\
100	&	7	&	0.4	&	1 & \textbf{1.2} & 256 & \textbf{52.18}	&	1.2 & 1.8 & 959 & 53.23	&	599.3 (0) & 46.34\% & 605536 & 52.21	\\
100	&	7	&	0.6	&	1.2 & 1.8 & 591 & \textbf{35.97}	&	\textbf{1.1} & \textbf{1.2} & 442 & 37.38	&	599.2 (0) & 31.25\% & 677159 & 36.12	\\
100	&	7	&	0.8	&	1.5 & 2.3 & 1008 & \textbf{20.27}	&	1.5 & 2.4 & 1287 & 22.02	&	599.2 (0) & 16.78\% & 911147 & \textbf{20.27}	\\
100	&	10	&	0.4	&	1.7 & 4.1 & 580 & \textbf{50.66}	&	12.3 & 54.3 & 36234 & 51.75	&	599.2 (0) & 46.81\% & 352950 & 50.83	\\
100	&	10	&	0.6	&	\textbf{1.9} & 2.9 & 985 & \textbf{34.85}	&	3.9 & 16.7 & 5175 & 36.29	&	599.2 (0) & 31.96\% & 407864 & 35.23	\\
100	&	10	&	0.8	&	5.6 & 17.6 & \textbf{5364} & \textbf{20.2}	&	29.4 & 183 & 42168 & 21.97	&	599.2 (0) & 18.22\% & 844424 & 20.77	\\
225	&	4	&	0.4	&	2 & \textbf{2.9} & 990 & \textbf{55.09}	&	\textbf{1.9} & 3.3 & 991 & 55.48	&	433.8 (4) & 36.25\% & 332362 & 55.11	\\
225	&	4	&	0.6	&	\textbf{1.7} & 2.5 & 1239 & \textbf{38.57}	&	4.6 & 29.6 & 7582 & 39.1	&	112.7 (9) & 19.79\% & 81949 & \textbf{38.57}	\\
225	&	4	&	0.8	&	\textbf{1.9} & 3.7 & 1240 & \textbf{21.09}	&	\textbf{1.9} & 3.2 & 1408 & 21.77	&	384.2 (4) & 10.46\% & 335049 & \textbf{21.09}	\\
225	&	7	&	0.4	&	8.4 & 36 & 5617 & \textbf{52.34}	&	15.1 & 48.9 & 16006 & 52.73	&	599.3 (0) & 47.91\% & 180294 & 52.49	\\
225	&	7	&	0.6	&	9.7 & 18.3 & 6206 & \textbf{36.27}	&	9.8 & 19.9 & 6317 & 36.79	&	599.2 (0) & 32.65\% & 201172 & 36.52	\\
225	&	7	&	0.8	&	17.1 & \textbf{48.7} & 10746 & \textbf{20.32}	&	28.6 & 172.6 & 20230 & 20.97	&	599.3 (0) & 17.63\% & 217743 & 20.38	\\
225	&	10	&	0.4	&	7.4 & 12.4 & 2464 & \textbf{50.25}	&	116.1 (9) & 0.11\% & 126039 & 50.65	&	599.2 (0) & 48.35\% & 105064 & 50.73	\\
225	&	10	&	0.6	&	33.9 & 90.1 & 13763 & \textbf{34.56}	&	70.4 & 218.1 & 50035 & 35.09	&	599.3 (0) & 33.24\% & 127851 & 35.03	\\
225	&	10	&	0.8	&	338.7 (7) & 12.07\% & 130079 & 19.7	&	367.3 (7) & 0.33\% & 198458 & 20.27	&	599.4 (0) & 18.61\% & 171503 & 20.08	\\
400	&	4	&	0.4	&	6.3 & 15.5 & 2546 & \textbf{55.37}	&	6.8 & 17.7 & 3494 & 55.5	&	382 (7) & 47.5\% & 94550 & 55.42	\\
400	&	4	&	0.6	&	7.5 & 16.7 & 4044 & \textbf{39.04}	&	8.6 & 14.1 & 7036 & 39.2	&	422 (5) & 32.69\% & 121143 & 39.08	\\
400	&	4	&	0.8	&	11.4 & 44.9 & 6263 & \textbf{21.03}	&	14.1 & 31 & 10928 & 21.24	&	475.2 (3) & 15.57\% & 167062 & \textbf{21.03}	\\
400	&	7	&	0.4	&	48.9 & 257 & 15696 & \textbf{52.05}	&	65.9 & 155.3 & 29514 & 52.26	&	599.8 (0) & 49.42\% & 81090 & 52.49	\\
400	&	7	&	0.6	&	\textbf{74.5} & 209.5 & \textbf{26944} & \textbf{36.12}	&	115.3 & 317.5 & 50818 & 36.4	&	599.9 (0) & 33.26\% & 93022 & 36.48	\\
400	&	7	&	0.8	&	98.2 & 182.5 & 35369 & \textbf{20.19}	&	155.6 & 320.5 & 67682 & 20.55	&	599.9 (0) & 17.91\% & 109012 & 20.35	\\
400	&	10	&	0.4	&	86.7 & 387.2 & 17514 & \textbf{50.58}	&	467.1 (3) & 0.26\% & 178559 & 50.78	&	599.8 (0) & 49.32\% & 52303 & 51.48	\\
400	&	10	&	0.6	&	335 (9) & 0.24\% & 69457 & \textbf{34.54}	&	328.1 (7) & 0.31\% & 80674 & 34.81	&	599.8 (0) & 33.83\% & 60762 & 35.43	\\
400	&	10	&	0.8	&	599 (\textbf{1}) & 0.4\% & 93171 & \textbf{19.5}	&	583.1 (\textbf{1}) & 0.64\% & 108405 & 19.86	&	599.8 (0) & 18.49\% & 65760 & 19.98	\\

\hline
\end{tabular}
\end{scriptsize}
\end{center}
\caption{Results obtained for the OWAP formulations with the Perfect Matching Problem}
%\label{tab_exact_resolution}
\end{table}
%-----------------------------------------------------------------------------------------------------------------------
%=======================================================================================================================
\newpage
%=======================================================================================================================
%\subsubsection*{PMP: Formulations $F_0^{s}$, $F^{s}$, $F_{R1}^{s}$, $F_{R2}^{s}$, $F_{R3}^{s}$}
%-----------------------------------------------------------------------------------------------------------------------
\begin{table}[h!]
\renewcommand{\esp}{@{\hspace{0.15cm}}}
\begin{center}
\begin{scriptsize}
\begin{tabular}{|@{}\esp c \esp c \esp c\esp |c\esp c\esp c\esp c\esp |c\esp c\esp c\esp c\esp |c\esp c\esp c\esp c\esp  @{}|}
\hline
\multicolumn{3}{|c|}{Inst} &
\multicolumn{4}{|c|}{$F_{}^{s}$} &
\multicolumn{4}{|c|}{$F_{R1}^{s}$} &
\multicolumn{4}{|c|}{$F_{R2}^{s}$}   \\
$|V|$ & $p$ & $\alpha$ &
$t(\#)$ & $t^*/gap^*$ & $\#nodes$ & $gap_{LR}$ &
$t(\#)$ & $t^*/gap^*$ & $\#nodes$ & $gap_{LR}$ &
$t(\#)$ & $t^*/gap^*$ & $\#nodes$ & $gap_{LR}$\\
\hline
100	&	4	&	0.4	&	0.7 & 1 & 126 & \textbf{55.44}	&	\textbf{0.5} & \textbf{0.7} & 147 & \textbf{55.44}	&	0.6 & 0.8 & 154 & \textbf{55.44}	\\
100	&	4	&	0.6	&	0.7 & 0.8 & 156 & \textbf{38.98}	&	\textbf{0.6} & 0.8 & 175 & \textbf{38.98}	&	\textbf{0.6} & 0.8 & 148 & \textbf{38.98}	\\
100	&	4	&	0.8	&	0.7 & 1 & 223 & \textbf{21.53}	&	\textbf{0.6} & \textbf{0.7} & 157 & \textbf{21.53}	&	\textbf{0.6} & \textbf{0.7} & 254 & \textbf{21.53}	\\
100	&	7	&	0.4	&	1.7 & 2.8 & 223 & \textbf{52.18}	&	1 & \textbf{1.2} & \textbf{205} & \textbf{52.18}	&	\textbf{0.9} & \textbf{1.2} & 260 & \textbf{52.18}	\\
100	&	7	&	0.6	&	2.1 & 2.8 & 464 & \textbf{35.97}	&	1.2 & 1.6 & 529 & \textbf{35.97}	&	\textbf{1.1} & 1.5 & 407 & \textbf{35.97}	\\
100	&	7	&	0.8	&	2.9 & 4.5 & 931 & \textbf{20.27}	&	\textbf{1.3} & \textbf{1.8} & 1075 & \textbf{20.27}	&	1.5 & 2.2 & 1000 & \textbf{20.27}	\\
100	&	10	&	0.4	&	63.8 (9) & 17.02\% & 17496 & 50.7	&	\textbf{1.5} & 1.9 & 333 & \textbf{50.66}	&	\textbf{1.5} & \textbf{1.7} & 363 & \textbf{50.66}	\\
100	&	10	&	0.6	&	5.8 & 9.3 & 870 & \textbf{34.85}	&	2 & 2.8 & 922 & \textbf{34.85}	&	\textbf{1.9} & 2.7 & 763 & \textbf{34.85}	\\
100	&	10	&	0.8	&	21.6 & 81.3 & 6559 & \textbf{20.2}	&	6.1 & 19.8 & 7018 & \textbf{20.2}	&	6 & 17.7 & 5723 & \textbf{20.2}	\\
225	&	4	&	0.4	&	2.4 & 4.2 & 1254 & \textbf{55.09}	&	\textbf{1.9} & 4.1 & 1095 & \textbf{55.09}	&	2.1 & 3.6 & \textbf{878} & \textbf{55.09}	\\
225	&	4	&	0.6	&	2.2 & 2.9 & 1289 & \textbf{38.57}	&	\textbf{1.7} & \textbf{2.2} & \textbf{982} & \textbf{38.57}	&	2 & 3 & 1544 & \textbf{38.57}	\\
225	&	4	&	0.8	&	2.4 & 4.8 & \textbf{819} & \textbf{21.09}	&	2 & 3.6 & 1221 & \textbf{21.09}	&	2 & 3.7 & 1467 & \textbf{21.09}	\\
225	&	7	&	0.4	&	24 & 103.1 & 9052 & \textbf{52.34}	&	8.7 & 29.3 & 6308 & \textbf{52.34}	&	9 & 30.7 & 8426 & \textbf{52.34}	\\
225	&	7	&	0.6	&	21.7 & 42.5 & 10157 & \textbf{36.27}	&	\textbf{8.8} & \textbf{15.9} & \textbf{5432} & \textbf{36.27}	&	9.6 & 20.4 & 7793 & \textbf{36.27}	\\
225	&	7	&	0.8	&	36.5 & 142.9 & 12890 & \textbf{20.32}	&	\textbf{14.7} & 50.1 & \textbf{9525} & \textbf{20.32}	&	17.9 & 78.1 & 10571 & \textbf{20.32}	\\
225	&	10	&	0.4	&	21.3 & 42 & \textbf{1510} & \textbf{50.25}	&	7.8 & 15.5 & 2265 & \textbf{50.25}	&	\textbf{6.7} & \textbf{12.1} & 1662 & \textbf{50.25}	\\
225	&	10	&	0.6	&	114.1 & 260.6 & 19000 & \textbf{34.56}	&	31.5 & 70.1 & 15465 & \textbf{34.56}	&	32.7 & 84.9 & 13290 & \textbf{34.56}	\\
225	&	10	&	0.8	&	565.3 (2) & 0.5\% & 55981 & 19.63	&	344.7 (8) & 0.33\% & 133025 & \textbf{19.62}	&	313.3 (8) & 0.29\% & 116767 & 19.63	\\
400	&	4	&	0.4	&	8.9 & 25.7 & \textbf{2510} & \textbf{55.37}	&	\textbf{6.1} & \textbf{9.6} & 2777 & \textbf{55.37}	&	\textbf{6.1} & 12.6 & 2993 & \textbf{55.37}	\\
400	&	4	&	0.6	&	9.8 & 29.3 & \textbf{3726} & \textbf{39.04}	&	8.7 & 25.3 & 6589 & \textbf{39.04}	&	8.7 & 23.5 & 5234 & \textbf{39.04}	\\
400	&	4	&	0.8	&	13.4 & 25.5 & 5008 & \textbf{21.03}	&	9.2 & 19.4 & 5194 & \textbf{21.03}	&	\textbf{8.8} & \textbf{15.5} & 5075 & \textbf{21.03}	\\
400	&	7	&	0.4	&	106 (9) & 0.53\% & 11087 & 52.06	&	37.7 & 218.2 & 11164 & \textbf{52.05}	&	48.3 & 265.6 & 15438 & \textbf{52.05}	\\
400	&	7	&	0.6	&	228.7 (9) & 0.11\% & 38373 & \textbf{36.12}	&	78.7 & \textbf{185.1} & 28692 & \textbf{36.12}	&	102.8 & 331 & 54127 & \textbf{36.12}	\\
400	&	7	&	0.8	&	222.6 (9) & 0.1\% & 37955 & \textbf{20.19}	&	92.6 & 206.4 & 34328 & \textbf{20.19}	&	89.4 & \textbf{175} & 35596 & \textbf{20.19}	\\
400	&	10	&	0.4	&	275.4 (9) & 0.35\% & 13254 & 50.59	&	91.9 & 407.6 & 19024 & \textbf{50.58}	&	93.8 & 404.1 & 20114 & \textbf{50.58}	\\
400	&	10	&	0.6	&	515.7 (2) & 0.43\% & 25175 & 34.56	&	285.7 & \textbf{563.5} & 59428 & \textbf{34.54}	&	\textbf{253.3} (9) & 0.08\% & 50438 & \textbf{34.54}	\\
400	&	10	&	0.8	&	599.9 (0) & 0.49\% & 24012 & 19.52	&	\textbf{577} (\textbf{1}) & 0.43\% & 97258 & 19.52	&	599.9 (0) & \textbf{0.35}\% & 98257 & \textbf{19.5}	\\

\hline
\end{tabular}
\end{scriptsize}
\end{center}
\caption{Results obtained for the OWAP formulations with the Perfect Matching Problem}
%\label{tab_exact_resolution}
\end{table}
\begin{table}[h!]
\renewcommand{\esp}{@{\hspace{0.15cm}}}
\begin{center}
\begin{scriptsize}
\begin{tabular}{|@{}\esp c \esp c \esp c\esp |c\esp c\esp c\esp c\esp   @{}|}
\hline
\multicolumn{3}{|c|}{Inst} &
\multicolumn{4}{|c|}{$F_{R3}^{s}$} \\
$|V|$ & $p$ & $\alpha$ &
$t(\#)$ & $t^*/gap^*$ & $\#nodes$ & $gap_{LR}$\\
\hline
100	&	4	&	0.4	&	0.8 & 0.9 & 176 & \textbf{55.44}	\\
100	&	4	&	0.6	&	0.7 & 0.8 & \textbf{119} & \textbf{38.98}	\\
100	&	4	&	0.8	&	0.7 & 0.9 & \textbf{137} & \textbf{21.53}	\\
100	&	7	&	0.4	&	1.3 & 1.4 & 270 & \textbf{52.18}	\\
100	&	7	&	0.6	&	1.4 & 1.9 & 463 & \textbf{35.97}	\\
100	&	7	&	0.8	&	1.7 & 2.3 & \textbf{813} & \textbf{20.27}	\\
100	&	10	&	0.4	&	1.8 & 2.2 & 342 & \textbf{50.66}	\\
100	&	10	&	0.6	&	2 & 2.7 & \textbf{748} & \textbf{34.85}	\\
100	&	10	&	0.8	&	\textbf{5} & \textbf{11.9} & 5754 & \textbf{20.2}	\\
225	&	4	&	0.4	&	2.4 & 6.3 & 1118 & \textbf{55.09}	\\
225	&	4	&	0.6	&	2.3 & 3.1 & 1159 & \textbf{38.57}	\\
225	&	4	&	0.8	&	2.2 & \textbf{3.1} & 822 & \textbf{21.09}	\\
225	&	7	&	0.4	&	9.9 & 31 & 6591 & \textbf{52.34}	\\
225	&	7	&	0.6	&	13.6 & 31.1 & 8822 & \textbf{36.27}	\\
225	&	7	&	0.8	&	18.6 & 52.3 & 10352 & \textbf{20.32}	\\
225	&	10	&	0.4	&	7.6 & 15 & 3537 & \textbf{50.25}	\\
225	&	10	&	0.6	&	\textbf{30.2} & \textbf{62.7} & \textbf{12952} & \textbf{34.56}	\\
225	&	10	&	0.8	&	\textbf{244.6} & \textbf{533.1} & 102387 & \textbf{19.62}	\\
400	&	4	&	0.4	&	7.9 & 15.8 & 3787 & \textbf{55.37}	\\
400	&	4	&	0.6	&	8.9 & 17.2 & 4354 & \textbf{39.04}	\\
400	&	4	&	0.8	&	11.3 & 27.7 & 5626 & \textbf{21.03}	\\
400	&	7	&	0.4	&	41.7 & 116 & 22421 & \textbf{52.05}	\\
400	&	7	&	0.6	&	91.4 & 312.2 & 30629 & \textbf{36.12}	\\
400	&	7	&	0.8	&	106.4 & 249.3 & 35269 & \textbf{20.19}	\\
400	&	10	&	0.4	&	\textbf{51.8} & 204.4 & 13716 & \textbf{50.58}	\\
400	&	10	&	0.6	&	268.8 (9) & 0.12\% & 76615 & \textbf{34.54}	\\
400	&	10	&	0.8	&	577.6 (\textbf{1}) & 0.36\% & 114216 & \textbf{19.5}	\\

\hline
\end{tabular}
\end{scriptsize}
\end{center}
\caption{Results obtained for the OWAP formulations with the Perfect Matching Problem}
%\label{tab_exact_resolution}
\end{table}
%-----------------------------------------------------------------------------------------------------------------------

\newpage
%=======================================================================================================================
%\subsubsection*{SP con DV:  $F_{R2}^{z}$ + (25.1), (25.2), (26.1), (26.2), (29.1)}
%-----------------------------------------------------------------------------------------------------------------------
\begin{table}[h!]
\renewcommand{\esp}{@{\hspace{0.15cm}}}
\begin{center}
\begin{scriptsize}
\begin{tabular}{|@{}\esp c \esp c \esp c\esp |c\esp c\esp c\esp c\esp |c\esp c\esp c\esp c\esp |c\esp c\esp c\esp c\esp  @{}|}
\hline
\multicolumn{3}{|c|}{Inst} &
\multicolumn{4}{|c|}{$F_{R2}^{z}$} &
\multicolumn{4}{|c|}{$(\ref{cotai_inf}.1)$} &
\multicolumn{4}{|c|}{$(\ref{cotai_inf}.2)$}   \\
$|V|$ & $p$ & $\alpha$ &
$t(\#)$ & $t^*/gap^*$ & $\#nodes$ & $gap_{LR}$ &
$t(\#)$ & $t^*/gap^*$ & $\#nodes$ & $gap_{LR}$ &
$t(\#)$ & $t^*/gap^*$ & $\#nodes$ & $gap_{LR}$\\
\hline
100	&	4	&	0.4	&	0.5 & 0.6 & 15 & 55.79	&	0.4 & 0.6 & \textbf{12} & 55.83	&	0.5 & 0.7 & 18 & 55.83	\\
100	&	4	&	0.6	&	0.5 & 0.6 & 41 & 40.17	&	\textbf{0.4} & \textbf{0.5} & 22 & 40.22	&	0.5 & 0.6 & 29 & 40.22	\\
100	&	4	&	0.8	&	0.4 & 0.5 & 61 & 24.26	&	0.4 & 0.6 & 130 & 24.33	&	0.4 & 0.5 & \textbf{54} & 24.33	\\
100	&	7	&	0.4	&	\textbf{0.6} & \textbf{0.7} & 200 & 52.77	&	\textbf{0.6} & \textbf{0.7} & \textbf{60} & 52.78	&	0.7 & 1 & 252 & 52.78	\\
100	&	7	&	0.6	&	0.7 & 0.8 & 360 & 38.19	&	\textbf{0.6} & \textbf{0.7} & \textbf{110} & 38.21	&	0.8 & 1.2 & 239 & 38.21	\\
100	&	7	&	0.8	&	0.9 & 1.6 & 839 & 23.76	&	\textbf{0.7} & \textbf{1.1} & \textbf{355} & 23.78	&	1.3 & 2.6 & 1557 & 23.78	\\
100	&	10	&	0.4	&	2.1 & 4.2 & 6658 & 51.98	&	\textbf{0.7} & \textbf{0.8} & \textbf{118} & 52.05	&	1 & 1.6 & 505 & 52.05	\\
100	&	10	&	0.6	&	4.1 & 13.2 & 16386 & 37.89	&	\textbf{0.9} & \textbf{1.1} & \textbf{252} & 37.98	&	1.8 & 3.4 & 1928 & 37.98	\\
100	&	10	&	0.8	&	5.5 & 27.9 & 23353 & 24.83	&	\textbf{2} & \textbf{9.5} & \textbf{4929} & 24.95	&	4.5 & 11 & 8541 & 24.95	\\
225	&	4	&	0.4	&	0.8 & 1 & 48 & 55.77	&	0.9 & 1.1 & \textbf{31} & 55.78	&	0.9 & 1.1 & 41 & 55.78	\\
225	&	4	&	0.6	&	0.8 & \textbf{1} & \textbf{44} & 39.42	&	0.9 & \textbf{1} & 55 & 39.43	&	0.9 & 1.2 & 78 & 39.43	\\
225	&	4	&	0.8	&	\textbf{0.8} & 1.1 & 95 & 22.13	&	\textbf{0.8} & 1.1 & 87 & 22.14	&	0.9 & 1.2 & 87 & 22.14	\\
225	&	7	&	0.4	&	\textbf{1.2} & \textbf{1.3} & 99 & 52.61	&	1.3 & 1.4 & \textbf{49} & 52.66	&	1.4 & 1.6 & 180 & 52.66	\\
225	&	7	&	0.6	&	3.3 & 8.8 & 1554 & 37.63	&	\textbf{1.5} & \textbf{1.7} & \textbf{151} & 37.69	&	4.1 & 9.8 & 1703 & 37.69	\\
225	&	7	&	0.8	&	4.6 & 22.1 & 3082 & 22.76	&	\textbf{2} & \textbf{3.9} & \textbf{728} & 22.83	&	5.1 & 22.6 & 3163 & 22.83	\\
225	&	10	&	0.4	&	9.1 & 62.7 & 6427 & 51.68	&	\textbf{2.2} & \textbf{2.9} & \textbf{168} & 51.75	&	4.5 & 11.1 & 2675 & 51.75	\\
225	&	10	&	0.6	&	15.2 & 56.6 & 10148 & 37.07	&	\textbf{2.9} & \textbf{6.1} & \textbf{743} & 37.16	&	15.7 & 74.5 & 10538 & 37.16	\\
225	&	10	&	0.8	&	38.1 & 147.8 & 41223 & 23.12	&	\textbf{13.3} & 94.5 & 14078 & 23.23	&	24.6 & 106.9 & 23284 & 23.23	\\
400	&	4	&	0.4	&	1.4 & 1.8 & 57 & 55.07	&	1.6 & 2.1 & \textbf{53} & 55.1	&	1.6 & 2.3 & 58 & 55.1	\\
400	&	4	&	0.6	&	1.6 & 2 & 95 & 38.71	&	1.9 & 2.5 & 107 & 38.75	&	1.7 & 2 & \textbf{87} & 38.75	\\
400	&	4	&	0.8	&	1.8 & 2.9 & 182 & 21.57	&	\textbf{1.7} & \textbf{2.2} & \textbf{146} & 21.62	&	2.1 & 3.7 & 252 & 21.62	\\
400	&	7	&	0.4	&	6.5 & 41.1 & 1102 & 52.72	&	\textbf{2.9} & \textbf{3.8} & \textbf{147} & 52.74	&	8.5 & 36.7 & 1761 & 52.74	\\
400	&	7	&	0.6	&	\textbf{9.4} & 62.6 & \textbf{2952} & 37.41	&	63 (9) & 3.96\% & 14032 & 37.44	&	34 & 275.8 & 16442 & 37.44	\\
400	&	7	&	0.8	&	8.1 & 30.2 & 1994 & 21.87	&	\textbf{4.7} & \textbf{13.6} & \textbf{1333} & 21.91	&	15.2 & 63.6 & 4701 & 21.91	\\
400	&	10	&	0.4	&	158.5 (9) & 1.09\% & 100979 & 51.8	&	\textbf{4.8} & \textbf{8.1} & \textbf{316} & 51.83	&	56.5 & 427.7 & 19284 & 51.83	\\
400	&	10	&	0.6	&	61.8 & 121.5 & 37448 & 36.48	&	66 (9) & 9.45\% & 13244 & 36.55	&	\textbf{16} & \textbf{74.8} & \textbf{2806} & 36.52	\\
400	&	10	&	0.8	&	229.9 (8) & 0.61\% & 143034 & 21.82	&	\textbf{7.3} & \textbf{10.1} & \textbf{1561} & 21.87	&	33.7 & 103 & 8891 & 21.87	\\

\hline
\end{tabular}
\end{scriptsize}
\end{center}
\caption{Results obtained for the OWAP formulations with the Shortest Path Problem and valid inequalities}
%\label{tab_exact_resolution}
\end{table}
%-----------------------------------------------------------------------------------------------------------------------

%-----------------------------------------------------------------------------------------------------------------------
\begin{table}[h!]
\renewcommand{\esp}{@{\hspace{0.15cm}}}
\begin{center}
\begin{scriptsize}
\begin{tabular}{|@{}\esp c \esp c \esp c\esp |c\esp c\esp c\esp c\esp |c\esp c\esp c\esp c\esp |c\esp c\esp c\esp c\esp  @{}|}
\hline
\multicolumn{3}{|c|}{Inst} &
\multicolumn{4}{|c|}{$(\ref{cotai_inf_ord}.1)$} &
\multicolumn{4}{|c|}{$(\ref{cotai_inf_ord}.2)$} &
\multicolumn{4}{|c|}{$(\ref{cotai_uij_max}.1)$}   \\
$|V|$ & $p$ & $\alpha$ &
$t(\#)$ & $t^*/gap^*$ & $\#nodes$ & $gap_{LR}$ &
$t(\#)$ & $t^*/gap^*$ & $\#nodes$ & $gap_{LR}$ &
$t(\#)$ & $t^*/gap^*$ & $\#nodes$ & $gap_{LR}$\\
\hline
100	&	4	&	0.4	&	0.7 & 3.2 & 755 & 15.59	&	0.6 & 0.8 & 114 & 55.83	&	0.5 & 0.9 & 17 & 55.83	\\
100	&	4	&	0.6	&	0.6 & 2.8 & 620 & 15.9	&	0.5 & 0.6 & 89 & 40.22	&	\textbf{0.4} & 0.6 & \textbf{16} & 40.22	\\
100	&	4	&	0.8	&	0.5 & 0.9 & 256 & 12.76	&	0.5 & 0.6 & 93 & 24.33	&	0.4 & 0.5 & 57 & 24.33	\\
100	&	7	&	0.4	&	149.3 (9) & 4.93\% & 348607 & 17.9	&	2.3 & 9.2 & 5053 & 52.78	&	\textbf{0.6} & 0.8 & 259 & 52.78	\\
100	&	7	&	0.6	&	13.9 & 87.5 & 32996 & 17.78	&	3 & 10.9 & 6989 & 38.21	&	0.8 & 2 & 716 & 38.21	\\
100	&	7	&	0.8	&	6.5 & 19.7 & 13584 & 14.33	&	2.3 & 4.8 & 4222 & 23.78	&	1 & 1.7 & 1251 & 23.78	\\
100	&	10	&	0.4	&	109.5 (9) & 6.11\% & 203304 & 16.28	&	114.8 (9) & 0.51\% & 360306 & 52.05	&	1.5 & 3.3 & 2706 & 52.05	\\
100	&	10	&	0.6	&	70.1 (9) & 7.7\% & 146569 & 17.4	&	86.9 & 265.4 & 297886 & 37.98	&	4 & 22.9 & 14402 & 37.98	\\
100	&	10	&	0.8	&	22.8 & 90.6 & 54826 & 15.59	&	300.9 (6) & 3.83\% & 906920 & 24.95	&	9.2 & 28.8 & 33580 & 24.95	\\
225	&	4	&	0.4	&	180.4 (7) & 9.39\% & 118948 & 17.47	&	1.1 & 1.4 & 158 & 55.78	&	0.8 & 1.2 & 56 & 55.78	\\
225	&	4	&	0.6	&	2.4 & 9.3 & 1046 & 16.08	&	1.2 & 1.9 & 162 & 39.43	&	0.9 & 1.1 & 50 & 39.43	\\
225	&	4	&	0.8	&	12.8 & 94.4 & 10292 & 10.89	&	1.3 & 2.1 & 322 & 22.14	&	\textbf{0.8} & 1.2 & \textbf{78} & 22.14	\\
225	&	7	&	0.4	&	219.8 (7) & 8.78\% & 129203 & 16.2	&	12 & 35 & 11614 & 52.66	&	1.3 & 1.5 & 305 & 52.66	\\
225	&	7	&	0.6	&	204.1 (7) & 9.46\% & 117810 & 16.35	&	24.9 & 125.4 & 25656 & 37.69	&	5.8 & 38.7 & 4909 & 37.69	\\
225	&	7	&	0.8	&	118 (9) & 2.4\% & 84612 & 12.92	&	21.9 & 86.3 & 20357 & 22.83	&	8 & 54.2 & 4975 & 22.83	\\
225	&	10	&	0.4	&	433.4 (3) & 14.03\% & 201310 & 15.38	&	431.2 (5) & 2.2\% & 397948 & 51.75	&	4.2 & 17.7 & 2961 & 51.75	\\
225	&	10	&	0.6	&	441.3 (3) & 11.29\% & 226890 & 16.16	&	349.9 (6) & 4.53\% & 286952 & 37.21	&	66.7 (9) & 24.23\% & 33790 & 37.16	\\
225	&	10	&	0.8	&	176.9 (8) & 5.05\% & 119062 & 13.5	&	465.5 (3) & 3.5\% & 338320 & 23.23	&	18.6 & 76 & 16194 & 23.23	\\
400	&	4	&	0.4	&	64.4 (9) & 11.11\% & 15532 & 16.39	&	3.2 & 4.9 & 377 & 55.1	&	1.4 & 1.8 & 83 & 55.1	\\
400	&	4	&	0.6	&	130.3 (8) & 10.07\% & 35381 & 15.26	&	2.8 & 7.6 & 504 & 38.75	&	1.6 & 2.1 & 118 & 38.75	\\
400	&	4	&	0.8	&	183.1 (7) & 5.01\% & 52444 & 10.33	&	2.9 & 5.5 & 420 & 21.62	&	1.8 & 3.9 & 267 & 21.62	\\
400	&	7	&	0.4	&	224.9 (7) & 12.7\% & 46806 & 14.9	&	51 & 127.8 & 20041 & 52.74	&	7.4 & 22.3 & 2715 & 52.74	\\
400	&	7	&	0.6	&	205.6 (8) & 9.81\% & 72507 & 15.04	&	50.7 & 344.7 & 19488 & 37.44	&	65.8 (9) & 31.74\% & 12801 & 37.44	\\
400	&	7	&	0.8	&	114.5 (9) & 5.76\% & 36406 & 11.47	&	146.9 & 533.1 & 54415 & 21.91	&	18.8 & 66.4 & 8856 & 21.91	\\
400	&	10	&	0.4	&	367.7 (6) & 15.03\% & 102646 & 15.45	&	599.8 (0) & 3.87\% & 140360 & 51.86	&	368 (6) & 21.32\% & 244981 & 51.87	\\
400	&	10	&	0.6	&	284.8 (8) & 1.36\% & 151282 & 15.02	&	418.1 (4) & 2.04\% & 86835 & 36.55	&	76.9 & 502.4 & 52942 & 36.52	\\
400	&	10	&	0.8	&	546 (3) & 5.96\% & 253528 & 11.94	&	455.1 (3) & 4.23\% & 73085 & 21.87	&	150.9 (9) & 1.67\% & 85416 & 21.87	\\

\hline
\end{tabular}
\end{scriptsize}
\end{center}
\caption{Results obtained for the OWAP formulations with the Shortest Path Problem and valid inequalities}
%\label{tab_exact_resolution}
\end{table}
%-----------------------------------------------------------------------------------------------------------------------
%=======================================================================================================================
\newpage
%=======================================================================================================================
%\subsubsection*{SP con DV:  $F_{R2}^{z}$ + (29.2), (30.1), (30.2), (31), (33.1), (33.2)}
%-----------------------------------------------------------------------------------------------------------------------
\begin{table}[h!]
\renewcommand{\esp}{@{\hspace{0.15cm}}}
\begin{center}
\begin{scriptsize}
\begin{tabular}{|@{}\esp c \esp c \esp c\esp |c\esp c\esp c\esp c\esp |c\esp c\esp c\esp c\esp |c\esp c\esp c\esp c\esp  @{}|}
\hline
\multicolumn{3}{|c|}{Inst} &
\multicolumn{4}{|c|}{$(\ref{cotai_uij_max}.2)$} &
\multicolumn{4}{|c|}{$(\ref{cotaj_uij_max}.1)$} &
\multicolumn{4}{|c|}{$(\ref{cotaj_uij_max}.2)$}   \\
$|V|$ & $p$ & $\alpha$ &
$t(\#)$ & $t^*/gap^*$ & $\#nodes$ & $gap_{LR}$ &
$t(\#)$ & $t^*/gap^*$ & $\#nodes$ & $gap_{LR}$ &
$t(\#)$ & $t^*/gap^*$ & $\#nodes$ & $gap_{LR}$\\
\hline
100	&	4	&	0.4	&	0.5 & 0.9 & 17 & 55.83	&	\textbf{0.3} & \textbf{0.4} & \textbf{12} & \textbf{14.96}	&	0.6 & 0.8 & 93 & 55.83	\\
100	&	4	&	0.6	&	\textbf{0.4} & 0.6 & \textbf{16} & 40.22	&	\textbf{0.4} & \textbf{0.5} & 24 & \textbf{15.52}	&	0.6 & 1 & 141 & 40.22	\\
100	&	4	&	0.8	&	0.4 & 0.5 & 57 & 24.33	&	\textbf{0.3} & \textbf{0.4} & 55 & 12.58	&	0.5 & 0.7 & 83 & 24.33	\\
100	&	7	&	0.4	&	\textbf{0.6} & \textbf{0.7} & 259 & 52.78	&	\textbf{0.6} & \textbf{0.7} & 148 & \textbf{16.81}	&	2.2 & 8.4 & 5035 & 52.78	\\
100	&	7	&	0.6	&	0.8 & 2.2 & 716 & 38.21	&	0.7 & 1.4 & 403 & \textbf{17.18}	&	2.8 & 7.4 & 5790 & 38.21	\\
100	&	7	&	0.8	&	1 & 1.7 & 1251 & 23.78	&	0.8 & 1.4 & 840 & 14.04	&	2.6 & 4.2 & 4556 & 23.78	\\
100	&	10	&	0.4	&	1.4 & 3.3 & 2706 & 52.05	&	0.8 & 1.2 & 468 & \textbf{15.76}	&	190.5 (9) & 0.61\% & 588985 & 52.05	\\
100	&	10	&	0.6	&	4.2 & 25.3 & 14402 & 37.98	&	1.4 & 4.5 & 2084 & \textbf{17.05}	&	221.7 (8) & 1\% & 692638 & 37.98	\\
100	&	10	&	0.8	&	9 & 27.7 & 33580 & 24.95	&	3.7 & 22.9 & 10167 & \textbf{15.44}	&	369.2 (4) & 3.12\% & 1059031 & 24.95	\\
225	&	4	&	0.4	&	0.8 & 1.1 & 56 & 55.78	&	\textbf{0.7} & \textbf{0.8} & 44 & \textbf{16.8}	&	1.3 & 1.8 & 248 & 55.78	\\
225	&	4	&	0.6	&	0.8 & 1.1 & 50 & 39.43	&	\textbf{0.7} & 1.1 & 134 & \textbf{15.66}	&	1.3 & 1.7 & 165 & 39.43	\\
225	&	4	&	0.8	&	\textbf{0.8} & \textbf{1} & \textbf{78} & 22.14	&	\textbf{0.8} & \textbf{1} & 167 & 10.68	&	1.5 & 2.5 & 434 & 22.14	\\
225	&	7	&	0.4	&	1.3 & 1.7 & 305 & 52.66	&	1.3 & 1.9 & 233 & \textbf{15.72}	&	10 & 32.3 & 8756 & 52.66	\\
225	&	7	&	0.6	&	5.8 & 37.8 & 4909 & 37.69	&	1.6 & 2.1 & 615 & \textbf{16.06}	&	20.7 & 99.8 & 20776 & 37.69	\\
225	&	7	&	0.8	&	8 & 54.1 & 4975 & 22.83	&	3.4 & 7.8 & 2644 & 12.78	&	29.2 & 98.8 & 29227 & 22.83	\\
225	&	10	&	0.4	&	4.1 & 17.1 & 2961 & 51.75	&	2.4 & 5 & 823 & \textbf{14.64}	&	551.5 (4) & 3.06\% & 540678 & 51.75	\\
225	&	10	&	0.6	&	66.8 (9) & 24.23\% & 33837 & 37.16	&	3 & 7.7 & 1231 & \textbf{15.67}	&	404.2 (5) & 3.9\% & 360472 & 37.16	\\
225	&	10	&	0.8	&	18.5 & 74.8 & 16194 & 23.23	&	14.3 & 63.4 & 13236 & \textbf{13.38}	&	388 (4) & 4.9\% & 271631 & 23.23	\\
400	&	4	&	0.4	&	1.4 & 1.9 & 83 & 55.1	&	\textbf{1.2} & \textbf{1.5} & 55 & \textbf{15.91}	&	3.1 & 4.9 & 335 & 55.1	\\
400	&	4	&	0.6	&	1.6 & 2.1 & 118 & 38.75	&	\textbf{1.4} & \textbf{1.7} & 99 & \textbf{14.95}	&	2.7 & 3.4 & 200 & 38.75	\\
400	&	4	&	0.8	&	1.8 & 3.8 & 267 & 21.62	&	\textbf{1.7} & 2.7 & 269 & 10.19	&	3.5 & 9.2 & 624 & 21.62	\\
400	&	7	&	0.4	&	7.4 & 22.7 & 2715 & 52.74	&	63.4 (9) & 3.1\% & 15386 & \textbf{14.27}	&	86.8 & 375 & 37215 & 52.74	\\
400	&	7	&	0.6	&	65.8 (9) & 31.74\% & 12805 & 37.44	&	64.8 (9) & 6.76\% & 16454 & \textbf{14.79}	&	86.4 (9) & 0.13\% & 34712 & 37.44	\\
400	&	7	&	0.8	&	18.8 & 66.2 & 8856 & 21.91	&	10.6 & 43.8 & 4798 & 11.3	&	116.5 (9) & 1.3\% & 38204 & 21.91	\\
400	&	10	&	0.4	&	367.9 (6) & 21.32\% & 244648 & 51.87	&	50.8 & 208.8 & 35166 & \textbf{14.75}	&	599.4 (0) & 4.92\% & 129856 & 51.83	\\
400	&	10	&	0.6	&	76.8 & 501.9 & 52942 & 36.52	&	116.6 (9) & 0.63\% & 75261 & \textbf{14.79}	&	471.6 (3) & 1.05\% & 111124 & 36.52	\\
400	&	10	&	0.8	&	150.9 (9) & 1.66\% & 85478 & 21.87	&	143.7 & 431 & 80175 & \textbf{11.83}	&	538.8 (2) & 2.29\% & 97542 & 21.96	\\

\hline
\end{tabular}
\end{scriptsize}
\end{center}
\caption{Results obtained for the OWAP formulations with the Shortest Path Problem and valid inequalities}
%\label{tab_exact_resolution}
\end{table}
%-----------------------------------------------------------------------------------------------------------------------

%-----------------------------------------------------------------------------------------------------------------------
\begin{table}[h!]
\renewcommand{\esp}{@{\hspace{0.15cm}}}
\begin{center}
\begin{scriptsize}
\begin{tabular}{|@{}\esp c \esp c \esp c\esp |c\esp c\esp c\esp c\esp |c\esp c\esp c\esp c\esp |c\esp c\esp c\esp c\esp  @{}|}
\hline
\multicolumn{3}{|c|}{Inst} &
\multicolumn{4}{|c|}{$(\ref{cotazy})$} &
\multicolumn{4}{|c|}{$(\ref{validsubsets}.1)$} &
\multicolumn{4}{|c|}{$(\ref{validsubsets}.2)$}   \\
$|V|$ & $p$ & $\alpha$ &
$t(\#)$ & $t^*/gap^*$ & $\#nodes$ & $gap_{LR}$ &
$t(\#)$ & $t^*/gap^*$ & $\#nodes$ & $gap_{LR}$ &
$t(\#)$ & $t^*/gap^*$ & $\#nodes$ & $gap_{LR}$\\
\hline
100	&	4	&	0.4	&	0.7 & 3.9 & 590 & 55.83	&	2.9 & 8.2 & 4269 & 53.72	&	0.5 & 0.9 & 17 & 55.83	\\
100	&	4	&	0.6	&	0.5 & 0.6 & 91 & 40.22	&	3.3 & 10.8 & 4845 & 37.39	&	\textbf{0.4} & 0.6 & \textbf{16} & 40.22	\\
100	&	4	&	0.8	&	0.5 & 0.8 & 345 & 24.33	&	6.7 & 18.5 & 11571 & 20.79	&	0.4 & 0.5 & 57 & 24.33	\\
100	&	7	&	0.4	&	4.2 & 5.4 & 11077 & 52.78	&	122.3 (8) & 33.45\% & 386092 & 51.11	&	\textbf{0.6} & \textbf{0.7} & 269 & 52.78	\\
100	&	7	&	0.6	&	8.1 & 17.5 & 23353 & 38.21	&	4.1 & 10.6 & 6101 & 36.01	&	0.9 & 2.3 & 692 & 38.21	\\
100	&	7	&	0.8	&	19.2 & 29.2 & 57034 & 23.78	&	15.9 & 95.2 & 31791 & 21.08	&	1 & 1.7 & 1250 & 23.78	\\
100	&	10	&	0.4	&	307.1 (9) & 4.38\% & 576477 & 52.15	&	69.7 (9) & 27.74\% & 183041 & 49.29	&	2 & 5.1 & 6573 & 52.05	\\
100	&	10	&	0.6	&	557.5 (3) & 10.25\% & 717764 & 39.6	&	131.3 (8) & 23.08\% & 200947 & 34.41	&	5.4 & 26.1 & 21118 & 37.98	\\
100	&	10	&	0.8	&	599.3 (0) & 15.85\% & 592851 & 29.22	&	40.2 & 287.4 & 86591 & 20.61	&	13.4 & 35.8 & 53741 & 24.95	\\
225	&	4	&	0.4	&	1.3 & 1.8 & 129 & 55.78	&	9 & 15.2 & 4120 & 54.96	&	0.8 & 1.2 & 56 & 55.78	\\
225	&	4	&	0.6	&	1.3 & 2 & 240 & 39.43	&	96.5 (9) & 31.94\% & 49866 & 38.31	&	0.9 & 1.2 & 50 & 39.43	\\
225	&	4	&	0.8	&	1.5 & 2.3 & 473 & 22.14	&	168.2 (8) & 14.5\% & 87717 & 20.7	&	\textbf{0.8} & 1.1 & \textbf{78} & 22.14	\\
225	&	7	&	0.4	&	32.7 & 173.1 & 29035 & 52.66	&	140.1 (8) & 47.83\% & 59842 & 51.23	&	1.4 & 2.1 & 433 & 52.66	\\
225	&	7	&	0.6	&	40.9 & 142.1 & 48835 & 37.69	&	329.2 (5) & 30.8\% & 154197 & 35.76	&	3.8 & 17.3 & 2557 & 37.69	\\
225	&	7	&	0.8	&	96.6 & 191.5 & 111329 & 22.83	&	558.4 (1) & 14.28\% & 295442 & 20.44	&	3.7 & 11.7 & 2477 & 22.83	\\
225	&	10	&	0.4	&	597 (2) & 8.76\% & 252282 & 52.84	&	433.3 (3) & 47.63\% & 181402 & 50.44	&	4.3 & 18 & 3034 & 51.75	\\
225	&	10	&	0.6	&	599.7 (0) & 14.19\% & 197775 & 40.63	&	563.7 (1) & 31.68\% & 245679 & 35.75	&	10.5 & 40.5 & 5930 & 37.16	\\
225	&	10	&	0.8	&	599.7 (0) & 16.4\% & 228494 & 28.9	&	599.8 (0) & 14.98\% & 275289 & 20.85	&	18.6 & 70.1 & 16060 & 23.23	\\
400	&	4	&	0.4	&	3.2 & 4.4 & 221 & 55.1	&	193.1 (7) & 53.05\% & 40227 & 54.59	&	1.4 & 1.8 & 83 & 55.1	\\
400	&	4	&	0.6	&	3.6 & 5.7 & 347 & 38.75	&	135.5 (8) & 35.68\% & 27211 & 37.89	&	1.5 & 2 & 118 & 38.75	\\
400	&	4	&	0.8	&	5.8 & 13.7 & 1656 & 21.62	&	324.8 (5) & 17.63\% & 70658 & 20.47	&	1.8 & 3.6 & 274 & 21.62	\\
400	&	7	&	0.4	&	128.5 (9) & 2.03\% & 46011 & 52.77	&	494.4 (2) & 52.93\% & 79670 & 52.78	&	6.6 & 31.5 & 1675 & 52.74	\\
400	&	7	&	0.6	&	156 & 373.3 & 70111 & 37.44	&	444.5 (3) & 35.54\% & 82438 & 36.56	&	64.3 (9) & 31.7\% & 11699 & 37.44	\\
400	&	7	&	0.8	&	317.9 (8) & 0.73\% & 164129 & 21.91	&	551.3 (1) & 17.6\% & 96906 & 20.33	&	13.1 & 49.8 & 5172 & 21.91	\\
400	&	10	&	0.4	&	599.5 (0) & 8.61\% & 108900 & 53.5	&	480.2 (2) & 51.52\% & 72765 & 51.52	&	370 (6) & 21.32\% & 240940 & 51.87	\\
400	&	10	&	0.6	&	599.7 (0) & 13.96\% & 133291 & 40.26	&	427.5 (3) & 35.57\% & 70047 & 35.78	&	86.6 & 596.7 & 62085 & 36.52	\\
400	&	10	&	0.8	&	599.7 (0) & 28.77\% & 115387 & 27.86	&	426.2 (3) & 18.35\% & 65582 & 20.18	&	151 (9) & 1.19\% & 88335 & 21.87	\\

\hline
\end{tabular}
\end{scriptsize}
\end{center}
\caption{Results obtained for the OWAP formulations with the Shortest Path Problem and valid inequalities}
%\label{tab_exact_resolution}
\end{table}
%-----------------------------------------------------------------------------------------------------------------------
%=======================================================================================================================
\newpage
%=======================================================================================================================
%\subsubsection*{SP con DV:  $F_{R2}^{z}$ + (33.3), (33.4), (34), (35), (36), (37)}
%-----------------------------------------------------------------------------------------------------------------------
\begin{table}[h!]
\renewcommand{\esp}{@{\hspace{0.15cm}}}
\begin{center}
\begin{scriptsize}
\begin{tabular}{|@{}\esp c \esp c \esp c\esp |c\esp c\esp c\esp c\esp |c\esp c\esp c\esp c\esp |c\esp c\esp c\esp c\esp  @{}|}
\hline
\multicolumn{3}{|c|}{Inst} &
\multicolumn{4}{|c|}{$(\ref{validsubsets}.3)$} &
\multicolumn{4}{|c|}{$(\ref{validsubsets}.4)$} &
\multicolumn{4}{|c|}{$(\ref{vi:owa2eq})$}   \\
$|V|$ & $p$ & $\alpha$ &
$t(\#)$ & $t^*/gap^*$ & $\#nodes$ & $gap_{LR}$ &
$t(\#)$ & $t^*/gap^*$ & $\#nodes$ & $gap_{LR}$ &
$t(\#)$ & $t^*/gap^*$ & $\#nodes$ & $gap_{LR}$\\
\hline
100	&	4	&	0.4	&	0.5 & 0.9 & 17 & 55.83	&	0.5 & 0.9 & 17 & 55.83	&	0.4 & 0.6 & 212 & 41.06	\\
100	&	4	&	0.6	&	\textbf{0.4} & 0.6 & \textbf{16} & 40.22	&	\textbf{0.4} & 0.6 & \textbf{16} & 40.22	&	0.5 & 0.7 & 332 & 20.22	\\
100	&	4	&	0.8	&	0.4 & \textbf{0.4} & 57 & 24.33	&	0.4 & 0.6 & 57 & 24.33	&	0.6 & 1 & 793 & \textbf{5.31}	\\
100	&	7	&	0.4	&	\textbf{0.6} & \textbf{0.7} & 259 & 52.78	&	\textbf{0.6} & \textbf{0.7} & 259 & 52.78	&	4.4 & 6.2 & 14452 & 44.9	\\
100	&	7	&	0.6	&	0.9 & 2.2 & 760 & 38.21	&	0.9 & 2.1 & 735 & 38.21	&	11.9 & 19.1 & 40976 & 27.89	\\
100	&	7	&	0.8	&	1 & 1.6 & 1268 & 23.78	&	1 & 1.6 & 1197 & 23.78	&	28.1 & 45 & 93501 & \textbf{11.05}	\\
100	&	10	&	0.4	&	1.4 & 3.2 & 2675 & 52.05	&	1.5 & 3.2 & 3082 & 52.05	&	486 (6) & 5.58\% & 1335239 & 46.98	\\
100	&	10	&	0.6	&	4.7 & 20.4 & 14983 & 37.98	&	4.1 & 22.9 & 14784 & 37.98	&	578.5 (1) & 13.53\% & 1402582 & 33.7	\\
100	&	10	&	0.8	&	9.9 & 35.3 & 39048 & 24.95	&	10.5 & 33.3 & 41792 & 24.95	&	599.2 (0) & 19.81\% & 1621965 & 22.08	\\
225	&	4	&	0.4	&	0.8 & 1.1 & 56 & 55.78	&	0.8 & 1.1 & 56 & 55.78	&	1.2 & 2 & 583 & 41.02	\\
225	&	4	&	0.6	&	0.9 & 1.2 & 50 & 39.43	&	0.9 & 1.2 & 50 & 39.43	&	1.3 & 2.1 & 888 & 19.23	\\
225	&	4	&	0.8	&	\textbf{0.8} & 1.1 & \textbf{78} & 22.14	&	\textbf{0.8} & \textbf{1} & \textbf{78} & 22.14	&	1.8 & 4 & 1621 & \textbf{2.65}	\\
225	&	7	&	0.4	&	1.3 & 1.6 & 321 & 52.66	&	1.4 & 2.3 & 462 & 52.66	&	19.3 & 31.1 & 30039 & 44.72	\\
225	&	7	&	0.6	&	3.2 & 12.8 & 1753 & 37.69	&	19.9 & 178.6 & 18733 & 37.69	&	61.9 & 165.6 & 86213 & 27.24	\\
225	&	7	&	0.8	&	5.7 & 31.2 & 3487 & 22.83	&	43.2 & 407.7 & 43205 & 22.83	&	185.3 & 576.3 & 234029 & \textbf{9.89}	\\
225	&	10	&	0.4	&	4.2 & 18.3 & 3135 & 51.75	&	4.5 & 13.5 & 3276 & 51.75	&	599.8 (0) & 11.07\% & 669308 & 47.51	\\
225	&	10	&	0.6	&	66.1 (9) & 20.67\% & 36169 & 37.16	&	69.5 (9) & 23.5\% & 38304 & 37.16	&	599.8 (0) & 14.89\% & 695871 & 33.05	\\
225	&	10	&	0.8	&	77 (9) & 9.92\% & 62188 & 23.31	&	18.3 & 68.9 & 15974 & 23.23	&	599.8 (0) & 19.27\% & 726358 & 20.83	\\
400	&	4	&	0.4	&	1.4 & 1.8 & 83 & 55.1	&	1.4 & 2 & 67 & 55.1	&	2.1 & 3.9 & 522 & 40.1	\\
400	&	4	&	0.6	&	1.5 & 2 & 118 & 38.75	&	1.5 & 1.9 & 111 & 38.75	&	4.1 & 9.3 & 1439 & 18.28	\\
400	&	4	&	0.8	&	1.8 & 3.9 & 267 & 21.62	&	\textbf{1.7} & 3.5 & 246 & 21.62	&	7.8 & 15.4 & 5076 & \textbf{1.96}	\\
400	&	7	&	0.4	&	8.1 & 24.3 & 3123 & 52.74	&	6.2 & 28 & 1589 & 52.74	&	54.9 & 86.1 & 41042 & 44.84	\\
400	&	7	&	0.6	&	65.8 (9) & 31.92\% & 12983 & 37.48	&	64.4 (9) & 32.2\% & 11723 & 37.49	&	153.2 & 302.4 & 109497 & 26.98	\\
400	&	7	&	0.8	&	20.1 & 64.4 & 9184 & 21.91	&	13.3 & 60.2 & 4992 & 21.91	&	343.2 (8) & 1.15\% & 229819 & \textbf{8.85}	\\
400	&	10	&	0.4	&	373.7 (6) & 21.29\% & 247604 & 51.87	&	216.4 (8) & 51.31\% & 132409 & 51.87	&	599.8 (0) & 8.13\% & 342591 & 47.35	\\
400	&	10	&	0.6	&	84.1 & 567.4 & 55975 & 36.52	&	47.7 & 310.6 & 29281 & 36.52	&	599.9 (0) & 30.46\% & 303178 & 32.37	\\
400	&	10	&	0.8	&	167.2 (9) & 1.57\% & 93684 & 21.87	&	114.5 (9) & 1.06\% & 66020 & 21.87	&	599.8 (0) & 11.5\% & 312739 & 17.19	\\

\hline
\end{tabular}
\end{scriptsize}
\end{center}
\caption{Results obtained for the OWAP formulations with the Shortest Path Problem and valid inequalities}
%\label{tab_exact_resolution}
\end{table}
%-----------------------------------------------------------------------------------------------------------------------

%-----------------------------------------------------------------------------------------------------------------------
\begin{table}[h!]
\renewcommand{\esp}{@{\hspace{0.15cm}}}
\begin{center}
\begin{scriptsize}
\begin{tabular}{|@{}\esp c \esp c \esp c\esp |c\esp c\esp c\esp c\esp |c\esp c\esp c\esp c\esp |c\esp c\esp c\esp c\esp  @{}|}
\hline
\multicolumn{3}{|c|}{Inst} &
\multicolumn{4}{|c|}{$(\ref{vi:cotayydis1})$} &
\multicolumn{4}{|c|}{$(\ref{vi:cotayydis2})$} &
\multicolumn{4}{|c|}{$(\ref{vi:yyrel1})$}   \\
$|V|$ & $p$ & $\alpha$ &
$t(\#)$ & $t^*/gap^*$ & $\#nodes$ & $gap_{LR}$ &
$t(\#)$ & $t^*/gap^*$ & $\#nodes$ & $gap_{LR}$ &
$t(\#)$ & $t^*/gap^*$ & $\#nodes$ & $gap_{LR}$\\
\hline
100	&	4	&	0.4	&	0.5 & 0.7 & 18 & 55.83	&	0.5 & 0.9 & 164 & 55.83	&	0.5 & 0.7 & 147 & 55.83	\\
100	&	4	&	0.6	&	\textbf{0.4} & 0.6 & 31 & 40.22	&	0.5 & 1 & 392 & 40.22	&	0.7 & 0.9 & 475 & 40.22	\\
100	&	4	&	0.8	&	\textbf{0.3} & \textbf{0.4} & 62 & 24.33	&	0.5 & 1 & 499 & 24.33	&	0.8 & 1.8 & 1018 & 24.33	\\
100	&	7	&	0.4	&	\textbf{0.6} & 0.9 & 252 & 52.78	&	\textbf{0.6} & 0.8 & 523 & 52.78	&	10 & 30.1 & 15811 & 52.78	\\
100	&	7	&	0.6	&	0.8 & 1.9 & 616 & 38.21	&	1.2 & 2.1 & 2583 & 38.21	&	31.3 & 84.6 & 53345 & 38.21	\\
100	&	7	&	0.8	&	1.3 & 3 & 1599 & 23.78	&	3.2 & 7.4 & 9255 & 23.78	&	139.6 & 255.7 & 224889 & 23.78	\\
100	&	10	&	0.4	&	2.6 & 11.9 & 8116 & 52.05	&	2.7 & 12.8 & 9498 & 52.05	&	576.9 (1) & 8.58\% & 423512 & 52.11	\\
100	&	10	&	0.6	&	4.3 & 20 & 17641 & 37.98	&	23.2 & 76.7 & 87121 & 37.98	&	599.2 (0) & 8.39\% & 364546 & 37.98	\\
100	&	10	&	0.8	&	12 & 38.1 & 47989 & 24.95	&	95.6 & 287.8 & 332536 & 24.95	&	599.2 (0) & 10.03\% & 335942 & 25	\\
225	&	4	&	0.4	&	0.8 & 1 & 68 & 55.78	&	1.3 & 4.7 & 914 & 55.78	&	1.1 & 1.7 & 436 & 55.78	\\
225	&	4	&	0.6	&	0.9 & 1.1 & 61 & 39.43	&	1.2 & 3 & 952 & 39.43	&	3.5 & 15.2 & 3308 & 39.43	\\
225	&	4	&	0.8	&	\textbf{0.8} & 1.3 & 122 & 22.14	&	1.4 & 3.4 & 1098 & 22.14	&	2.8 & 5.9 & 2086 & 22.14	\\
225	&	7	&	0.4	&	1.3 & 2.1 & 282 & 52.66	&	5 & 23.2 & 6785 & 52.66	&	60.4 & 247.1 & 41320 & 52.66	\\
225	&	7	&	0.6	&	2.5 & 7.1 & 945 & 37.69	&	26.5 & 98 & 36980 & 37.69	&	239.6 & 495 & 152146 & 37.69	\\
225	&	7	&	0.8	&	6.2 & 27.3 & 4202 & 22.83	&	48.1 & 96.2 & 62740 & 22.83	&	520.5 (4) & 2.27\% & 306598 & 22.83	\\
225	&	10	&	0.4	&	4.4 & 19.6 & 3389 & 51.75	&	7.8 & 17 & 9596 & 51.75	&	599.9 (0) & 7.45\% & 124631 & 51.75	\\
225	&	10	&	0.6	&	8.2 & 33.8 & 4540 & 37.16	&	72.5 & 164.9 & 97785 & 37.16	&	599.9 (0) & 10.16\% & 107265 & 37.29	\\
225	&	10	&	0.8	&	14.2 & \textbf{58.2} & \textbf{12717} & 23.23	&	461.6 (6) & 3.04\% & 523762 & 23.23	&	599.9 (0) & 11.14\% & 101873 & 23.78	\\
400	&	4	&	0.4	&	1.5 & 1.8 & 74 & 55.1	&	2.1 & 4.7 & 547 & 55.1	&	5.2 & 23.5 & 2214 & 55.1	\\
400	&	4	&	0.6	&	1.6 & 2.1 & 99 & 38.75	&	8.1 & 41.1 & 5644 & 38.75	&	7.4 & 14.3 & 2689 & 38.75	\\
400	&	4	&	0.8	&	1.8 & 3.1 & 231 & 21.62	&	25 & 201.5 & 12344 & 21.62	&	9.6 & 20.5 & 3762 & 21.62	\\
400	&	7	&	0.4	&	7.4 & 35.8 & 2015 & 52.74	&	74.9 (9) & 0.77\% & 48469 & 52.74	&	293.1 (9) & 1.83\% & 74737 & 52.74	\\
400	&	7	&	0.6	&	10.8 & \textbf{39.4} & 4737 & 37.44	&	195.7 (9) & 0.41\% & 107708 & 37.44	&	554.2 (3) & 5.74\% & 107728 & 37.44	\\
400	&	7	&	0.8	&	69.7 (9) & 14.84\% & 18239 & 21.92	&	344.1 (6) & 2.37\% & 157933 & 21.91	&	599.4 (0) & 6.75\% & 101173 & 21.92	\\
400	&	10	&	0.4	&	308.8 (7) & 30.8\% & 205157 & 51.83	&	65.9 & 310.5 & 35458 & 51.83	&	600 (0) & 9.25\% & 51387 & 51.83	\\
400	&	10	&	0.6	&	143.7 (8) & 33.04\% & 73262 & 36.52	&	357.7 (6) & 2.89\% & 133142 & 36.53	&	599.6 (0) & 8.63\% & 47052 & 36.65	\\
400	&	10	&	0.8	&	190.8 (8) & 3.62\% & 90019 & 21.91	&	558 (1) & 4.04\% & 164201 & 21.87	&	599.7 (0) & 9.28\% & 45950 & 22.21	\\

\hline
\end{tabular}
\end{scriptsize}
\end{center}
\caption{Results obtained for the OWAP formulations with the Shortest Path Problem and valid inequalities}
%\label{tab_exact_resolution}
\end{table}
%-----------------------------------------------------------------------------------------------------------------------
%=======================================================================================================================
\newpage
%=======================================================================================================================
%\subsubsection*{SP con DV:  $F_{R2}^{z}$ + (38), (39)}
%-----------------------------------------------------------------------------------------------------------------------
\begin{table}[h!]
\renewcommand{\esp}{@{\hspace{0.15cm}}}
\begin{center}
\begin{scriptsize}
\begin{tabular}{|@{}\esp c \esp c \esp c\esp |c\esp c\esp c\esp c\esp |c\esp c\esp c\esp c\esp  @{}|}
\hline
\multicolumn{3}{|c|}{Inst} &
\multicolumn{4}{|c|}{$(\ref{vi:yyrel2})$} &
\multicolumn{4}{|c|}{$(\ref{vi:yyrel3})$}  \\
$|V|$ & $p$ & $\alpha$ &
$t(\#)$ & $t^*/gap^*$ & $\#nodes$ & $gap_{LR}$ &
$t(\#)$ & $t^*/gap^*$ & $\#nodes$ & $gap_{LR}$\\
\hline
100	&	4	&	0.4	&	0.4 & 0.6 & 84 & 55.83	&	0.4 & 0.5 & 97 & 55.83\\
100	&	4	&	0.6	&	0.6 & 1 & 425 & 40.22	&	0.5 & 0.7 & 195 & 40.22\\
100	&	4	&	0.8	&	0.6 & 0.9 & 515 & 24.33	&	0.6 & 1.3 & 675 & 24.33\\
100	&	7	&	0.4	&	4.4 & 11.2 & 10737 & 52.78	&	3.9 & 5.6 & 9681 & 52.78\\
100	&	7	&	0.6	&	20.6 & 51.1 & 54008 & 38.21	&	17.5 & 37 & 48030 & 38.21\\
100	&	7	&	0.8	&	58.9 & 100.2 & 158817 & 23.78	&	56.3 & 114.2 & 141438 & 23.78\\
100	&	10	&	0.4	&	404.9 (4) & 5.23\% & 780957 & 52.05	&	413.3 (5) & 2.24\% & 816185 & 52.05\\
100	&	10	&	0.6	&	543.8 (3) & 6.64\% & 959614 & 37.98	&	572.5 (1) & 6.16\% & 982457 & 37.98\\
100	&	10	&	0.8	&	599.2 (0) & 8.34\% & 971640 & 24.95	&	599.1 (0) & 8.47\% & 955372 & 24.95\\
225	&	4	&	0.4	&	0.8 & 1.1 & 281 & 55.78	&	1 & 1.4 & 398 & 55.78\\
225	&	4	&	0.6	&	2 & 7.3 & 1776 & 39.43	&	2 & 7.4 & 1515 & 39.43\\
225	&	4	&	0.8	&	1.5 & 2.5 & 997 & 22.14	&	2.5 & 8.5 & 1928 & 22.14\\
225	&	7	&	0.4	&	46.5 & 139.9 & 50060 & 52.66	&	43.8 & 143.6 & 46575 & 52.66\\
225	&	7	&	0.6	&	164.9 & 462.8 & 154252 & 37.69	&	140.4 & 391.4 & 122888 & 37.69\\
225	&	7	&	0.8	&	390.3 (9) & 0.75\% & 336722 & 22.83	&	323.1 & 592.7 & 272833 & 22.83\\
225	&	10	&	0.4	&	599.8 (0) & 5.65\% & 365109 & 51.75	&	599.8 (0) & 4.9\% & 344234 & 51.75\\
225	&	10	&	0.6	&	599.9 (0) & 9.46\% & 272912 & 37.29	&	599.9 (0) & 8.45\% & 269151 & 37.23\\
225	&	10	&	0.8	&	599.9 (0) & 10.47\% & 225292 & 23.67	&	599.9 (0) & 8.38\% & 219993 & 23.44\\
400	&	4	&	0.4	&	2.8 & 5.4 & 873 & 55.1	&	2.7 & 5.6 & 849 & 55.1\\
400	&	4	&	0.6	&	4.2 & 15.7 & 1746 & 38.75	&	6.9 & 20.1 & 3066 & 38.75\\
400	&	4	&	0.8	&	28.6 & 229 & 15015 & 21.62	&	6.4 & 11.4 & 2615 & 21.62\\
400	&	7	&	0.4	&	197.7 & 463.3 & 76117 & 52.74	&	138 & 275.6 & 45881 & 52.74\\
400	&	7	&	0.6	&	489.1 (5) & 4.08\% & 136882 & 37.46	&	442.8 (5) & 3.08\% & 108475 & 37.44\\
400	&	7	&	0.8	&	596.5 (1) & 5.97\% & 140394 & 21.95	&	590.4 (1) & 4.87\% & 124476 & 21.91\\
400	&	10	&	0.4	&	599.7 (0) & 7.47\% & 107208 & 51.86	&	599.9 (0) & 7.43\% & 98167 & 51.86\\
400	&	10	&	0.6	&	599.7 (0) & 7.64\% & 86331 & 36.53	&	599.7 (0) & 8\% & 80915 & 36.76\\
400	&	10	&	0.8	&	599.8 (0) & 8.46\% & 72226 & 22.27	&	599.4 (0) & 6.85\% & 71397 & 22.15\\
\hline
\end{tabular}
\end{scriptsize}
\end{center}
\caption{Results obtained for the OWAP formulations with the Shortest Path Problem and valid inequalities}
%\label{tab_exact_resolution}
\end{table}
%-----------------------------------------------------------------------------------------------------------------------
%=======================================================================================================================
\newpage
%=======================================================================================================================
%\subsubsection*{PMP con DV:  $F_{R1}^{z}$ + (25.1), (25.2), (26.1), (26.2), (29.1)}
%-----------------------------------------------------------------------------------------------------------------------
\begin{table}[h!]
\renewcommand{\esp}{@{\hspace{0.15cm}}}
\begin{center}
\begin{scriptsize}
\begin{tabular}{|@{}\esp c \esp c \esp c\esp |c\esp c\esp c\esp c\esp |c\esp c\esp c\esp c\esp |c\esp c\esp c\esp c\esp  @{}|}
\hline
\multicolumn{3}{|c|}{Inst} &
\multicolumn{4}{|c|}{$F_{R1}^{z}$} &
\multicolumn{4}{|c|}{$(\ref{cotai_inf}.1)$} &
\multicolumn{4}{|c|}{$(\ref{cotai_inf}.2)$}   \\
$|V|$ & $p$ & $\alpha$ &
$t(\#)$ & $t^*/gap^*$ & $\#nodes$ & $gap_{LR}$ &
$t(\#)$ & $t^*/gap^*$ & $\#nodes$ & $gap_{LR}$ &
$t(\#)$ & $t^*/gap^*$ & $\#nodes$ & $gap_{LR}$\\
\hline
100	&	4	&	0.4	&	\textbf{0.6} & \textbf{0.7} & 186 & 55.44	&	\textbf{0.6} & 0.8 & 190 & 55.44	&	\textbf{0.6} & 0.8 & 113 & 55.44	\\
100	&	4	&	0.6	&	\textbf{0.6} & 0.7 & 152 & 38.98	&	0.7 & 0.8 & 315 & 38.98	&	0.7 & 0.8 & 156 & 38.98	\\
100	&	4	&	0.8	&	\textbf{0.6} & \textbf{0.7} & 302 & 21.53	&	0.7 & 0.9 & 438 & 21.53	&	0.7 & 0.8 & 281 & 21.53	\\
100	&	7	&	0.4	&	\textbf{1} & 1.3 & 236 & 52.18	&	\textbf{1} & 1.5 & 367 & 52.18	&	1.2 & 1.5 & 220 & 52.18	\\
100	&	7	&	0.6	&	\textbf{1.1} & 1.4 & 480 & 35.97	&	\textbf{1.1} & 1.5 & 728 & 35.97	&	1.3 & 1.8 & 488 & 35.97	\\
100	&	7	&	0.8	&	1.4 & 2 & 965 & 20.27	&	1.4 & \textbf{1.9} & 1189 & 20.27	&	1.7 & 2.5 & 1118 & 20.27	\\
100	&	10	&	0.4	&	1.5 & \textbf{1.9} & 299 & 50.66	&	\textbf{1.4} & 2.1 & 417 & 50.66	&	1.8 & 2.3 & 290 & 50.66	\\
100	&	10	&	0.6	&	1.9 & 2.6 & 963 & 34.85	&	1.9 & \textbf{2.4} & 996 & 34.85	&	2.3 & 2.9 & 883 & 34.85	\\
100	&	10	&	0.8	&	6 & 19.4 & 6329 & 20.2	&	5.6 & \textbf{16.2} & 6617 & 20.2	&	7.8 & 28.4 & 5982 & 20.2	\\
225	&	4	&	0.4	&	2.1 & 4.4 & 1188 & 55.09	&	2.3 & 3.6 & 1343 & 55.09	&	2.3 & 3.9 & 1183 & 55.09	\\
225	&	4	&	0.6	&	1.7 & 2.9 & 1236 & 38.57	&	2.3 & 3.5 & 1664 & 38.57	&	1.9 & 3.7 & 1162 & 38.57	\\
225	&	4	&	0.8	&	\textbf{1.9} & 3.2 & 1101 & 21.09	&	2.4 & 5.5 & 1639 & 21.09	&	2.1 & 3.9 & 1032 & 21.09	\\
225	&	7	&	0.4	&	7.1 & 22.8 & 9208 & 52.34	&	8.7 & 24.8 & 6369 & 52.34	&	10.4 & 39.6 & 7204 & 52.34	\\
225	&	7	&	0.6	&	10 & 16 & \textbf{6038} & 36.27	&	10.7 & 16.4 & 7790 & 36.27	&	14.9 & 30.6 & 7280 & 36.27	\\
225	&	7	&	0.8	&	17.2 & 62.5 & 10491 & 20.32	&	16.8 & \textbf{46.1} & 11458 & 20.32	&	21.3 & 78.5 & 9834 & 20.32	\\
225	&	10	&	0.4	&	7.5 & 13.2 & 2136 & 50.25	&	\textbf{6.9} & 13.1 & 2592 & 50.25	&	9.8 & 28.8 & 1992 & 50.25	\\
225	&	10	&	0.6	&	32.4 & 123.2 & 15537 & 34.56	&	\textbf{29.5} & 104.5 & 15139 & 34.56	&	36.3 & \textbf{72} & 12789 & 34.56	\\
225	&	10	&	0.8	&	295 (8) & 0.32\% & 114029 & 19.62	&	292.3 (\textbf{9}) & 0.26\% & 143909 & 19.62	&	380.8 (\textbf{9}) & 0.39\% & 149239 & 19.62	\\
400	&	4	&	0.4	&	7.3 & 22.3 & 3345 & 55.37	&	6.7 & 13.5 & 3972 & 55.37	&	8.7 & 15.8 & 3276 & 55.37	\\
400	&	4	&	0.6	&	6.7 & 11.9 & 4103 & 39.04	&	7.4 & 12.2 & 4921 & 39.04	&	9.6 & 26.1 & 6014 & 39.04	\\
400	&	4	&	0.8	&	9 & 22.1 & 5397 & 21.03	&	10.9 & 21.7 & 6517 & 21.03	&	12.1 & 28.1 & 5661 & 21.03	\\
400	&	7	&	0.4	&	34.4 & 144.4 & 10464 & 52.05	&	31.7 & 100.7 & 10750 & 52.05	&	61.7 & 258 & 14738 & 52.05	\\
400	&	7	&	0.6	&	83.4 & 250.9 & 27604 & 36.12	&	98.8 & 207.7 & 41491 & 36.12	&	100.9 & 301.5 & 28857 & 36.12	\\
400	&	7	&	0.8	&	84.4 & 187.6 & 28762 & 20.19	&	105.2 & 236.7 & 40724 & 20.19	&	118.6 & 245.3 & 31177 & 20.19	\\
400	&	10	&	0.4	&	68.4 & 197.4 & 13777 & 50.58	&	\textbf{54.1} & \textbf{149.6} & 13423 & 50.58	&	92.5 & 215.7 & 13456 & 50.58	\\
400	&	10	&	0.6	&	289.4 (9) & 0.11\% & 61886 & 34.54	&	273.6 (9) & 0.14\% & 70779 & 34.54	&	395.2 (7) & 0.19\% & 57614 & 34.54	\\
400	&	10	&	0.8	&	583.5 (\textbf{1}) & 0.42\% & 97022 & 19.5	&	596.8 (\textbf{1}) & 0.32\% & 121574 & 19.49	&	599.8 (0) & 0.48\% & 90232 & 19.51	\\
\hline
\end{tabular}
\end{scriptsize}
\end{center}
\caption{Results obtained for the OWAP formulations with the Perfect Matching Problem and valid inequalities}
%\label{tab_exact_resolution}
\end{table}
%-----------------------------------------------------------------------------------------------------------------------

%-----------------------------------------------------------------------------------------------------------------------
\begin{table}[h!]
\renewcommand{\esp}{@{\hspace{0.15cm}}}
\begin{center}
\begin{scriptsize}
\begin{tabular}{|@{}\esp c \esp c \esp c\esp |c\esp c\esp c\esp c\esp |c\esp c\esp c\esp c\esp |c\esp c\esp c\esp c\esp  @{}|}
\hline
\multicolumn{3}{|c|}{Inst} &
\multicolumn{4}{|c|}{$(\ref{cotai_inf_ord}.1)$} &
\multicolumn{4}{|c|}{$(\ref{cotai_inf_ord}.2)$} &
\multicolumn{4}{|c|}{$(\ref{cotai_uij_max}.1)$}   \\
$|V|$ & $p$ & $\alpha$ &
$t(\#)$ & $t^*/gap^*$ & $\#nodes$ & $gap_{LR}$ &
$t(\#)$ & $t^*/gap^*$ & $\#nodes$ & $gap_{LR}$ &
$t(\#)$ & $t^*/gap^*$ & $\#nodes$ & $gap_{LR}$\\
\hline
100	&	4	&	0.4	&	0.7 & 0.8 & 178 & 13.97	&	0.7 & 1 & 153 & 55.44	&	\textbf{0.6} & 0.8 & 186 & 55.44	\\
100	&	4	&	0.6	&	2.4 & 16.7 & 6672 & 13.72	&	0.7 & 0.8 & 180 & 38.98	&	\textbf{0.6} & 0.7 & 152 & 38.98	\\
100	&	4	&	0.8	&	0.9 & 1.7 & 389 & 9.35	&	0.7 & 0.9 & 210 & 21.53	&	\textbf{0.6} & \textbf{0.7} & 302 & 21.53	\\
100	&	7	&	0.4	&	130.3 (8) & 11.47\% & 465836 & 12.34	&	1.8 & 2.4 & 275 & 52.18	&	\textbf{1} & 1.3 & 236 & 52.18	\\
100	&	7	&	0.6	&	180.8 (7) & 10.6\% & 942362 & 12.68	&	2 & 3.3 & 426 & 35.97	&	\textbf{1.1} & 1.4 & 480 & 35.97	\\
100	&	7	&	0.8	&	65.1 (9) & 5.98\% & 310921 & 9.18	&	2.6 & 3.8 & 1075 & 20.27	&	1.4 & \textbf{1.9} & \textbf{955} & 20.27	\\
100	&	10	&	0.4	&	421.9 (3) & 15.05\% & 1669042 & 13.35	&	4.1 & 5.5 & 381 & 50.66	&	1.5 & \textbf{1.9} & 299 & 50.66	\\
100	&	10	&	0.6	&	421.8 (3) & 12.88\% & 1692879 & 12.23	&	5.4 & 6.9 & 802 & 34.85	&	1.9 & 2.6 & 963 & 34.85	\\
100	&	10	&	0.8	&	309.1 (5) & 6.57\% & 994974 & 9.62	&	19.9 & 67.5 & 5947 & 20.2	&	6 & 19.4 & 6329 & 20.2	\\
225	&	4	&	0.4	&	2 & \textbf{2.6} & 1203 & 13.18	&	2.6 & 4.5 & 1099 & 55.09	&	2.1 & 4.4 & 1188 & 55.09	\\
225	&	4	&	0.6	&	1.9 & 2.5 & \textbf{721} & 13.07	&	2.3 & 3.7 & 1349 & 38.57	&	1.7 & 2.9 & 1236 & 38.57	\\
225	&	4	&	0.8	&	2.2 & 3.2 & \textbf{689} & 8.81	&	2.6 & 4.8 & 1098 & 21.09	&	\textbf{1.9} & 3.3 & 1101 & 21.09	\\
225	&	7	&	0.4	&	9.5 & 39.3 & 10622 & 11.31	&	18.8 & 54.7 & 6961 & 52.34	&	\textbf{6.7} & 19 & 7774 & 52.34	\\
225	&	7	&	0.6	&	\textbf{8.3} & \textbf{15.6} & 6182 & 11.88	&	25.7 & 70.2 & 10691 & 36.27	&	10.1 & 16.3 & \textbf{6038} & 36.27	\\
225	&	7	&	0.8	&	\textbf{14.7} & 55.4 & \textbf{8735} & 8.88	&	40.6 & 134.7 & 11264 & 20.32	&	17.1 & 62.5 & 10388 & 20.32	\\
225	&	10	&	0.4	&	72.3 (9) & 10.76\% & 54622 & 10.5	&	33 & 61.7 & 3162 & 50.25	&	7.6 & 13.5 & 2136 & 50.25	\\
225	&	10	&	0.6	&	93.8 (9) & 10.2\% & 38952 & 11.33	&	104 & 242.7 & 14785 & 34.56	&	32.2 & 124.2 & 15219 & 34.56	\\
225	&	10	&	0.8	&	\textbf{272.5} (\textbf{9}) & 8.13\% & 121097 & 8.93	&	569.6 (2) & 0.57\% & 51468 & 19.64	&	295.3 (8) & 0.32\% & 113983 & 19.62	\\
400	&	4	&	0.4	&	5.6 & 11.2 & \textbf{2224} & 12.99	&	10.6 & 33 & 3711 & 55.37	&	7.3 & 22.5 & 3345 & 55.37	\\
400	&	4	&	0.6	&	\textbf{6} & 10.2 & \textbf{2637} & 13.3	&	12 & 26.6 & 5024 & 39.04	&	6.8 & 12 & 4103 & 39.04	\\
400	&	4	&	0.8	&	9.7 & 17.2 & 4292 & 8.52	&	13.6 & 33.5 & 6133 & 21.03	&	9.1 & 22.2 & 5397 & 21.03	\\
400	&	7	&	0.4	&	28.6 & 104.3 & 13958 & 11.29	&	102.4 & 413 & 11179 & 52.05	&	34.4 & 144.7 & 10489 & 52.05	\\
400	&	7	&	0.6	&	\textbf{54.1} & \textbf{147} & 30203 & 11.98	&	236 & 565.2 & 60207 & 36.12	&	75.4 & 182.8 & \textbf{25950} & 36.12	\\
400	&	7	&	0.8	&	\textbf{68.2} & \textbf{133.8} & \textbf{28042} & 8.88	&	202.3 & 395.8 & 28993 & 20.19	&	85 & 188.2 & 28762 & 20.19	\\
400	&	10	&	0.4	&	107.4 (9) & 10.68\% & 33049 & 10.27	&	222.5 (9) & 0.38\% & \textbf{9768} & 50.59	&	68.6 & 198.8 & 13777 & 50.58	\\
400	&	10	&	0.6	&	\textbf{205.3} & \textbf{540.3} & 87754 & 10.73	&	523.4 (2) & 0.35\% & 21252 & 34.55	&	290.8 (9) & 0.11\% & 61637 & 34.54	\\
400	&	10	&	0.8	&	586.8 (\textbf{1}) & 0.36\% & 154947 & 8.51	&	599.7 (0) & 0.73\% & 38344 & 19.54	&	583.5 (\textbf{1}) & 0.42\% & 96613 & 19.5	\\

\hline
\end{tabular}
\end{scriptsize}
\end{center}
\caption{Results obtained for the OWAP formulations with the Perfect Matching Problem and valid inequalities}
%\label{tab_exact_resolution}
\end{table}
%-----------------------------------------------------------------------------------------------------------------------
%=======================================================================================================================
\newpage
%=======================================================================================================================
%\subsubsection*{PMP con DV:  $F_{R1}^{z}$ + (29.2), (30.1), (30.2), (31), (33.1), (33.2)}
%-----------------------------------------------------------------------------------------------------------------------
\begin{table}[h!]
\renewcommand{\esp}{@{\hspace{0.15cm}}}
\begin{center}
\begin{scriptsize}
\begin{tabular}{|@{}\esp c \esp c \esp c\esp |c\esp c\esp c\esp c\esp |c\esp c\esp c\esp c\esp |c\esp c\esp c\esp c\esp  @{}|}
\hline
\multicolumn{3}{|c|}{Inst} &
\multicolumn{4}{|c|}{$(\ref{cotai_uij_max}.2)$} &
\multicolumn{4}{|c|}{$(\ref{cotaj_uij_max}.1)$} &
\multicolumn{4}{|c|}{$(\ref{cotaj_uij_max}.2)$}   \\
$|V|$ & $p$ & $\alpha$ &
$t(\#)$ & $t^*/gap^*$ & $\#nodes$ & $gap_{LR}$ &
$t(\#)$ & $t^*/gap^*$ & $\#nodes$ & $gap_{LR}$ &
$t(\#)$ & $t^*/gap^*$ & $\#nodes$ & $gap_{LR}$\\
\hline
100	&	4	&	0.4	&	\textbf{0.6} & \textbf{0.7} & 186 & 55.44	&	\textbf{0.6} & \textbf{0.7} & \textbf{106} & \textbf{13.23}	&	0.7 & 1 & 166 & 55.44	\\
100	&	4	&	0.6	&	\textbf{0.6} & 0.7 & 152 & 38.98	&	0.7 & 0.8 & \textbf{129} & \textbf{13.27}	&	0.7 & 0.9 & 156 & 38.98	\\
100	&	4	&	0.8	&	\textbf{0.6} & \textbf{0.7} & 302 & 21.53	&	0.7 & 0.8 & \textbf{144} & 9.13	&	0.8 & 0.9 & 269 & 21.53	\\
100	&	7	&	0.4	&	\textbf{1} & 1.3 & 236 & 52.18	&	\textbf{1} & 1.3 & \textbf{125} & \textbf{11.75}	&	2 & 3.3 & 392 & 52.18	\\
100	&	7	&	0.6	&	\textbf{1.1} & 1.4 & 480 & 35.97	&	1.2 & 1.7 & 438 & \textbf{11.91}	&	2 & 2.8 & 470 & 35.97	\\
100	&	7	&	0.8	&	1.4 & \textbf{1.9} & \textbf{955} & 20.27	&	\textbf{1.3} & 2 & 984 & 9.03	&	3.1 & 4.8 & 1230 & 20.27	\\
100	&	10	&	0.4	&	1.5 & \textbf{1.9} & 299 & 50.66	&	1.6 & 2 & \textbf{267} & \textbf{10.97}	&	4.4 & 5.2 & 535 & 50.66	\\
100	&	10	&	0.6	&	1.9 & 2.6 & 963 & 34.85	&	\textbf{1.8} & 2.5 & \textbf{608} & \textbf{11.56}	&	6.4 & 8.2 & 1222 & 34.85	\\
100	&	10	&	0.8	&	6 & 19.5 & 6329 & 20.2	&	\textbf{5.5} & 22.2 & 7908 & \textbf{9.5}	&	23.6 & 87 & 6403 & 20.2	\\
225	&	4	&	0.4	&	2.1 & 4.4 & 1188 & 55.09	&	\textbf{1.7} & 3.4 & 1116 & \textbf{12.67}	&	2.9 & 6.6 & 1324 & 55.09	\\
225	&	4	&	0.6	&	1.7 & 2.9 & 1236 & 38.57	&	\textbf{1.6} & \textbf{2.2} & 1038 & \textbf{12.76}	&	2.2 & 3.4 & 1025 & 38.57	\\
225	&	4	&	0.8	&	\textbf{1.9} & 3.2 & 1101 & 21.09	&	2.1 & \textbf{3} & 995 & 8.66	&	2.3 & 3.9 & 1083 & 21.09	\\
225	&	7	&	0.4	&	\textbf{6.7} & 18.9 & 7774 & 52.34	&	9.5 & 38.8 & 6128 & \textbf{10.97}	&	24.1 & 83.4 & 7365 & 52.34	\\
225	&	7	&	0.6	&	10.1 & 16.2 & \textbf{6038} & 36.27	&	11.6 & 30 & 7394 & \textbf{11.67}	&	36.1 & 81.2 & 9748 & 36.27	\\
225	&	7	&	0.8	&	17.2 & 62.7 & 10388 & 20.32	&	19.8 & 84.1 & 11054 & 8.78	&	49.9 & 179.9 & 11402 & 20.32	\\
225	&	10	&	0.4	&	7.5 & 13.4 & 2136 & 50.25	&	7 & 12.9 & 1527 & \textbf{9.94}	&	31.8 & 61.8 & 2363 & 50.25	\\
225	&	10	&	0.6	&	32.2 & 123.9 & 15219 & 34.56	&	36.5 & 99.5 & 15362 & \textbf{11}	&	142 & 375.1 & 14273 & 34.56	\\
225	&	10	&	0.8	&	295.2 (8) & 0.32\% & 114087 & 19.62	&	356.5 (7) & 0.32\% & 148693 & \textbf{8.78}	&	594.5 (1) & 0.69\% & \textbf{47129} & 19.66	\\
400	&	4	&	0.4	&	7.3 & 22.5 & 3345 & 55.37	&	\textbf{4.9} & \textbf{9.7} & 2506 & \textbf{12.74}	&	12.2 & 34.2 & 4227 & 55.37	\\
400	&	4	&	0.6	&	6.8 & 11.9 & 4103 & 39.04	&	\textbf{6} & \textbf{10.1} & 2800 & \textbf{13.15}	&	11.8 & 25.2 & 7289 & 39.04	\\
400	&	4	&	0.8	&	9.1 & 22.3 & 5397 & 21.03	&	\textbf{8.7} & \textbf{15} & 4689 & 8.45	&	14.2 & 31.3 & 7456 & 21.03	\\
400	&	7	&	0.4	&	34.5 & 144.6 & 10489 & 52.05	&	\textbf{26.2} & \textbf{70} & 7118 & \textbf{11.04}	&	63.5 & 189 & \textbf{6801} & 52.05	\\
400	&	7	&	0.6	&	75.3 & 182 & \textbf{25950} & 36.12	&	89.4 & 255.1 & 29574 & \textbf{11.84}	&	294.7 (8) & 0.22\% & 38975 & 36.12	\\
400	&	7	&	0.8	&	84.8 & 187.7 & 28762 & 20.19	&	86 & 139.6 & 35249 & 8.82	&	224.6 & 449.8 & 29823 & 20.19	\\
400	&	10	&	0.4	&	68.6 & 198.5 & 13777 & 50.58	&	58 & 178.9 & 10826 & \textbf{9.98}	&	312.7 (8) & 0.48\% & 12545 & 50.59	\\
400	&	10	&	0.6	&	290.4 (9) & 0.11\% & 61365 & 34.54	&	299.2 (9) & 0.08\% & 60221 & \textbf{10.63}	&	517.2 (3) & 0.45\% & \textbf{18359} & 34.55	\\
400	&	10	&	0.8	&	583.5 (\textbf{1}) & 0.42\% & 96664 & 19.5	&	598.3 (\textbf{1}) & 0.34\% & 93994 & \textbf{8.47}	&	599.9 (0) & 0.59\% & \textbf{20318} & 19.57	\\

\hline
\end{tabular}
\end{scriptsize}
\end{center}
\caption{Results obtained for the OWAP formulations with the Perfect Matching Problem and valid inequalities}
%\label{tab_exact_resolution}
\end{table}
%-----------------------------------------------------------------------------------------------------------------------

%-----------------------------------------------------------------------------------------------------------------------
\begin{table}[h!]
\renewcommand{\esp}{@{\hspace{0.15cm}}}
\begin{center}
\begin{scriptsize}
\begin{tabular}{|@{}\esp c \esp c \esp c\esp |c\esp c\esp c\esp c\esp |c\esp c\esp c\esp c\esp |c\esp c\esp c\esp c\esp  @{}|}
\hline
\multicolumn{3}{|c|}{Inst} &
\multicolumn{4}{|c|}{$(\ref{cotazy})$} &
\multicolumn{4}{|c|}{$(\ref{validsubsets}.1)$} &
\multicolumn{4}{|c|}{$(\ref{validsubsets}.2)$}   \\
$|V|$ & $p$ & $\alpha$ &
$t(\#)$ & $t^*/gap^*$ & $\#nodes$ & $gap_{LR}$ &
$t(\#)$ & $t^*/gap^*$ & $\#nodes$ & $gap_{LR}$ &
$t(\#)$ & $t^*/gap^*$ & $\#nodes$ & $gap_{LR}$\\
\hline
100	&	4	&	0.4	&	1 & 1.4 & 357 & 55.44	&	\textbf{0.6} & 0.8 & 140 & 55.44	&	0.7 & 1 & 186 & 55.44	\\
100	&	4	&	0.6	&	0.9 & 1.3 & 406 & 38.98	&	\textbf{0.6} & \textbf{0.6} & 159 & 38.98	&	\textbf{0.6} & 0.7 & 152 & 38.98	\\
100	&	4	&	0.8	&	1 & 1.4 & 372 & 21.53	&	\textbf{0.6} & \textbf{0.7} & 242 & 21.53	&	\textbf{0.6} & \textbf{0.7} & 302 & 21.53	\\
100	&	7	&	0.4	&	3 & 5.8 & 3307 & 52.18	&	\textbf{1} & 1.4 & 214 & 52.18	&	\textbf{1} & 1.3 & 236 & 52.18	\\
100	&	7	&	0.6	&	4.3 & 7.8 & 5187 & 35.97	&	\textbf{1.1} & 1.4 & \textbf{402} & 35.97	&	\textbf{1.1} & 1.4 & 480 & 35.97	\\
100	&	7	&	0.8	&	6.5 & 14.1 & 9016 & 20.27	&	1.5 & 2.2 & 1165 & 20.27	&	1.4 & 2 & \textbf{955} & 20.27	\\
100	&	10	&	0.4	&	131.5 & 493.7 & 56509 & 50.66	&	1.6 & 2 & 363 & 50.66	&	1.6 & 2 & 299 & 50.66	\\
100	&	10	&	0.6	&	348.2 (8) & 5.36\% & 145975 & 35.16	&	1.9 & 2.9 & 757 & 34.85	&	1.9 & 2.6 & 963 & 34.85	\\
100	&	10	&	0.8	&	573.4 (2) & 7.81\% & 148813 & 22.74	&	6.1 & 23.6 & 6411 & 20.2	&	6 & 19.5 & 6329 & 20.2	\\
225	&	4	&	0.4	&	3.1 & 7.3 & 1760 & 55.09	&	2 & 3.7 & \textbf{897} & 55.09	&	2.1 & 4.4 & 1188 & 55.09	\\
225	&	4	&	0.6	&	3.3 & 5.6 & 2229 & 38.57	&	1.8 & 2.7 & 950 & 38.57	&	1.7 & 2.9 & 1151 & 38.57	\\
225	&	4	&	0.8	&	3.6 & 8.3 & 2023 & 21.09	&	\textbf{1.9} & 3.3 & 1142 & 21.09	&	\textbf{1.9} & 3.1 & 1101 & 21.09	\\
225	&	7	&	0.4	&	91.2 (9) & 0.57\% & 53230 & 52.36	&	9.6 & 39.8 & 6398 & 52.34	&	6.9 & 20.3 & 8260 & 52.34	\\
225	&	7	&	0.6	&	163.5 & 432.3 & 92543 & 36.27	&	10.6 & 20.8 & 6258 & 36.27	&	10.1 & 16.2 & \textbf{6038} & 36.27	\\
225	&	7	&	0.8	&	155.9 & 479.1 & 68304 & 20.32	&	17.1 & 52.1 & 9917 & 20.32	&	17.3 & 62.7 & 10635 & 20.32	\\
225	&	10	&	0.4	&	599.9 (0) & 21.12\% & 58860 & 55.15	&	7.8 & 13.2 & 2118 & 50.25	&	7.6 & 13.6 & 2136 & 50.25	\\
225	&	10	&	0.6	&	599.9 (0) & 17.98\% & 67172 & 38.68	&	31.9 & 78.6 & \textbf{12701} & 34.56	&	32.7 & 124.5 & 15743 & 34.56	\\
225	&	10	&	0.8	&	599.8 (0) & 11.01\% & 53762 & 23.51	&	349.6 (\textbf{9}) & \textbf{0.13}\% & 138048 & 19.62	&	295.8 (\textbf{9}) & 0.32\% & 116525 & 19.62	\\
400	&	4	&	0.4	&	9 & 22.2 & 6168 & 55.37	&	6.6 & 14 & 3260 & 55.37	&	7.3 & 22.5 & 3345 & 55.37	\\
400	&	4	&	0.6	&	10.4 & 29.5 & 5976 & 39.04	&	10.5 & 30.5 & 6922 & 39.04	&	6.8 & 11.9 & 4172 & 39.04	\\
400	&	4	&	0.8	&	13.8 & 30.6 & 9083 & 21.03	&	11.6 & 31.7 & 10499 & 21.03	&	8.8 & 22.4 & 5203 & 21.03	\\
400	&	7	&	0.4	&	181.9 & 587.3 & 46354 & 52.05	&	49.3 & 284.4 & 14528 & 52.05	&	35.4 & 145 & 10976 & 52.05	\\
400	&	7	&	0.6	&	499.6 (4) & 1.07\% & 113087 & 36.31	&	82 & 258.9 & 28264 & 36.12	&	83.7 & 251.1 & 27775 & 36.12	\\
400	&	7	&	0.8	&	401.9 (6) & 0.54\% & 90350 & 20.25	&	79.5 & 142.5 & 32783 & 20.19	&	84.9 & 188.3 & 28762 & 20.19	\\
400	&	10	&	0.4	&	600 (0) & 19.09\% & 59670 & 53.93	&	86.8 & 342.2 & 17939 & 50.58	&	68.4 & 196.2 & 13777 & 50.58	\\
400	&	10	&	0.6	&	600 (0) & 12.78\% & 81380 & 39.41	&	279.9 (9) & 0.14\% & 55631 & 34.54	&	289.9 (9) & 0.11\% & 62061 & 34.54	\\
400	&	10	&	0.8	&	600 (0) & 9.53\% & 77123 & 23.49	&	596.1 (\textbf{1}) & \textbf{0.29}\% & 93226 & 19.49	&	583.5 (\textbf{1}) & 0.42\% & 96823 & 19.5	\\

\hline
\end{tabular}
\end{scriptsize}
\end{center}
\caption{Results obtained for the OWAP formulations with the Perfect Matching Problem and valid inequalities}
%\label{tab_exact_resolution}
\end{table}
%-----------------------------------------------------------------------------------------------------------------------
%=======================================================================================================================
\newpage
%=======================================================================================================================
%\subsubsection*{PMP con DV:  $F_{R1}^{z}$ + (33.3), (33.4), (34), (35), (36), (37)}
%-----------------------------------------------------------------------------------------------------------------------
\begin{table}[h!]
\renewcommand{\esp}{@{\hspace{0.15cm}}}
\begin{center}
\begin{scriptsize}
\begin{tabular}{|@{}\esp c \esp c \esp c\esp |c\esp c\esp c\esp c\esp |c\esp c\esp c\esp c\esp |c\esp c\esp c\esp c\esp  @{}|}
\hline
\multicolumn{3}{|c|}{Inst} &
\multicolumn{4}{|c|}{$(\ref{validsubsets}.3)$} &
\multicolumn{4}{|c|}{$(\ref{validsubsets}.4)$} &
\multicolumn{4}{|c|}{$(\ref{vi:owa2eq})$}   \\
$|V|$ & $p$ & $\alpha$ &
$t(\#)$ & $t^*/gap^*$ & $\#nodes$ & $gap_{LR}$ &
$t(\#)$ & $t^*/gap^*$ & $\#nodes$ & $gap_{LR}$ &
$t(\#)$ & $t^*/gap^*$ & $\#nodes$ & $gap_{LR}$\\
\hline
100	&	4	&	0.4	&	\textbf{0.6} & \textbf{0.7} & 186 & 55.44	&	\textbf{0.6} & \textbf{0.7} & 186 & 55.44	&	1.6 & 1.9 & 262 & 41.27	\\
100	&	4	&	0.6	&	\textbf{0.6} & 0.7 & 152 & 38.98	&	\textbf{0.6} & 0.7 & 152 & 38.98	&	1.3 & 1.8 & 353 & 19.56	\\
100	&	4	&	0.8	&	\textbf{0.6} & \textbf{0.7} & 302 & 21.53	&	\textbf{0.6} & \textbf{0.7} & 302 & 21.53	&	0.7 & 0.9 & 279 & \textbf{2.1}	\\
100	&	7	&	0.4	&	\textbf{1} & 1.3 & 236 & 52.18	&	\textbf{1} & 1.3 & 236 & 52.18	&	5.1 & 9.5 & 3106 & 45.38	\\
100	&	7	&	0.6	&	\textbf{1.1} & \textbf{1.3} & 480 & 35.97	&	\textbf{1.1} & 1.4 & 480 & 35.97	&	7.5 & 14.2 & 5210 & 26.87	\\
100	&	7	&	0.8	&	1.4 & \textbf{1.9} & \textbf{955} & 20.27	&	1.4 & 2 & \textbf{955} & 20.27	&	10.1 & 15.4 & 10804 & \textbf{8.93}	\\
100	&	10	&	0.4	&	1.5 & \textbf{1.9} & 299 & 50.66	&	1.5 & \textbf{1.9} & 299 & 50.66	&	35.7 & 109.5 & 36224 & 46.24	\\
100	&	10	&	0.6	&	1.9 & 2.5 & 963 & 34.85	&	1.9 & 2.6 & 963 & 34.85	&	99.1 & 292.7 & 88305 & 29.02	\\
100	&	10	&	0.8	&	6 & 19.5 & 6329 & 20.2	&	6 & 19.5 & 6329 & 20.2	&	341.4 (7) & 10.43\% & 225313 & 14.95	\\
225	&	4	&	0.4	&	2.1 & 4.4 & 1188 & 55.09	&	2.1 & 4.4 & 1188 & 55.09	&	5.1 & 11.6 & 2247 & 40.6	\\
225	&	4	&	0.6	&	1.8 & 3 & 1236 & 38.57	&	1.8 & 3 & 1236 & 38.57	&	4.8 & 9 & 1998 & 18.74	\\
225	&	4	&	0.8	&	\textbf{1.9} & 3.1 & 1101 & 21.09	&	\textbf{1.9} & 3.2 & 1101 & 21.09	&	2.3 & 5.4 & 781 & \textbf{1.44}	\\
225	&	7	&	0.4	&	7.1 & 22.4 & 9432 & 52.34	&	\textbf{6.7} & \textbf{18.8} & 7785 & 52.34	&	169.5 (9) & 0.28\% & 72554 & 44.84	\\
225	&	7	&	0.6	&	10.1 & 16.2 & \textbf{6038} & 36.27	&	10 & 16.2 & \textbf{6038} & 36.27	&	228.4 & 487.4 & 78835 & 26.22	\\
225	&	7	&	0.8	&	17.1 & 62.5 & 10499 & 20.32	&	17.2 & 62.5 & 10533 & 20.32	&	166.5 (9) & 0.6\% & 60230 & \textbf{7.79}	\\
225	&	10	&	0.4	&	7.5 & 13.3 & 2136 & 50.25	&	7.5 & 13.3 & 2136 & 50.25	&	567 (2) & 11.81\% & 118429 & 46.76	\\
225	&	10	&	0.6	&	32.2 & 124.2 & 15322 & 34.56	&	32.1 & 123.8 & 15124 & 34.56	&	596.9 (1) & 9.23\% & 128499 & 30.09	\\
225	&	10	&	0.8	&	291.1 (\textbf{9}) & 0.32\% & 114568 & 19.62	&	299.2 (8) & 0.32\% & 116112 & 19.62	&	599.4 (0) & 5.82\% & 130289 & 14.6	\\
400	&	4	&	0.4	&	7.3 & 22.5 & 3339 & 55.37	&	7.3 & 22.5 & 3345 & 55.37	&	15.9 & 55.5 & 4037 & 40.66	\\
400	&	4	&	0.6	&	6.7 & 11.9 & 4057 & 39.04	&	6.8 & 11.9 & 4103 & 39.04	&	16 & 31.2 & 5902 & 18.94	\\
400	&	4	&	0.8	&	9.1 & 22.3 & 5402 & 21.03	&	9 & 22.2 & 5397 & 21.03	&	10.2 & 24.2 & \textbf{3531} & \textbf{1.29}	\\
400	&	7	&	0.4	&	34.6 & 145.1 & 10594 & 52.05	&	34.6 & 144.3 & 10638 & 52.05	&	344.5 (8) & 0.7\% & 74878 & 44.34	\\
400	&	7	&	0.6	&	75.3 & 182 & \textbf{25950} & 36.12	&	83.6 & 252.2 & 27604 & 36.12	&	508.8 (3) & 1.14\% & 88507 & 26.06	\\
400	&	7	&	0.8	&	84.7 & 188.1 & 28762 & 20.19	&	84.7 & 188.6 & 28762 & 20.19	&	499.6 (5) & 0.27\% & 121674 & \textbf{7.35}	\\
400	&	10	&	0.4	&	68.5 & 198.2 & 13777 & 50.58	&	68.5 & 197.7 & 13777 & 50.58	&	599.8 (0) & 18.51\% & 65891 & 48.51	\\
400	&	10	&	0.6	&	289.6 (9) & 0.11\% & 61601 & 34.54	&	289.5 (9) & 0.11\% & 61833 & 34.54	&	599.8 (0) & 9.38\% & 70265 & 30.73	\\
400	&	10	&	0.8	&	583.6 (\textbf{1}) & 0.42\% & 97099 & 19.5	&	\textbf{583.4} (\textbf{1}) & 0.42\% & 96958 & 19.5	&	599.9 (0) & 7.4\% & 96536 & 14.49	\\

\hline
\end{tabular}
\end{scriptsize}
\end{center}
\caption{Results obtained for the OWAP formulations with the Perfect Matching Problem and valid inequalities}
%\label{tab_exact_resolution}
\end{table}
%-----------------------------------------------------------------------------------------------------------------------

%-----------------------------------------------------------------------------------------------------------------------
\begin{table}[h!]
\renewcommand{\esp}{@{\hspace{0.15cm}}}
\begin{center}
\begin{scriptsize}
\begin{tabular}{|@{}\esp c \esp c \esp c\esp |c\esp c\esp c\esp c\esp |c\esp c\esp c\esp c\esp |c\esp c\esp c\esp c\esp  @{}|}
\hline
\multicolumn{3}{|c|}{Inst} &
\multicolumn{4}{|c|}{$(\ref{vi:cotayydis1})$} &
\multicolumn{4}{|c|}{$(\ref{vi:cotayydis2})$} &
\multicolumn{4}{|c|}{$(\ref{vi:yyrel1})$}   \\
$|V|$ & $p$ & $\alpha$ &
$t(\#)$ & $t^*/gap^*$ & $\#nodes$ & $gap_{LR}$ &
$t(\#)$ & $t^*/gap^*$ & $\#nodes$ & $gap_{LR}$ &
$t(\#)$ & $t^*/gap^*$ & $\#nodes$ & $gap_{LR}$\\
\hline
100	&	4	&	0.4	&	\textbf{0.6} & 0.8 & 139 & 55.44	&	1.2 & 1.3 & 159 & 55.44	&	1.2 & 1.3 & 236 & 55.44	\\
100	&	4	&	0.6	&	\textbf{0.6} & 0.7 & 140 & 38.98	&	1.2 & 1.3 & 156 & 38.98	&	1.2 & 1.5 & 283 & 38.98	\\
100	&	4	&	0.8	&	\textbf{0.6} & 0.8 & 329 & 21.53	&	1.1 & 1.3 & 269 & 21.53	&	1.2 & 1.3 & 239 & 21.53	\\
100	&	7	&	0.4	&	\textbf{1} & \textbf{1.2} & 256 & 52.18	&	1.9 & 2.2 & 247 & 52.18	&	4 & 4.9 & 779 & 52.18	\\
100	&	7	&	0.6	&	1.2 & 2 & 591 & 35.97	&	2.2 & 2.7 & 516 & 35.97	&	4.6 & 5.6 & 1298 & 35.97	\\
100	&	7	&	0.8	&	1.6 & 3 & 1008 & 20.27	&	2.7 & 3.3 & 1105 & 20.27	&	4.4 & 6.2 & 1643 & 20.27	\\
100	&	10	&	0.4	&	1.7 & 4.1 & 580 & 50.66	&	2.7 & 3 & 302 & 50.66	&	9.3 & 12.3 & 969 & 50.66	\\
100	&	10	&	0.6	&	1.9 & 2.9 & 985 & 34.85	&	3.5 & 4.3 & 982 & 34.85	&	14.9 & 44.1 & 2966 & 34.85	\\
100	&	10	&	0.8	&	5.6 & 17.7 & \textbf{5364} & 20.2	&	10.2 & 31.2 & 6769 & 20.2	&	48.2 & 177.8 & 22384 & 20.2	\\
225	&	4	&	0.4	&	2 & 2.9 & 990 & 55.09	&	3.9 & 6.7 & 998 & 55.09	&	3.8 & 5.8 & 1297 & 55.09	\\
225	&	4	&	0.6	&	1.7 & 2.5 & 1239 & 38.57	&	3.9 & 5.7 & 1179 & 38.57	&	4 & 5.5 & 1710 & 38.57	\\
225	&	4	&	0.8	&	\textbf{1.9} & 3.7 & 1240 & 21.09	&	3.5 & 5.8 & 961 & 21.09	&	4 & 6.8 & 1912 & 21.09	\\
225	&	7	&	0.4	&	8.3 & 35.8 & \textbf{5617} & 52.34	&	14.9 & 59.5 & 9584 & 52.34	&	29.4 & 69.5 & 8433 & 52.34	\\
225	&	7	&	0.6	&	9.7 & 18.1 & 6206 & 36.27	&	16.6 & 32 & 6614 & 36.27	&	37.8 & 97.9 & 12295 & 36.27	\\
225	&	7	&	0.8	&	17.2 & 48.3 & 10746 & 20.32	&	23.1 & 81.7 & 10207 & 20.32	&	39.4 & 161.5 & 15449 & 20.32	\\
225	&	10	&	0.4	&	7.4 & \textbf{12.4} & 2464 & 50.25	&	10.7 & 16.3 & \textbf{1516} & 50.25	&	42.3 & 77.3 & 5457 & 50.25	\\
225	&	10	&	0.6	&	33.8 & 89.7 & 13763 & 34.56	&	51.5 & 105.8 & 15061 & 34.56	&	129.3 & 340.5 & 20824 & 34.56	\\
225	&	10	&	0.8	&	338.5 (7) & 12.07\% & 130287 & 19.7	&	416.8 (7) & 0.38\% & 116869 & 19.63	&	545.3 (2) & 0.55\% & 94492 & 19.64	\\
400	&	4	&	0.4	&	6.2 & 15.3 & 2546 & 55.37	&	10.5 & 21.2 & 2513 & 55.37	&	12 & 16.4 & 4533 & 55.37	\\
400	&	4	&	0.6	&	7.5 & 17 & 4044 & 39.04	&	15 & 43.4 & 5642 & 39.04	&	17 & 28 & 6513 & 39.04	\\
400	&	4	&	0.8	&	11.4 & 44.8 & 6263 & 21.03	&	16.4 & 33.4 & 5405 & 21.03	&	19.5 & 35.9 & 6255 & 21.03	\\
400	&	7	&	0.4	&	48.8 & 256 & 15696 & 52.05	&	52.4 & 219.6 & 10383 & 52.05	&	120.6 (9) & 0.4\% & 18046 & 52.05	\\
400	&	7	&	0.6	&	74.4 & 209.9 & 26944 & 36.12	&	119.5 & 359.3 & 27775 & 36.12	&	201.5 & 485.5 & 45604 & 36.12	\\
400	&	7	&	0.8	&	98 & 182.7 & 35369 & 20.19	&	142.6 & 307.4 & 34805 & 20.19	&	249.1 & 440.3 & 53146 & 20.19	\\
400	&	10	&	0.4	&	86.4 & 384.5 & 17514 & 50.58	&	121.3 & 402.7 & 18138 & 50.58	&	325.4 (9) & 0.26\% & 24495 & 50.58	\\
400	&	10	&	0.6	&	334.7 (9) & 0.24\% & 69498 & 34.54	&	434.8 (5) & 0.24\% & 57897 & 34.54	&	511 (3) & 0.53\% & 63696 & 34.58	\\
400	&	10	&	0.8	&	598.7 (\textbf{1}) & 0.4\% & 93441 & 19.5	&	599.8 (0) & 0.46\% & 65597 & 19.51	&	599.6 (0) & 13.85\% & 36099 & 19.6	\\

\hline
\end{tabular}
\end{scriptsize}
\end{center}
\caption{Results obtained for the OWAP formulations with the Perfect Matching Problem and valid inequalities}
%\label{tab_exact_resolution}
\end{table}
%-----------------------------------------------------------------------------------------------------------------------
%=======================================================================================================================
\newpage
%=======================================================================================================================
%\subsubsection*{PMP con DV:  $F_{R1}^{z}$ + (38), (39)}
%-----------------------------------------------------------------------------------------------------------------------
\begin{table}[h!]
\renewcommand{\esp}{@{\hspace{0.15cm}}}
\begin{center}
\begin{scriptsize}
\begin{tabular}{|@{}\esp c \esp c \esp c\esp |c\esp c\esp c\esp c\esp |c\esp c\esp c\esp c\esp  @{}|}
\hline
\multicolumn{3}{|c|}{Inst} &
\multicolumn{4}{|c|}{$(\ref{vi:yyrel2})$} &
\multicolumn{4}{|c|}{$(\ref{vi:yyrel3})$} \\
$|V|$ & $p$ & $\alpha$ &
$t(\#)$ & $t^*/gap^*$ & $\#nodes$ & $gap_{LR}$ &
$t(\#)$ & $t^*/gap^*$ & $\#nodes$ & $gap_{LR}$\\
\hline
100	&	4	&	0.4	&	1.2 & 1.4 & 230 & 55.44	&	0.9 & 1.3 & 256 & 55.44\\
100	&	4	&	0.6	&	1.1 & 1.4 & 268 & 38.98	&	0.9 & 1.2 & 271 & 38.98\\
100	&	4	&	0.8	&	1 & 1.2 & 210 & 21.53	&	1 & 1.3 & 371 & 21.53\\
100	&	7	&	0.4	&	3.3 & 3.8 & 579 & 52.18	&	3 & 3.9 & 683 & 52.18\\
100	&	7	&	0.6	&	3.5 & 4.1 & 684 & 35.97	&	3.5 & 5.5 & 1014 & 35.97\\
100	&	7	&	0.8	&	4.1 & 5.5 & 1089 & 20.27	&	4.3 & 7.5 & 1544 & 20.27\\
100	&	10	&	0.4	&	8.6 & 27.8 & 1484 & 50.66	&	6.3 & 7.4 & 511 & 50.66\\
100	&	10	&	0.6	&	8.2 & 10.3 & 1389 & 34.85	&	9 & 12.4 & 1513 & 34.85\\
100	&	10	&	0.8	&	132.1 (8) & 11.79\% & 88337 & 20.29	&	27.8 & 92.7 & 6895 & 20.2\\
225	&	4	&	0.4	&	3.7 & 5.6 & 1335 & 55.09	&	3.4 & 4.5 & 1690 & 55.09\\
225	&	4	&	0.6	&	3.6 & 4.1 & 1185 & 38.57	&	3 & 5 & 1715 & 38.57\\
225	&	4	&	0.8	&	3.3 & 5.1 & 993 & 21.09	&	3.7 & 5.2 & 1925 & 21.09\\
225	&	7	&	0.4	&	24.8 & 94.3 & 7360 & 52.34	&	30.6 & 126.9 & 8058 & 52.34\\
225	&	7	&	0.6	&	32.8 & 98.1 & 13897 & 36.27	&	38.9 & 110.2 & 13943 & 36.27\\
225	&	7	&	0.8	&	32.8 & 113.5 & 12765 & 20.32	&	55.7 & 226.4 & 13267 & 20.32\\
225	&	10	&	0.4	&	33.1 & 64 & 5220 & 50.25	&	30.9 & 54.7 & 2591 & 50.25\\
225	&	10	&	0.6	&	98.2 & 302.5 & 20554 & 34.56	&	156.8 & 352.3 & 15493 & 34.56\\
225	&	10	&	0.8	&	520.4 (4) & 0.43\% & 99710 & 19.63	&	592.1 (1) & 0.65\% & 49121 & 19.65\\
400	&	4	&	0.4	&	13.2 & 25.5 & 3826 & 55.37	&	13.4 & 33.4 & 4517 & 55.37\\
400	&	4	&	0.6	&	18.4 & 39 & 7049 & 39.04	&	15.7 & 32.4 & 5985 & 39.04\\
400	&	4	&	0.8	&	17.9 & 54.4 & 5108 & 21.03	&	76.8 (9) & 13.06\% & 31435 & 21.03\\
400	&	7	&	0.4	&	116.5 (9) & 0.37\% & 13888 & 52.05	&	139.1 (9) & 0.1\% & 14289 & 52.05\\
400	&	7	&	0.6	&	207.1 & 504 & 34054 & 36.12	&	263 & 488.1 & 30940 & 36.12\\
400	&	7	&	0.8	&	249.1 & 530.2 & 48541 & 20.19	&	334.5 (8) & 0.07\% & 42947 & 20.19\\
400	&	10	&	0.4	&	191.8 (9) & 0.17\% & 22990 & 50.58	&	269.7 (9) & 0.32\% & 12771 & 50.58\\
400	&	10	&	0.6	&	515 (3) & 0.33\% & 53833 & 34.55	&	527.1 (3) & 0.45\% & 35374 & 34.57\\
400	&	10	&	0.8	&	599.7 (0) & 0.48\% & 42979 & 19.55	&	599.8 (0) & 0.73\% & 37615 & 19.58\\

\hline
\end{tabular}
\end{scriptsize}
\end{center}
\caption{Results obtained for the OWAP formulations with the Perfect Matching Problem and valid inequalities}
%\label{tab_exact_resolution}
\end{table}
%-----------------------------------------------------------------------------------------------------------------------

%=======================================================================================================================

%=======================================================================================================================
\section{Conclusions}\label{Sec:Conclusions}
In this work we have presented and revisited different mathematical programming formulations for the OWAP using different sets of decision variables. These formulations reinforced with appropriate constraints have shown to be rather promising for efficiently solving many medium size OWAP instances. However, from the obtained results it is also clear that for solving larger OWAP instances with more objective functions further improvements are needed. Our current research focuses on the design of more sophisticated elimination tests as well as from alternative formulations leading to tighter LP bounds.
%=======================================================================================================================

%=======================================================================================================================
\vspace{1.5cm}
\noindent{\large\textbf{Acknowledgements}}\newline
This research has been partially supported by the Spanish Ministry of Science and Education through grants MTM2009-14039-C06-05 and MTM2010-19576-C02-01, by Junta de Andalucía grant FQM5849 and by ERDF funds. This support is gratefully acknowledged.

%=======================================================================================================================
%=======================================================================================================================
\setlength{\parskip}{0pt}
\bibliographystyle{estilobibs1}
\begin{small}

\end{small}
%=======================================================================================================================

\end{document}